\documentclass[a4paper,fleqn]{cas-sc}  

\usepackage[square, comma, sort&compress, numbers]{natbib}

\def\tsc#1{\csdef{#1}{\textsc{\lowercase{#1}}\xspace}}
\tsc{WGM}
\tsc{QE}

\begin{document}
	\let\WriteBookmarks\relax
	\def\floatpagepagefraction{1}
	\def\textpagefraction{.001}
	
\shorttitle{A meshless GFDM for solving the Stokes-Darcy coupled problem in static and moving systems}    

\shortauthors{Y.N.Xing}  

\title [mode = title]{A meshless generalized finite difference method for solving the Stokes-Darcy coupled problem in static and moving systems}  
	
	
	%
	
	\author[a]{Yanan Xing}
	
	
	
	
	
	\credit{Methodology, Writing the initial draft and the final draft}
	
	\affiliation[a]{
		addressline={School of Mathematical Sciences, Key Laboratory of MEA(Ministry of Education), Shanghai Key Laboratory of PMMP, East China Normal University}, 
		city={shanghai},
		postcode={200241}, 
		country={China}}
	
	\author[a]{Haibiao Zheng}
	
	\fnmark[*]
	
	\ead{hbzheng@math.ecnu.edu.cn}
	
	
	\credit{Methodology, Revising the paper}
	

\cortext[a]{Corresponding author}



\begin{abstract}
In this paper, a meshless Generalized Finite Difference Method (GFDM) is proposed to deal with the Stokes-Darcy coupled problem with the Beavers-Joseph-Saffman (BJS) interface conditions. Some high order GFDMs are proposed to show the advantage of the high order GFDM for the Stokes-Darcy coupled problem, which is that the high order method has high order accuracy and the convergence order. Some Stokes-Darcy coupled problems with closed interfaces, which has more complex geometric shape, are given to show the advantage of the GFDM for the complex interface. The interface location has been changed to show the influence of the interface location for the Stokes-Darcy coupled problem. The BJS interface conditions has related to the partial derivatives of unknown variables and the GFDM has advantage in dealing with  the interface conditions with the jump of derivatives. Four numerical examples have been provided to verify the existence of the good performance of the GFDM for the Stokes-Darcy coupled problems, including that the simplicity, accuracy, and stability in static and moving systems. Especially, the GFDM has the tolerance of the large jump. The Neaumann boundary condition is used in numerical simulations. 
\end{abstract}

\maketitle

\section{Introduction}
The Stokes-Darcy coupled problem can accurately describe the motion of fluids in open spaces and the seepage process in porous media, and is widely used in the study of groundwater flow in aquifers, oil and gas extraction processes, and so on.
    
In the past few years, a number of efficient numerical methods have been proposed to solve the Stokes-Darcy coupled problem. For instance, the coupled finite element method[1,2], the stabilized mixed finite element method[3,4], the finite element based double or multiple mesh method[5-8], the Lagrangian multiplier method[9], the partitioned time stepping or operator splitting method[10-13], and the domain decomposition method[14-19]. In particular, Sun et al.[20,21] used the domain decomposition method to solve the fully-mixed Stokes-Darcy coupling problem, and further used the Robin type domain decomposition algorithm based on two-grid techniques to solve the Stokes-Darcy coupled problem. Chen et al.[22] analyzed the Stokes-Darcy coupled problem by using a parallel Robin-Robin type domain decomposition algorithm. Li et al.[23] used a stabilized finite element method based on two local Gauss integrations to solve the coupled Stokes-Darcy problem. Li et al.[24] used an augmented Cartsian grid method to solve Stokes-Darcy fluid-structure interaction problems with closed interfaces. These computational methods based on finite element methods  have contributed to the processing of Stokes-Darcy coupled models, but the generation and reconstruction of meshs have significantly  increased the computational workload significantly. Therefore, experts are actively exploring meshless methods to solve Stokes-Darcy coupled models and have made some progress.

There are several meshless methods that have been successfully proposed to solve the Stokes-Darcy coupled problem. For instance, Safarpoor et al. [25] proposed a coupled method called the Radial Basis Function-Finite Difference(RBF-FD) to solve the Stokes-Darcy equations. Pu et al. [26] proposed a meshless deep learning method based on a physics-informed neural network to solve coupled Stokes-Darcy equations with Bever-Joseph-Saffman interface conditions. Yue et al. [27] proposed a coupled deep neural networks to solve the time-dependent coupled Stokes-Darcy problems. Although these meshed and meshless methods have been developed for solving the Stokes-Darcy coupled problem, it is still a changing task to accurately, efficiently and stably solve this classical Stokes-Darcy coupled problem. In this work, we will consider the Generalized Finite Difference Method (GFDM) for solving the Stokes-Darcy coupled problems.

The GFDM, which is based on Taylor series expansions and weighted moving square (MLS) approximation, is truly free from mesh generation and numerical quadrature. The basic concept of the GFDM was proposed by Orkisz and Lizska [28] in the early eighties. The GFDM was improved and developed by Benito et al. [29,30] and then extended to various mechanical problems [31,32]. Later, the GFDM was developed by Fan et al. [33] and Gu et al. [34] to solve the ill-posed inverse problem. More recently, the GFDM has been applied to the interface problem in static and moving systems [35-38], the elasticity interface problem[39], the biharmonic interface problems[40], the Stokes equations[41] and the Stokes interface problem[42],  and so on. For especially, the Stokes equation and the Stokes interface problems used a mixed boundary condition (MBC) to deal with the pressure oscillations in the Stokes equations, and in this paper, we will adopt a pressure Poisson equation to reduce the pressure oscillations and use the GFDM to solve the Stokes-Darcy coupled problem.

The rest of the paper is structured as follows: Section 2 introduces the model of Stokes-Darcy coupled problem. Section 3 presents the procedure of the GFDM, the handling skill for the Stokes-Darcy coupled problem, and the GFDM for the Stokes-Darcy coupled problem. In Section 4, convergence analysis is proposed. In Section 5, four numerical examples are presented to verify the accuracy, high efficiency and stability of the proposed methods. Finally, a conclusion is presented in Section 6.
\section{The Stokes-Darcy coupled problem}
We consider the Stokes equations giverned by the fluid velocity $\textbf{u}_f(x,y)$ and pressure $p(x,y)$  in \ $\Omega_f$:
\begin{align}
	- \nabla \cdot T (\textbf{u}_f, p_f)&=\textbf{f}_f,   \quad \quad \quad in \ \Omega_f,\\
	 \nabla \cdot \textbf{u}_f &=  0, \quad \quad \quad  in \ \Omega_f,
\end{align}
where $T=-pI+2vD(\textbf{u}_f)$ denotes the stress tensor, and $D(\textbf{u}_f)=\frac{1}{2}(\nabla \textbf{u}_f+(\nabla \textbf{u}_f)^T)$ represents the deformation tensor, in which $v$ is the kinematic viscosity of the fluid, $\textbf{f}_f$ is the external body force.

 The Darcy equations governed by the Darcy velocity $\textbf{u}_p (x,y)$ and the piezometric head $\phi(x,y)$ in $\Omega_p$:
 \begin{align}
 	\textbf{u}_p &=-K \nabla \phi,   \quad \quad \quad in \ \Omega_p,\\
 	\nabla \cdot \textbf{u}_p &= \textbf{f}_p,  \quad \quad \quad \quad \quad  in \ \Omega_p,
 \end{align}
 in which $K$ denotes the hydraulic conductivity tensor, $\textbf{f}_p$ is the external body force. And the Beavers-Joseph-Saffman(-Jones) coupling conditions on the interface $\Gamma$:
  \begin{align}
 	\textbf{u}_f \cdot \textbf{n}_f+\textbf{u}_p \cdot \textbf{n}_p &=0, \quad \quad  \quad \quad \quad  \quad \quad \quad  \quad \quad \quad  \quad \quad \quad \quad on \ \Gamma,\\
 	\textbf{n}_f \cdot T(\textbf{u}_f, p_f) \cdot  \textbf{n}_f&=g\phi, \quad \quad \quad  \quad \quad \quad  \quad \quad \quad  \quad \quad \quad  \quad \quad \quad on \ \Gamma,\\
 	\textbf{n}_f \cdot T(\textbf{u}_f, p_f)\cdot \tau _i&=\frac{v \alpha \sqrt{d}}{\sqrt{trace(\Pi)}} \textbf{u}_f \cdot \tau_i, 1\leq i \leq d-1,  \quad \quad \quad on \ \Gamma,
 \end{align}
where $\textbf{n}_f$ and $\textbf{n}_p$ denote the unit outer normal to the fluid and the porous media regions at the interface $\Gamma$, respectively. $\tau_i(i=1,2,...,d-1)$ denote  mutually orthogonal unit tangential vectors to the interface $\Gamma$ and the constant parameter $\alpha$ depends on $v$ and $K$. For simplicity, we assume that the hydraulic head $\phi$ and the fluid velocity $\textbf{u}_f$ satisfy a homogeneous Dirichlet condition except on $\Gamma$, unless they are stated otherwise:
 \begin{align}
	\phi &=0,   \quad \quad \quad on \ \partial \Omega_p \textbackslash \Gamma,\\
	\textbf{u}_f &=0,   \quad \quad \quad on \ \partial \Omega_f \textbackslash \Gamma.
\end{align}
We can see that the Eq.(1) can be rewritten as the follows
\begin{eqnarray} 
	-\nabla \cdot(v(\nabla \mathbf{u}_f+(\nabla \mathbf{u}_f)^T))+\nabla\cdot p&=&\mathbf{f}_f,\  \  in \ \Omega_f.
\end{eqnarray}
Note that Eq.(2) doesn't have any information about the pressure $p$. Therefore, we use the classical pressure poisson equation from [43], which finds divergence on both sides of the equal sign for Eq.(10), then
\begin{eqnarray} 
	-\nabla \cdot(v(\Delta\mathbf{u}_f+(\Delta \mathbf{u}_f)^T))+	\Delta p&=&\nabla \cdot\mathbf{f}_f,\  \  in \ \Omega_f.
\end{eqnarray}
Due to Eq.(2), the above equation can be simplified as follows
\begin{eqnarray}
	\Delta p &=&\nabla \cdot \mathbf{f}_f,\quad \quad in \ \Omega_f,
\end{eqnarray}
Therefore, the Eq.(2) can be exchange into the following part
\begin{align}
	\Delta p &=\nabla \cdot \mathbf{f}_f,\quad \quad in \ \Omega_f,\\
	\nabla \cdot \mathbf{u}_f&=0,\quad\quad\quad\quad on \  \partial \Omega_f. 
\end{align}
In the boundary part of the above equations, we adopt the following scheme
\begin{eqnarray}
	\nabla \cdot \mathbf{u}_f+p-p&=&0,\quad\quad\quad\quad on \  \partial \Omega_f, 
\end{eqnarray}
then a new boundary condition can be obtained
\begin{eqnarray}
	\nabla \cdot \mathbf{u}_f+p&=&p,\quad\quad\quad\quad on \  \partial \Omega_f.
\end{eqnarray}

For simplicity, the governing equation and the interface condition can be described by matrix components as follows:\\
Stokes equation:
\begin{align}
	-2v\frac{\partial^2 u_1(x,y)}{\partial x^2}-v\frac{\partial^2 u_1(x,y)}{\partial y^2}-v\frac{\partial^2 u_2(x,y)}{\partial x \partial y}+\frac{\partial p(x,y)}{\partial x}&=f_1(x,y),\quad \quad \quad in \  \Omega_f,\\
	-v\frac{\partial^2 u_2(x,y)}{\partial x^2}-2v\frac{\partial^2 u_2(x,y)}{\partial y^2}-v\frac{\partial^2 u_1(x,y)}{\partial x \partial y}+\frac{\partial p(x,y)}{\partial y}&=f_2(x,y),\quad \quad \quad in \  \Omega_f,\\
	\frac{\partial^2 p(x,y)}{\partial x^2}+\frac{\partial^2 p(x,y)}{\partial y^2}&=0, \quad \quad \quad \quad \quad \quad in \  \Omega_f,
\end{align}
Darcy equation:
 \begin{align}
	-K\frac{\partial^2 \phi(x,y)}{\partial x^2}-K\frac{\partial^2 \phi(x,y)}{\partial y^2}={f}_p(x,y),   \quad \quad \quad in \ \Omega_p,
\end{align}
with interface condition:
 \begin{align}
u_1(x,y) n_{f_1}+u_2(x,y) n_{f_2}-K\frac{\partial \phi(x,y)}{\partial x} n_{p_1}-K\frac{\partial \phi(x,y)}{\partial y} n_{p_2}=0, \quad \quad \quad \quad \quad \quad \quad \quad \quad \quad on \ \Gamma,\\
p(x,y)-2v\frac{\partial u_1(x,y)}{\partial x} (n_{f_1})^2-2v\frac{\partial u_1(x,y)}{\partial y} n_{f_2}n_{f_1}-2v\frac{\partial u_2(x,y)}{\partial x} n_{f_1}n_{f_2}-2v\frac{\partial u_2(x,y)}{\partial y}(n_{f_2})^2-g\phi(x,y)
\end{align}
$$=0\quad \quad \quad \quad \quad \quad \quad \quad \quad \quad \quad \quad \quad \quad \quad \quad \quad \quad  \quad \quad \quad \quad \quad \quad \quad \quad \quad  \quad \quad \quad \quad \quad \quad \quad \quad \quad on \ \Gamma,$$
 \begin{equation}
-2v\frac{\partial u_1(x,y)}{\partial x} n_{f_1} \tau_{i_1}-v\frac{\partial u_1(x,y)}{\partial y} n_{f_1}\tau_{i_2}-v\frac{\partial u_1(x,y)}{\partial y} n_{f_2}\tau_{i_1}-v\frac{\partial u_2(x,y)}{\partial x}n_{f_2}\tau_{i_1}-v\frac{\partial u_2(x,y)}{\partial x}n_{f_1}\tau_{i_2}
\end{equation}
$$-2v\frac{\partial u_2(x,y)}{\partial y} n_{f_2} \tau_{i_2}-\frac{v \alpha \sqrt{d}}{\sqrt{trace(\Pi)}}u_1(x,y)\tau_{i_1}-\frac{v \alpha \sqrt{d}}{\sqrt{trace(\Pi)}}u_2(x,y)\tau_{i_2}=0,\quad \quad \quad \quad \quad \quad \quad \quad on \ \Gamma,$$
and the Dirichlet boundary condition:
 \begin{align}
	\phi(x,y) &=0,   \quad \quad \quad on \ \partial \Omega_p \textbackslash \Gamma,\\
	u_1(x,y) &=0,   \quad \quad \quad on \ \partial \Omega_f \textbackslash \Gamma,\\
	u_2(x,y) &=0,   \quad \quad \quad on \ \partial \Omega_f \textbackslash \Gamma,\\
    \frac{\partial u_1(x,y)}{\partial x}+\frac{\partial u_2(x,y)}{\partial y}+p(x,y)&=p(x,y), \quad on \ \partial \Omega_f,
\end{align}
\begin{figure}
	\centering
	\includegraphics[scale=.5]{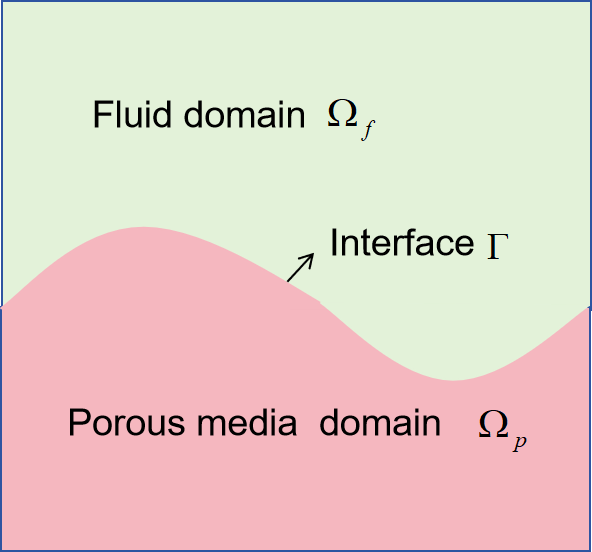} 
	   \hspace{1cm}
	\includegraphics[scale=.5]{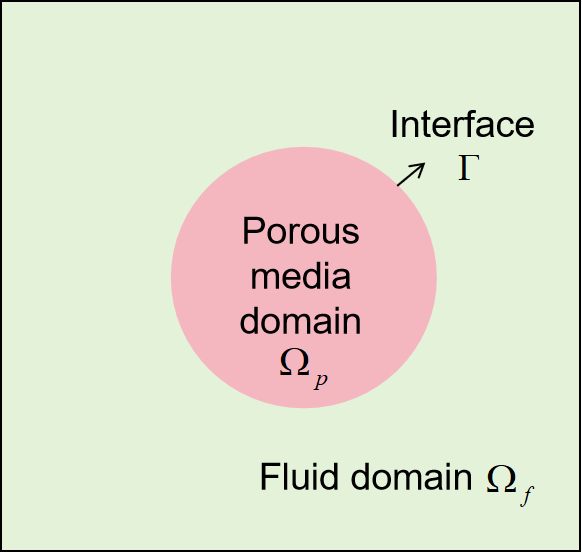}    
	\caption{ The distributions of the Stokes-Darcy coupled problem with unclosed interface(left) and closed interface(right).}
\end{figure}


\section{Numerical Schemes }
In this section, we briefly describe the numerical scheme of the GFDM, the  Domain Decomposition scheme and the GFDM for the Stokes-Darcy coupled problem are provided to detailed descripte the process of the numerical scheme.   
\subsection{The generalized finite difference method}
For a given interior node $(x_0,y_0)$,we select $m$ nearest nodes$(x_k,y_k)(k=1,2,\cdots,m)$and make these $m+1$ nodes form a star, as shown in Fig.2. Taking a 2D case for example, let $u_0=u(x_0,y_0)$ is the function value at the node $(x_0,y_0)$ and $u_k=u(x_k,y_k)(k=1,2,\cdots,m)$ are function values at nodes $(x_k,y_k)(k=1,2,\cdots,m)$. To expand the values of $u_k=u(x_k,y_k)(k=1,2,\cdots,m)$ around the given node $(x_0,y_0)$ using the Taylor series:
\begin{equation}
	u_k=u_0+h_k\frac{\partial u_0}{\partial x}+l_k\frac{\partial u_0}{\partial y}+\frac{{h_k}^2}{2}\frac{\partial^2 u_0}{\partial x^2}+\frac{{l_k}^2}{2}\frac{\partial^2 u_0}{\partial y^2}+h_kl_k\frac{\partial^2 u_0}{\partial x\partial y}
	+O(\rho^3),k=1,\cdots,m,
\end{equation}
where $h_k$ and $l_k$ are the distance between the node $(x_k,y_k)$ and $(x_0,y_0)$:
\begin{equation}
	h_k=x_k-x_0,k=1,2,\cdots,m,
\end{equation}
\begin{equation}
	l_k=y_k-y_0,k=1,2,\cdots,m.
\end{equation}
$\frac{\partial u_0}{\partial x},\frac{\partial u_0}{\partial y},\frac{\partial^2 u_0}{\partial x^2},\frac{\partial^2 u_0}{\partial y^2},\frac{\partial^2 u_0}{\partial x\partial y}$ are partial derivatives of each order at the node$(x_0,y_0)$. By truncating the Taylor series expression (Eq.(28)) after the second order derivatives, a second order GFDM can be obtained, then a residual function $B_2(u)$ can be obtained,
\begin{equation}
	B_2(u)=\sum_{k=1}^m[(u_0-u_k+h_k\frac{\partial u_0}{\partial x}+l_k\frac{\partial u_0}{\partial y}+\frac{{h_k}^2}{2}\frac{\partial^2 u_0}{\partial x^2}+\frac{{l_k}^2}{2}\frac{\partial^2 u_0}{\partial y^2}+h_kl_k\frac{\partial^2 u_0}{\partial x\partial y})\omega_k]^2,
\end{equation}
in which $\omega_k(k=1,2,\cdots,m)$ is the weighting coefficient at $(x_k,y_k)$. the weight function is taken as follows:
\begin{equation}
	\omega_k=\left\{
	\begin{array}{cl}
		{1-6{(\frac{d_k}{d_m})}^2+8{(\frac{d_k}{d_m})}^3-3{(\frac{d_k}{d_m})}^4},{d_k\leq d_m},\\{0},{d_k\geq d_m}.
	\end{array}
	\right.
\end{equation}
Here $d_k=\sqrt{(x_k-x_0)^2+(y_k-y_0)^2}$ is the distance between nodes $(x_k,y_k)$ and $(x_0,y_0)$, and $d_m$ is the maximum value of $d_k(k=1,2,\cdots,m)$. From the above weight function, we can see that the weighting coefficient is inversely proportional to the distance from the corresponding node $(x_k,y_k)$ to $(x_0,y_0)$.
\begin{figure}
	\centering
	\includegraphics[width=0.3\textwidth]{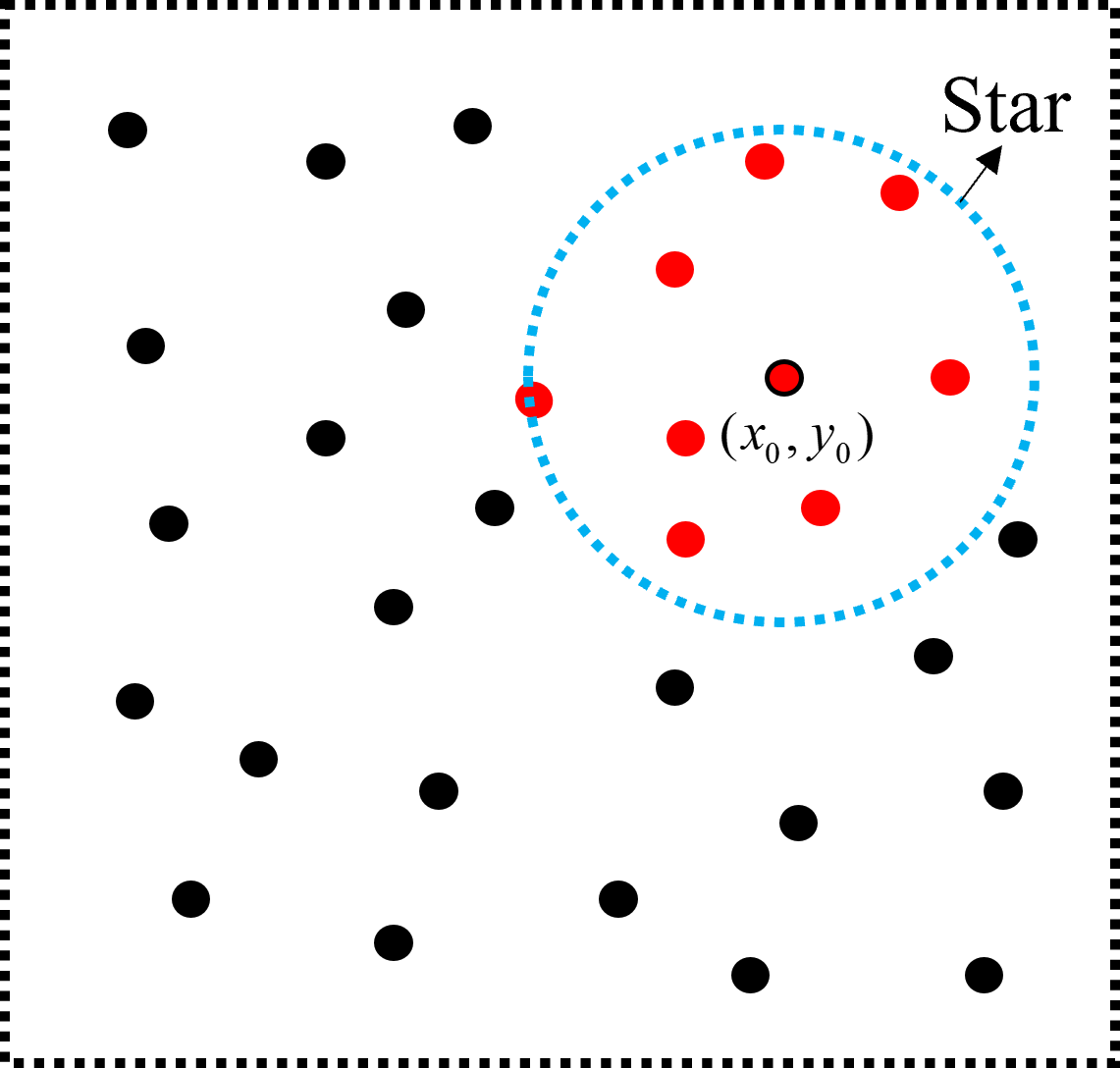}
	\caption{ The star of a central node $(x_0,y_0)$ for $2D$ case.}
\end{figure}
To minimize the above residual function concerning the partial derivatives $D(u)={(\frac{\partial u_0}{\partial x},\frac{\partial u_0}{\partial y},\frac{\partial^2 u_0}{\partial x^2},\frac{\partial^2 u_0}{\partial y^2},\frac{\partial^2 u_0}{\partial x\partial y})}^T$, that is, let
\begin{equation}
	\frac{\partial B_2(u)}{\partial {D_u}}=0,
\end{equation}
then the following linear equation system can be obtained,
\begin{equation}
	AD_u=b,
\end{equation}
where $A$ is a symmetric matrix given by
\begin{equation}
	A=\sum_{{k=1}} ^m {S_{k}^TS_{k}},
\end{equation}
here $S_{k}$ is a diagonal matrice. and its diagonal elements can defined as
\begin{equation}
	S_{k}=\omega_k(h_k,l_k,\frac{h_k^2}{2},\frac{l_k^2}{2},h_kl_k)^T,
\end{equation}
and
\begin{equation}
	b=BU=(-\sum_{k=1}^m \omega_kS_{(k)},\omega_1S_{(1)},\omega_2S_{(2)},\cdots,\omega_mS_{(m)})_{5\times(m+1)}(u_0,u_1,\cdots,u_m)^T,
\end{equation}
in which $U=(u_0,u_1,\cdots,u_m)^T$ are function values of all nodes inside the star. According to Eq. (34) and Eq. (37), the partial derivative vector $D_u$ can be expressed as follows:
\begin{equation}
	D_u=A^{-1}b=A^{-1}(BU)=(A^{-1}B)U=EU,
\end{equation}
where
\begin{equation}
	E=A^{-1}B,
\end{equation}
is generally obtained by using the moving least-squares method.
Based on the above analysis, the derivative at $(x_0,y_0)$ can be expressed as a linear combination of the function values for itself and its $m$ neighbouring nodes, as shown in Eq. (38). For example,
\begin{equation}
	\frac{\partial^2 u_0}{\partial y^2}=\sum_{j=1}^m E_{4j}u_j^0.
\end{equation}
It means that the second order derivative of the function value at the $(x_0,y_0)$ can to be written as the product of the fourth row of the matrix  formed by Eq. (39) and the matrix formed by the function values of the nodes inside the star. Repeating the same procedure for each node within the computational domain, we can obtain the algebraic equations for other interior points.
\subsection{The handling skill for the Stokes-Darcy coupled problem}
In this subsection, the idea of the Domain Decomposition is used to divide the transformed Stokes-Darcy coupled problem. By this scheme, the original Stokes/Darcy coupled problem (Eqs.(1),(3)-(9),(13),(16)) can be divided into two non-interface subproblems:
\begin{align}
	- \nabla \cdot T (\textbf{u}_f, p_f)&=\textbf{f}_f,\quad   \quad in \ \Omega_f,\\
	\Delta p &=\nabla \cdot \mathbf{f}_f, \quad in \ \Omega_f,\\
	  \textbf{u}_f &=0,   \quad \quad \quad on \ \partial \Omega_f \textbackslash \Gamma,\\
	\mathbf{u}_f+p(x,y)&=p(x,y),\quad on \  \partial \Omega_f,\\
    \nabla \cdot  
    \textbf{u}_f \cdot \textbf{n}_f+\textbf{u}_p \cdot \textbf{n}_p &=0,  \quad \quad \quad on \ \Gamma,\\
    \textbf{n}_f \cdot T(\textbf{u}_f, p_f) \cdot  \textbf{n}_f&=g\phi, \quad \quad on \ \Gamma, 
\end{align}
and 
\begin{align}
	\textbf{u}_p &=-K \nabla \phi,   \quad \quad \quad in \ \Omega_p,\\
	\nabla \cdot \textbf{u}_p &= \textbf{f}_p,  \quad \quad \quad \quad \quad  in \ \Omega_p,\\
		\phi &=0,   \quad \quad \quad \quad \quad on \ \partial \Omega_p \textbackslash \Gamma,\\
			\textbf{n}_f \cdot T(\textbf{u}_f, p_f)\cdot \tau _i&=\frac{v \alpha \sqrt{d}}{\sqrt{trace(\Pi)}} \textbf{u}_f \cdot \tau_i, 1\leq i \leq d-1,  \quad \quad \quad on \ \Gamma.
\end{align}
Due to the above problems are related to the matrix, these problems can be described by matrix components:
\begin{align}
	-2v\frac{\partial^2 u_1(x,y)}{\partial x^2}-v\frac{\partial^2 u_1(x,y)}{\partial y^2}-v\frac{\partial^2 u_2(x,y)}{\partial x \partial y}+\frac{\partial p(x,y)}{\partial x}=f_1(x,y),\quad \quad \quad in \  \Omega_f,\\
	-v\frac{\partial^2 u_2(x,y)}{\partial x^2}-2v\frac{\partial^2 u_2(x,y)}{\partial y^2}-v\frac{\partial^2 u_1(x,y)}{\partial x \partial y}+\frac{\partial p(x,y)}{\partial y}=f_2(x,y),\quad \quad \quad in \  \Omega_f,\\
	\frac{\partial^2 p(x,y)}{\partial x^2}+\frac{\partial^2 p(x,y)}{\partial y^2}=0, \quad \quad \quad \quad \quad \quad in \  \Omega_f,\\
		u_1(x,y) =0,   \quad \quad \quad on \ \partial \Omega_f \textbackslash \Gamma,\\
	u_2(x,y) =0,   \quad \quad \quad on \ \partial \Omega_f \textbackslash \Gamma,\\
	\frac{\partial u_1(x,y)}{\partial x}+\frac{\partial u_2(x,y)}{\partial y}+p(x,y)=p(x,y), \quad on \ \partial \Omega_f,\\
	u_1(x,y) n_{f_1}+u_2(x,y) n_{f_2}-K\frac{\partial \phi(x,y)}{\partial x} n_{p_1}-K\frac{\partial \phi(x,y)}{\partial y} n_{p_2}=0, \quad \quad \quad \quad \quad \quad \quad \quad \quad \quad on \ \Gamma,\\
	p(x,y)-2v\frac{\partial u_1(x,y)}{\partial x} (n_{f_1})^2-2v\frac{\partial u_1(x,y)}{\partial y} n_{f_2}n_{f_1}-2v\frac{\partial u_2(x,y)}{\partial x} n_{f_1}n_{f_2}-2v\frac{\partial u_2(x,y)}{\partial y}(n_{f_2})^2-g\phi(x,y)
	\end{align}
	$$=0\quad \quad \quad \quad \quad \quad \quad \quad \quad \quad \quad \quad \quad \quad \quad \quad \quad \quad  \quad \quad \quad \quad \quad \quad \quad \quad \quad  \quad \quad \quad \quad \quad \quad \quad \quad \quad on \ \Gamma,$$
and
 \begin{align}
	-K\frac{\partial^2 \phi(x,y)}{\partial x^2}-K\frac{\partial^2 \phi(x,y)}{\partial y^2}={f}_p(x,y),   \quad in \ \Omega_p,\\
		\phi(x,y) = 0,   \quad \quad \quad \quad \quad  on \  \partial \Omega_p \textbackslash \Gamma,
		 \end{align}
	 \begin{equation}
		-2v\frac{\partial u_1(x,y)}{\partial x} n_{f_1} \tau_{i_1}-v\frac{\partial u_1(x,y)}{\partial y} n_{f_1}\tau_{i_2}-v\frac{\partial u_1(x,y)}{\partial y} n_{f_2}\tau_{i_1}-v\frac{\partial u_2(x,y)}{\partial x}n_{f_2}\tau_{i_1}-v\frac{\partial u_2(x,y)}{\partial x}n_{f_1}\tau_{i_2}
	\end{equation}
	$$-2v\frac{\partial u_2(x,y)}{\partial y} n_{f_2} \tau_{i_2}-\frac{v \alpha \sqrt{d}}{\sqrt{trace(\Pi)}}u_1(x,y)\tau_{i_1}-\frac{v \alpha \sqrt{d}}{\sqrt{trace(\Pi)}}u_2(x,y)\tau_{i_2}=0,\quad \quad \quad \quad \quad \quad \quad \quad on \ \Gamma,$$
 From the above analysis, these subproblem (Eq.s(51)-(61)) are equivalent to the original interface problem (Eqs.(1),(3)-(9),(13),(16)). 
\subsection{The GFDM for the Stokes-Darcy coupled problem}
In this part, we use the GFDM to analyze the Stokes-Darcy coupled problem Eqs.(51)-(61). Firstly, $N^+_{inp}+N^-_{inp}, N^+_{bp}+N^-_{bp}$ and$ N^+_{\Gamma}+N^-_{\Gamma}$  scattered nodes can be obtained by discreting the computational domain $\Omega=\Omega^+\cup\Omega^-$, the domain boundary $\partial \Omega = \partial \Omega^+\cup \partial \Omega^-$ and the interface $\Gamma =\Gamma ^+\cup \Gamma ^-$, respectively. In particular, the computational domain $\Omega=\Omega^+(t)\cup\Omega^-(t)$ , the domain boundary $\partial \Omega(t) = \partial \Omega^+(t)\cup \partial \Omega^-(t)$  and the interface  $\Gamma(t) =\Gamma ^+(t)\cup \Gamma ^-(t)$ are considered when the interface is moving. $\Gamma(t)$ is a sufficiently smooth closed curve. For simplicity, we assume that $\Gamma(t)$ does not touch the boundary. We assume that time $t$ is a uniform distribution of points and discrete the time $t$  with $N_j^t$ nodes. That is let $dt=(T-0)./N_j^t$ then $t=t_0+(j-1)dt$. If we fix on $t, t=t_j, \Gamma(t)=\Gamma _j$, the $j_{th}$ location of the interface can be obtained. This means that we can get the exact interface information at $t_j$. Specially, the selection of interface points and the form of the star of the interface points can refer the article[36].

To enforce that the interior node satisfies the governing equation, the boundary node satisfies the boundary condition and the interface node satisfies the interface condition, $N^+_{inp}+N^-_{inp}+N^+_{bp}+N^-_{bp}+ N^+_{\Gamma}+N^-_{\Gamma}$linear algebraic equations can be obtained. To solve this algebraic system of equations, the result can be obtained. For simplicity, the final system of linear equations can be rewritten in matrix form:
	\begin{equation}
		\hat{A} X=b,
	\end{equation}
where
\begin{eqnarray}
\hat{A}=\left (
	\begin{array}{l|l}
		{F} &{D_1}\\
		\hline
        {D_2} &{D}
	\end{array}
	\right),
	U=\left (
	\begin{array}{ll}
		U^+\\ 
		U^-
	\end{array}
	\right),
	b=\left (
	\begin{array}{ll}
		b^+ \\
		b^-
	\end{array}
	\right),
\end{eqnarray}
Here
\begin{eqnarray}
		F=\left (
	\begin{array}{ll}
		{G} &{P_1}\\
		{P_2} &{P_p}
	\end{array}
	\right),
	D_1=\left (
	\begin{array}{l}
		{D_{11}}\\
		{D_{12}}\\
		{D_{13}}
	\end{array}
	\right),
	D_2=\left (
	\begin{array}{lll}
		{D_{21}} &{D_{22}}&{D_{23}}
	\end{array}
	\right),
	U^+=\left (
\begin{array}{ll}
	U_1^+\\ 
	U_2^+\\
	P\\
\end{array}
\right),
	b^+=\left (
\begin{array}{ll}
	b_1^+ \\
	b_2^+\\
	b_3^+
\end{array}
\right),
\end{eqnarray}
where
\begin{eqnarray}
		G=\left (
\begin{array}{llll}
	{G_{11}} &{G_{12}}\\
	{G_{21}} &{G_{22}}
\end{array}
\right),	
P_1=\left (
\begin{array}{l}
	{P_{11}}\\
	{P_{12}}
\end{array}
\right),
P_2=\left (
\begin{array}{ll}
	{P_{21}} &{P_{22}}
\end{array}
\right),
\end{eqnarray}
\begin{center}
	\begin{equation}
		$$
		U_1^+=(u^+_{1,1},\cdots, u^+_{1,N^+_{bp}},u^+_{1,N^+_{bp}+1},\cdots, u^+_{1,N^+_{bp}+N^+_{\Gamma}},u^+_{1,N^+_{bp}+N^+_{\Gamma}+1},\cdots, u^+_{1,N^+_{bp}+N^+_{\Gamma}+N^+_{inp}})^T,
		$$
	\end{equation}
\end{center}
\begin{center}
	\begin{equation}
		$$
		U_2^+=(u^+_{2,1},\cdots, u^+_{2,N^+_{bp}},u^+_{2,N^+_{bp}+1},\cdots u^+_{2,N^+_{bp}+N^+_{\Gamma}},u^+_{2,N^+_{bp}+N^+_{\Gamma}+1},\cdots, u^+_{2,N^+_{bp}+N^+_{\Gamma}+N^+_{inp}})^T,
		$$
	\end{equation}
\end{center}
\begin{center}
	\begin{equation}
		$$
		P=(p_{1},\cdots, p_{N^+_{bp}},p_{N^+_{bp}+1},\cdots p_{N^+_{bp}+N^+_{\Gamma}},p_{N^+_{bp}+N^+_{\Gamma}+1},\cdots, p_{N^+_{bp}+N^+_{\Gamma}+N^+_{inp}})^T,
		$$
	\end{equation}
\end{center}
\begin{center}
	\begin{equation}
		$$
		U^-=(\phi_1,\cdots, \phi_{N^-_{bp}},\phi_{N^-_{bp}+1},\cdots, \phi_{N^-_{bp}+N^-_{\Gamma}},\phi_{N^-_{bp}+N^-_{\Gamma}+1},\cdots, \phi_{N^-_{bp}+N^-_{\Gamma}+N^-_{inp}})^T,
		$$
	\end{equation}
\end{center}
\begin{center}
	\begin{equation}
		$$
		b_1^+=(0,\cdots, 0,0,\cdots, 0,f_{1,N^+_{bp}+N^+_{\Gamma}+1},\cdots, f_{1,N^+_{bp}+N^+_{\Gamma}+N^+_{inp}})^T,
		$$
	\end{equation}
\end{center}
\begin{center}
	\begin{equation}
		$$
		b_2^+=(0,\cdots, 0,0,\cdots, 0,f_{2,N^+_{bp}+N^+_{\Gamma}+1},\cdots, f_{2,N^+_{bp}+N^+_{\Gamma}+N^+_{inp}})^T,
		$$
	\end{equation}
\end{center}
\begin{center}
	\begin{equation}
		$$
		b_3^+=(p_{1},\cdots, p_{N^+_{bp}},p_{N^+_{bp}+1},\cdots p_{N^+_{bp}+N^+_{\Gamma}},0,\cdots, 0)^T,
		$$
	\end{equation}
\end{center}
\begin{center}
	\begin{equation}
		$$
		b^-=(0,\cdots, 0,0,\cdots, 0,f_{p,N^-_{bp}+N^-_{\Gamma}+1},\cdots, f_{p,N^-_{bp}+N^-_{\Gamma}+N^-_{inp}})^T,
		$$
	\end{equation}
\end{center}
and
\begin{eqnarray}
	G_{11}=\left (
\begin{array}{lll}
	{E_{N^+_{bp}\times{N^+_{bp}}}} &{0}&{0}\\
	{0} &{K_{1,N^+_{\Gamma}\times{N^+_{\Gamma}}}} &{0}\\
	{0} &{0} &{I_{1,N^+_{inp}\times N^+_{inp}}}
	\end{array}
	\right),
	G_{12}=\left (
	\begin{array}{lll}
		{0} &{0}&{0}\\
		{0} &{K_{2,N^+_{\Gamma}\times{N^+_{\Gamma}}}} &{0}\\
		{0} &{0} &{I_{2,N^+_{inp}\times N^+_{inp}}}
		\end{array}
		\right),
			\end{eqnarray}
\begin{eqnarray}
		G_{21}=\left (
\begin{array}{lll}
	{0} &{0}&{0}\\
	{0} &{K_{3,N^+_{\Gamma}\times{N^+_{\Gamma}}}} &{0}\\
	{0} &{0} &{I_{3,N^+_{inp}\times N^+_{inp}}}
	\end{array}
	\right),
	G_{22}=\left (
	\begin{array}{lll}
		{E_{N^+_{bp}\times{N^+_{bp}}}} &{0}&{0}\\
		{0} &{K_{4,N^+_{\Gamma}\times{N^+_{\Gamma}}}} &{0}\\
		{0} &{0} &{I_{4,N^+_{inp}\times N^+_{inp}}}
		\end{array}
		\right),
\end{eqnarray}
\begin{eqnarray}
	P_{11}=\left (
	\begin{array}{lll}
		{L_{1,N^+_{bp}\times{N^+_{bp}}}} &{0}&{0}\\
		{0} &{L_{2,N^+_{\Gamma}\times{N^+_{\Gamma}}}} &{0}\\
		{0} &{0} &{0}
	\end{array}
	\right),
	P_{12}=\left (
	\begin{array}{lll}
		{L_{3,N^+_{bp}\times{N^+_{bp}}}} &{0}&{0}\\
		{0} &{L_{4,N^+_{\Gamma}\times{N^+_{\Gamma}}}} &{0}\\
		{0} &{0} &{0}
	\end{array}
	\right),
\end{eqnarray}
\begin{eqnarray}
	P_{21}=\left (
	\begin{array}{lll}
		{0} &{0}&{0}\\
		{0} &{0} &{0}\\
		{0} &{0} &{R_{1,N^+_{inp}\times N^+_{inp}}}
	\end{array}
	\right),
	P_{22}=\left (
	\begin{array}{lll}
		{0} &{0}&{0}\\
		{0} &{E_{N^+_{\Gamma}\times{N^+_{\Gamma}}}} &{0}\\
		{0} &{0} &{R_{2,N^+_{inp}\times N^+_{inp}}}
	\end{array}
	\right),
\end{eqnarray}
\begin{eqnarray}
	P_{p}=\left (
	\begin{array}{lll}
		{E_{N^+_{bp}\times{N^+_{bp}}}} &{0}&{0}\\
		{0} &{E_{N^+_{\Gamma}\times{N^+_{\Gamma}}}} &{0}\\
		{0} &{0} &{P_{1,N^+_{inp}\times N^+_{inp}}}
	\end{array}
	\right),
\end{eqnarray}
\begin{eqnarray}
	D_{11}=\left (
	\begin{array}{lll}
		{0} &{0}&{0}\\
		{0} &{M_{1,N^-_{\Gamma}\times{N^-_{\Gamma}}}} &{0}\\
		{0} &{0} &{0}
	\end{array}
	\right),
	D_{12}=\left (
	\begin{array}{lll}
		{0} &{0}&{0}\\
		{0} &{M_{2,N^-_{\Gamma}\times{N^-_{\Gamma}}}} &{0}\\
		{0} &{0} &{0}
	\end{array}
	\right),
	D_{13}=D_{23}=0,
\end{eqnarray}
\begin{eqnarray}
	D_{21}=\left (
	\begin{array}{lll}
		{0} &{0}&{0}\\
		{0} &{M_{3,N^-_{\Gamma}\times{N^-_{\Gamma}}}} &{0}\\
		{0} &{0} &{0}
	\end{array}
	\right),
	D_{22}=\left (
	\begin{array}{lll}
		{0} &{0}&{0}\\
		{0} &{M_{4,N^-_{\Gamma}\times{N^-_{\Gamma}}}} &{0}\\
		{0} &{0} &{0}
	\end{array}
	\right),
		D=\left (
	\begin{array}{lll}
		{E_{N^-_{bp}\times{N^-_{bp}}}} &{0}&{0}\\
		{0} &{0} &{0}\\
		{0} &{0} &{M_{5,N^-_{inp}\times N^-_{inp}}}
	\end{array}
	\right),
\end{eqnarray}
in which
\begin{eqnarray}
	K_{1,N^+_{\Gamma}\times{N^+_{\Gamma}}}(i,j)&=&n_{f_1},\\
	K_{2,N^+_{\Gamma}\times{N^+_{\Gamma}}}(i,j)&=&n_{f_2},\\
    K_{3,N^+_{\Gamma}\times{N^+_{\Gamma}}}(i,j)&=&-2vE(1,i)n_{f_1}^2,\\
    K_{4,N^+_{\Gamma}\times{N^+_{\Gamma}}}(i,j)&=&-2vE(1,i)n_{f_1}n_{f_2},\\
	I_{1,N^+_{inp}\times N^+_{inp}}(i,j)&=&-2vE(3,i)-vE(4,i),\\
	I_{2,N^+_{inp}\times N^+_{inp}}(i,j)&=&-vE(5,i),\\
	I_{3,N^+_{inp}\times N^+_{inp}}(i,j)&=&-vE(5,i),\\
	I_{4,N^+_{inp}\times N^+_{inp}}(i,j)&=&-vE(3,i)-2vE(4,i),\\
	L_{1,N^+_{bp}\times{N^+_{bp}}}(i,j)&=&E(1,i),\\
	L_{3,N^+_{bp}\times{N^+_{bp}}}(i,j)&=&E(2,i),\\
	L_{2,N^+_{\Gamma}\times{N^+_{\Gamma}}}(i,j)&=&E(1,i),\\
	L_{4,N^+_{\Gamma}\times{N^+_{\Gamma}}}(i,j)&=&E(2,i),\\
	R_{1,N^+_{inp}\times N^+_{inp}}(i,j)&=&E(1,i),\\
	R_{2,N^+_{inp}\times N^+_{inp}}(i,j)&=&E(2,i),\\
	P_{1,N^+_{inp}\times N^+_{inp}}(i,j)&=&E(3,i)+E(4,i),\\
	M_{1,N^-_{\Gamma}\times{N^-_{\Gamma}}}(i,j)&=&-2vE(1,i)n_{f_1}\tau_{i_1}-vE(2,i)n_{f_1}\tau_{i_2}-vE(2,i)n_{f_2}\tau_{i_1}-\frac{v\alpha\sqrt{d}}{\sqrt{trace(\Pi)}}\tau_{i_1},\\
	M_{2,N^-_{\Gamma}\times{N^-_{\Gamma}}}(i,j)&=&-vE(1,i)n_{f_2}\tau_{i_1}-vE(1,i)n_{f_1}\tau_{i_2}-2vE(2,i)n_{f_2}\tau_{i_2}-\frac{v\alpha\sqrt{d}}{\sqrt{trace(\Pi)}}\tau_{i_2},\\
	M_{3,N^-_{\Gamma}\times{N^-_{\Gamma}}}(i,j)&=&-KE(1,i)n_{p_1}-KE(2,i)n_{p_2},\\
	M_{4,N^-_{\Gamma}\times{N^-_{\Gamma}}}(i,j)&=&-g,\\
	M_{5,N^-_{inp}\times N^-_{inp}}(i,j)&=&-KE(3,i)-KE(4,i).
\end{eqnarray}
 The matrices $\hat{A}$ is formed according to the idea  of the domain decomposition, $F$ and $D$ are formed by the governing equations (Eq.(51)-(53),Eq.(59)), the Dirichlet boundary conditions (Eq.(54)-(56),Eq.(60)) and the interface condition (Eq.(57)-(58),Eq.(61)), respectively. The interface parts of $F, D_1,D_2,D$ are all created by the interface condition (Eq.(57),Eq.(58),Eq.(61)). Obviously, it is simple to use the above matrices to describe the interface conditions and the interface informations are only used on the interface points, the GFDM has advantage in dealing with the complex interface and the moving interface. $P_1$,$P_2$ obtain the pressure information in Eq.(51)-(52). $P_p$ is formed according to the governing equation (Eq.(53)). we can obtain different matrix if we use different form which is used to enrich the information about the pressure $p$. We can see that the matrix $\hat{A}$ is sparse.
\section{Convergence analysis}
Theorem 1: The even ($2i$)-order generalized finite difference method has an even ($2i$)-order convergence order, in which $i \in Z ^+$.

Proof:

For $i=1$, we consider the fourth order Taylar series expanding, the residual function can be obtained \\
\begin{equation}\label{Eq:25}
	\begin{aligned}
		B_2^*(u)&\left.{=\sum_{k=1}^m[(u_0-u_k+h_k\frac{\partial u_0}{\partial x}+l_k\frac{\partial u_0}{\partial y}+\frac{{h_k}^2}{2}\frac{\partial^2 u_0}{\partial x^2}+\frac{{l_k}^2}{2}\frac{\partial^2 u_0}{\partial y^2}+h_kl_k\frac{\partial^2 u_0}{\partial x\partial y}+\frac{1}{6}(h_k\frac{\partial u_0}{\partial x}+l_k\frac{\partial u_0}{\partial y})^3}\right.\\
		&\left.{+\frac{1}{24}(h_k\frac{\partial u_0}{\partial x}+l_k\frac{\partial u_0}{\partial y})^4+\cdots)\omega_k]^2,}\right.
	\end{aligned}
\end{equation}
Regarding the second-order partial derivative $D_{2, u}$, minimizing the extended residual function $B^*(u)$
\begin{equation}\label{Eq:26}
	\frac{\partial B_2^*(u)}{\partial {D_{2,u}}}=0,
\end{equation}
then the following linear equation system can be obtained,
\begin{align}\label{Eq:27}
	AD_{2,u}=b_2^*,\\
	b_2^*=b_2-b_2^{**},
\end{align}
\begin{eqnarray}\label{Eq:28}
	b_2^{**}=\left (
	\begin{array}{lll}
		{\sum_{k=1}^m[\frac{h_k}{6}(h_k\frac{\partial u_0}{\partial x}+l_k\frac{\partial u_0}{\partial y})^3\omega_k^2+\frac{h_k}{24}(h_k\frac{\partial u_0}{\partial x}+l_k\frac{\partial u_0}{\partial y})^4\omega_k^2+\cdots]}\\
		{\sum_{k=1}^m[\frac{l_k}{6}(h_k\frac{\partial u_0}{\partial x}+l_k\frac{\partial u_0}{\partial y})^3\omega_k^2+\frac{l_k}{24}(h_k\frac{\partial u_0}{\partial x}+l_k\frac{\partial u_0}{\partial y})^4\omega_k^2+\cdots]} \\
		{...}\\
		{\sum_{k=1}^m[\frac{h_k{l_k}}{6}(h_k\frac{\partial u_0}{\partial x}+l_k\frac{\partial u_0}{\partial y})^3\omega_k^2+\frac{h_k{l_k}}{24}(h_k\frac{\partial u_0}{\partial x}+l_k\frac{\partial u_0}{\partial y})^4\omega_k^2+\cdots]}\\
	\end{array}
	\right).
\end{eqnarray}
To define the following function 
\begin{equation}\label{Eq:29}
	\begin{aligned}
		f_2(u)&\left.{=(\beta,\beta,...,\beta)_{1\times 5}D_{2,u}=(\beta,\beta,...,\beta)_{1\times 5}A^{-1}b_2-(\beta,\beta,...,\beta)_{1\times 5}A^{-1}b_2^{**}}\right.
	\end{aligned}
\end{equation}
According to the explicit form of the generalized finite difference method(From Ref[44]),
\begin{equation}\label{Eq:100}
	f_2(u)=-m_0+\sum_{k=1}^m{m_ku_k}
\end{equation}
in which
\begin{equation}
	u_0=\frac{1}{m_0}\sum_{k=1}^m{m_ku_k}, \sum_{k=1}^m{m_k}=m_0. 
\end{equation}
then
\begin{equation}\label{Eq:30}
	f_2(u)-(\beta,\beta,...,\beta)_{1\times 5}D_{2,u}=-(\beta,\beta,...,\beta)_{1\times 5}A^{-1}b_2^{**}=\varepsilon_{2s}.
\end{equation}
Similaring to equation (108) and (109) in reference[44], and let $m_k=\frac{\mu_{k}}{h_k^2+l_k^2}, \mu_{k} \in R, $ in which the definition of $m_k$ can be seen in Appendix A,  then
\begin{equation}\label{Eq:31}
	\begin{aligned}
		\varepsilon_{2s}&\left.{=-(\beta,\beta,...,\beta)_{1\times 5}A^{-1}\left (
			\begin{array}{lll}
				{\sum_{k=1}^m[\frac{1}{6}(h_k\frac{\partial u_0}{\partial x}+l_k\frac{\partial u_0}{\partial y})^3\omega_k^2h_k+\frac{1}{24}(h_k\frac{\partial u_0}{\partial x}+l_k\frac{\partial u_0}{\partial y})^4\omega_k^2h_k+\cdots]}\\
				{\sum_{k=1}^m[\frac{1}{6}(h_k\frac{\partial u_0}{\partial x}+l_k\frac{\partial u_0}{\partial y})^3\omega_k^2l_k+\frac{1}{24}(h_k\frac{\partial u_0}{\partial x}+l_k\frac{\partial u_0}{\partial y})^4\omega_k^2l_k+\cdots]} \\
				{...}\\
				{\sum_{k=1}^m[\frac{1}{6}(h_k\frac{\partial u_0}{\partial x}+l_k\frac{\partial u_0}{\partial y})^3\omega_k^2h_k{l_k}+\frac{1}{24}(h_k\frac{\partial u_0}{\partial x}+l_k\frac{\partial u_0}{\partial y})^4\omega_k^2h_k{l_k}+\cdots]}\\
			\end{array}
			\right)}\right.\\
		&\left.{=-\sum_{k=1}^m\mu_{k}[\frac{1}{6(h_k^2+l_k^2)}(h_k\frac{\partial u_0}{\partial x}+l_k\frac{\partial u_0}{\partial y})^3+\frac{1}{24(h_k^2+l_k^2)}(h_k\frac{\partial u_0}{\partial x}+l_k\frac{\partial u_0}{\partial y})^4]+o(h_k^3,l_k^3).}\right.\\
	\end{aligned}
\end{equation}
Therefore, the conclusion can be obtained when $i=1$. It means that the second order generalized finite difference method has an second order convergence order.

For $i=2$, we consider the sixth order Taylor series expanding, the residual function can be obtained
\begin{equation}\label{Eq:32}
	\begin{aligned}
		B_4^*(u)&\left.{=\sum_{k=1}^m[(u_0-u_k+h_k\frac{\partial u_0}{\partial x}+l_k\frac{\partial u_0}{\partial y}+\frac{{h_k}^2}{2}\frac{\partial^2 u_0}{\partial x^2}+\frac{{l_k}^2}{2}\frac{\partial^2 u_0}{\partial y^2}+h_kl_k\frac{\partial^2 u_0}{\partial x\partial y}+\frac{{h_k}^3}{6}\frac{\partial^3 u_0}{\partial x^3}+\frac{{l_k}^3}{6}\frac{\partial^3 u_0}{\partial y^3}}\right.\\
		&\left.{+\frac{{h_k}^2l_k}{2}\frac{\partial^3 u_0}{\partial x^2\partial y}+\frac{h_k{l_k}^2}{2}\frac{\partial^3 u_0}{\partial x\partial y^2}+\frac{{h_k}^4}{24}\frac{\partial^4 u_0}{\partial x^4}+\frac{{l_k}^4}{24}\frac{\partial^4 u_0}{\partial y^4}+\frac{{h_k^3}l_k}{6} \frac{\partial^4 u_0}{\partial x^3\partial y}+\frac{{h_k}^2{l_k}^2}{4}\frac{\partial^4 u_0}{\partial x^2\partial y^2}+\frac{{l_k}^3h_k}{6}\frac{\partial^4 u_0}{\partial x\partial y^3}}\right.\\
		&\left.{+\frac{1}{120}(h_k\frac{\partial u_0}{\partial x}+l_k\frac{\partial u_0}{\partial y})^5+\frac{1}{720}(h_k\frac{\partial u_0}{\partial x}+l_k\frac{\partial u_0}{\partial y})^6\omega_k+\cdots]^2}\right.
	\end{aligned}
\end{equation}

Regarding the fourth order partial derivative $D_4(u)=(\frac{\partial u_0}{\partial x},\frac{\partial u_0}{\partial y},\frac{\partial^2 u_0}{\partial x^2},\frac{\partial^2 u_0}{\partial y^2},$
$\frac{\partial^2 u_0}{\partial x\partial y},\frac{\partial^3 u_0}{\partial x^3},\frac{\partial^3 u_0}{\partial y^3},\frac{\partial^3 u_0}{\partial x^2\partial y},\frac{\partial^3 u_0}{\partial x\partial y^2},\frac{\partial^4 u_0}{\partial x^4},\frac{\partial^4 u_0}{\partial y^4},$\\
$\frac{\partial^4 u_0}{\partial x^3\partial y},\frac{\partial^4 u_0}{\partial x^2\partial y^2},\frac{\partial^4 u_0}{\partial x\partial y^3})^T$, minimizing the extended residual function $B_4^*(u)$
\begin{equation}\label{Eq:33}
	\frac{\partial B_4^*(u)}{\partial {D_{4,u}}}=0,
\end{equation}
then the following linear equation system can be obtained,
\begin{align}\label{Eq:34}
	AD_{4,u}=b^*,\\
	b_4^*=b_4-b_4^{**},
\end{align}
in which
\begin{equation}\label{Eq:35}
	\begin{aligned}
		&\left.{	b_4^{**}=\left (
			\begin{array}{lll}
				{\sum_{k=1}^m[\frac{1}{120}(h_k\frac{\partial u_0}{\partial x}+l_k\frac{\partial u_0}{\partial y})^5\omega_k^2h_k+\frac{1}{720}(h_k\frac{\partial u_0}{\partial x}+l_k\frac{\partial u_0}{\partial y})^6\omega_k^2h_k+\cdots]}\\
				{\sum_{k=1}^m[\frac{1}{120}(h_k\frac{\partial u_0}{\partial x}+l_k\frac{\partial u_0}{\partial y})^5\omega_k^2l_k+\frac{1}{720}(h_k\frac{\partial u_0}{\partial x}+l_k\frac{\partial u_0}{\partial y})^6\omega_k^2l_k+\cdots]} \\
				{...}\\
				{\sum_{k=1}^m[\frac{1}{120}(h_k\frac{\partial u_0}{\partial x}+l_k\frac{\partial u_0}{\partial y})^5\omega_k^2h_k{l_k}^3+\frac{1}{720}(h_k\frac{\partial u_0}{\partial x}+l_k\frac{\partial u_0}{\partial y})^6\omega_k^2h_k{l_k}^3+\cdots]}\\
			\end{array}
			\right)}\right.\\
	\end{aligned}
\end{equation}
Therefore,
\begin{equation}\label{Eq:36}
	\varepsilon_{4s}=-(\beta,\beta,...,\beta)_{1\times 14}A^{-1}b^{**}
\end{equation}
It is similar to the process at$i=1$, let $m_k=\frac{\mu_{k}}{h_k^2+l_k^2}, \mu_{k} \in R, $ in which the definition of $m_k$ is similar to Appendix A, and then 
\begin{equation}\label{Eq:37}
	\begin{aligned}
		\varepsilon_{4s}&\left.{=-\sum_{k=1}^mm_k[\frac{1}{120(h_k^2+l_k^2)}(h_k\frac{\partial u_0}{\partial x}+l_k\frac{\partial u_0}{\partial y})^5+\frac{1}{720(h_k^2+l_k^2)}(h_k\frac{\partial u_0}{\partial x}+l_k\frac{\partial u_0}{\partial y})^6]+o(h_k^5,l_k^5).}\right.
	\end{aligned}
\end{equation}
Therefore, the conclusion can be obtained when $i=2$. It means that the fourth order generalized finite difference method has an fourth order convergence order.
and so forth...

To suppose the conclusion can be obtained when $i=n-1$
\begin{equation}\label{Eq:38}
	\varepsilon_{(2n-2)s}=-(\beta,\beta,...,\beta)_{1\times N_{Partial}}A^{-1}b^{**}
\end{equation}
let $m_k=\frac{\mu_{k}}{h_k^2+l_k^2}, \mu_{k} \in R, $ in which the definition of $m_k$ is similar to Appendix A, and then
\begin{equation}\label{Eq:39}
	\begin{aligned}
		\varepsilon_{(2n-2)s}&\left.{=-\sum_{k=1}^m[\frac{1}{(2n-1)!}(h_k\frac{\partial u_0}{\partial x}+l_k\frac{\partial u_0}{\partial y})^{2n-1}+\frac{1}{(2n)!}(h_k\frac{\partial u_0}{\partial x}+l_k\frac{\partial u_0}{\partial y})^{2n}]+o(h_k^{2n-1},l_k^{2n-1}).}\right.
	\end{aligned}
\end{equation}
To truncate the terms of the Taylor series after 2n+2 orders when $i=n$, construct an extended residual function, and minimize the extended residual function, the error can be obtained:
\begin{equation}\label{Eq:40}
	\begin{aligned}
		&\left.{	\varepsilon_{2ns}=-(\beta,\beta,...,\beta)_{1\times N_{Partial}}A^{-1}}\right.\\
		&\left.{\left (
			\begin{array}{lll}
				{\sum_{k=1}^m[\frac{1}{(2n+1)!}(h_k\frac{\partial u_0}{\partial x}+l_k\frac{\partial u_0}{\partial y})^{2n+1}\omega_k^2h_k+\frac{1}{(2n+2)!}(h_k\frac{\partial u_0}{\partial x}+l_k\frac{\partial u_0}{\partial y})^{2n+2}\omega_k^2h_k+\cdots]}\\
				{\sum_{k=1}^m[\frac{1}{(2n+1)!}(h_k\frac{\partial u_0}{\partial x}+l_k\frac{\partial u_0}{\partial y})^{2n+1}\omega_k^2l_k+\frac{1}{(2n+2)!}(h_k\frac{\partial u_0}{\partial x}+l_k\frac{\partial u_0}{\partial y})^{2n+2}\omega_k^2l_k+\cdots]} \\
				{...}\\
				{\sum_{k=1}^m[\frac{1}{(2n+1)!}(h_k\frac{\partial u_0}{\partial x}+l_k\frac{\partial u_0}{\partial y})^{2n+1}\omega_k^2h_k{l_k^{(2n-1)}}+\frac{1}{(2n+2)!}(h_k\frac{\partial u_0}{\partial x}+l_k\frac{\partial u_0}{\partial y})^{2n+2}\omega_k^2h_k{l_k^{(2n-1)}}+\cdots]}\\
			\end{array}
			\right)}\right.\\
	\end{aligned}
\end{equation}
let $m_k=\frac{\mu_{k}}{h_k^2+l_k^2}, \mu_{k} \in R,$ in which the definition of $m_k$ is similar to Appendix A, and then
	\begin{equation}\label{Eq:41}
	\begin{aligned}	
		\varepsilon_{2ns}&\left.{=-\sum_{k=1}^m[\frac{1}{(2n+1)!}(h_k\frac{\partial u_0}{\partial x}+l_k\frac{\partial u_0}{\partial y})^{2n+1}+\frac{1}{(2n+2)!}(h_k\frac{\partial u_0}{\partial x}+l_k\frac{\partial u_0}{\partial y})^{2n+2}]+o(h_k^{2n+1},l_k^{2n+1}).}\right.
	\end{aligned}
\end{equation}

Therefore, the conclusion can be obtained when $i=n$. It means that the 2nth-order generalized finite difference method has an 2nth-order convergence order. In which, $N_{Partial}$ denotes the total number of partial derivatives of the corresponding order.

According to mathematical induction, the above conclusion holds true, the $2i$th-order generalized finite difference method has an $2i$th-order convergence order, in which $i \in Z ^+$.

\section{Numerical simulations}
In this section, four examples are provided to show the accuracy, stability and high efficiency of the meshless GFDM for solving the Stokes-Darcy coupled problem with the linear interface in Example 1-Example 2, the closed complex interface in Example 3 and the moving interface in Example 4. In Example 1, we use the fixed kinematic viscosity $v$ of the fluid and the element of the hydraulic conductivity tensor $K$. In Example 2, the varying physical coefficients $v$ and $K$ are used. 

For simplicity, we defined the error norms as follows:
\begin{eqnarray}
	L_\infty=max|u_i-u(x_i)|(i=1,2,\cdots, N_{T}),
\end{eqnarray}

\begin{eqnarray}
	L_2=[\sum_{i=1}^{N_{T}}{\frac{(u_i-u(x_i))^2}{N_{T}}}]^{\frac{1}{2}}(i=1,2,\cdots, N_{T}),
\end{eqnarray}
\begin{eqnarray}
	H^1=[\sum_{i=1}^{N_{T}}{\frac{(\nabla u_i-\nabla u(x_i))^2}{N_{T}}}]^{\frac{1}{2}}(i=1,2,\cdots, N_{T}).
\end{eqnarray}
and the relative errors are defined as follows
\begin{eqnarray}
	L_{\infty,relative}=\frac{max|u_i-u(x_i)|}{max|u(x_i)|}(i=1,2,\cdots, N_{T}),
\end{eqnarray}
\begin{eqnarray}
	L_{2,relative}=\frac{[\sum_{i=1}^{N_{T}}{\frac{(u_i-u(x_i))^2}{N_{T}}}]^{\frac{1}{2}}}{[\sum_{i=1}^{N_{T}}{\frac{u(x_i)^2}{N_{T}}}]^{\frac{1}{2}}}(i=1,2,\cdots, N_{T}),
\end{eqnarray}
\begin{eqnarray}
	H^1_{relative}=\frac{[\sum_{i=1}^{N_{T}}{\frac{(\nabla u_i-\nabla u(x_i))^2}{N_{T}}}]^{\frac{1}{2}}}{[\sum_{i=1}^{N_{T}}{\frac{(\nabla u(x_i))^2}{N_{T}}}]^{\frac{1}{2}}}(i=1,2,\cdots, N_{T}).
\end{eqnarray}
Note that the $H^1$ and $H^1_{relative}$ of the veocity of Darcy $\textbf{u}_p$ contains the first and second order derivatives of $\phi$, which can deduce that the numerical errors can illustrate the advantage of the GFDM for the derivative functions. $u_i$ and $u(x_i)$ are the numerical and exact solution at point $x_i$, respectively. Let the domain $\Omega$ is a rectangle, $\Omega_f$ is the upper domain and  $\Omega_p$ is the lower domain for linear interface $\Gamma$.  $\Omega_f$ outside the interface $\Gamma$ and $\Omega_p=\Omega/(\Omega_f\cup\Gamma_t)$ inside closed interface $\Gamma$. $N_{T}$ is the number of all scattered nodes in $\Omega_f$, $\Omega_p$, $\partial\Omega_f$,$\partial\Omega_p$ and $\Gamma$. Namely, $N_{T}=N_{f,inp}+N_{p,inp}+N_{f,bp}+N_{p,bp}+N_{f,\Gamma}+N_{p,\Gamma}$. For convenience, we set \textbf{K}=KI, where K is a positive constant. The global domain $\Omega$ consists of two subdomains with free fluid flow region $\Omega_f$ and porous medium region $\Omega_p$, where the interface between the conduit and matrix region is $\Gamma$ and $\Omega=\Omega_f \cup \Omega_p \cup \Gamma.$

\subsection{Example 1: The Stokes-Darcy coupled problem with a fixed physical coefficients}
In this example, we consider the Stokes-Darcy coupled problem with a linear interface (see Fig.3). The exact solution are\\

 $\textbf{u}_f=(u_1,u_2)^T=
 \left(
 \begin{matrix}
 	[x^2(y-1)^2+y] \\
 	[-\frac{2}{3}x(y-1)^3]+[2-\pi sin(\pi x)] \\
 \end{matrix}
 \right),$\\
 
 $p=[2-\pi sin(\pi x)]sin(\frac{1}{2}\pi y),$\\ 
 
  $\phi=[2-\pi sin(\pi x)][1-y-cos(\pi y)],$\\
  
  $\textbf{u}_p=-K\nabla \phi.$\\
  
The fixed kinematic viscosity $v=1$ of fluid and the element of hydraulic conductivity tensor $K=1$.\\
Case 1: $\phi=0$,  on $\Gamma_p$, $\Omega_f=[0,1]\times [1,1.25]$, $\Omega_p=[0,1]\times [0.25,1]$, $\Gamma=(0,1)\times[1],$ $\frac{\alpha v\sqrt{d}}{\sqrt{trace(\Pi)}}=1$, $g=1$.\\
Case 2: $\textbf{u}_p\cdot \textbf{n}_p=0,$ on $\Gamma_p,$ $\Omega_f=[0,1]\times [1,2]$, $\Omega_p=[0,1]\times [0,1]$, $\Gamma=(0,1)\times[1],$ $\alpha=1,g=1.$
 \begin{figure}
	\centering
	\includegraphics[scale=.6]{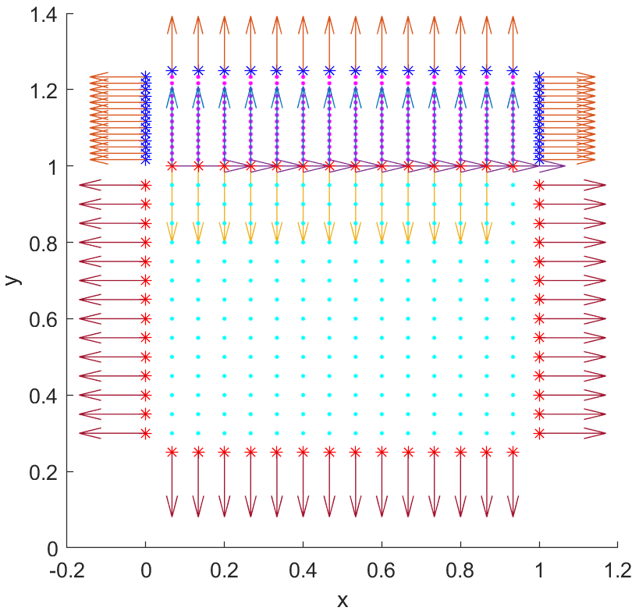}
		\raisebox{0.5\height}{\includegraphics[scale=.6]{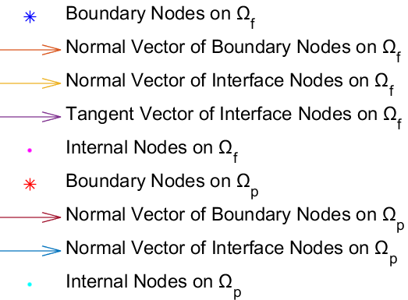}}
	\caption{ The point collocation for Example 1 Case 1.}
\end{figure}
\begin{figure}
	\centering
	\includegraphics[scale=.4]{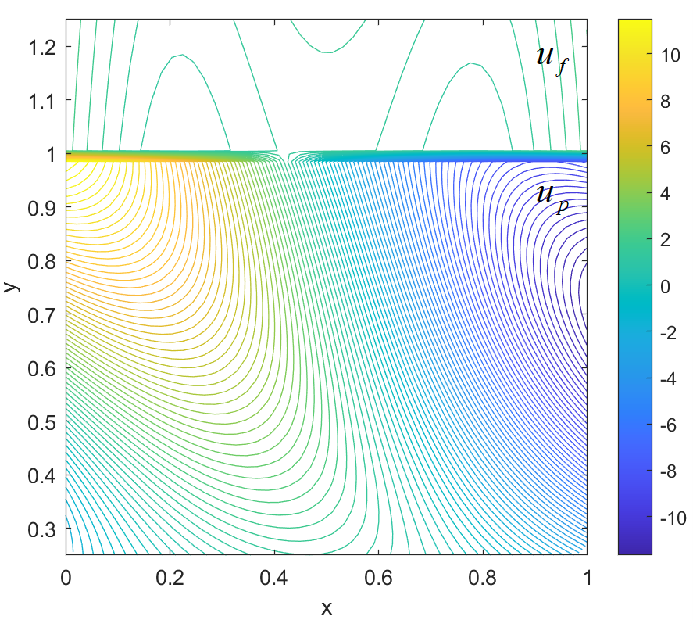}
	\hspace{1cm}
	\includegraphics[scale=.4]{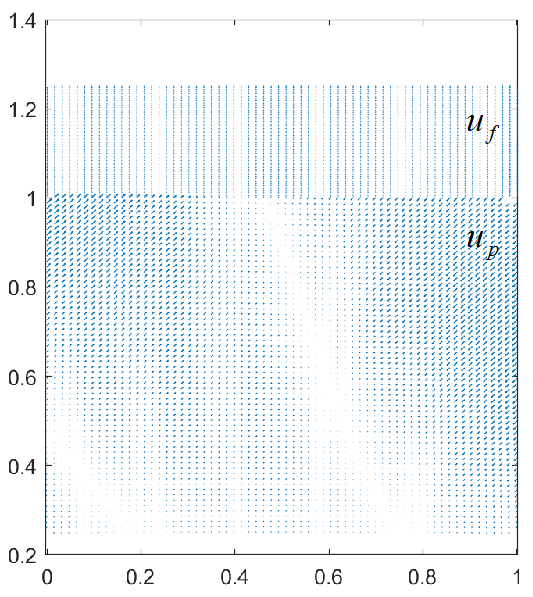}

	\caption{ The contour (left) and the vector (right) of the numerical solution for Example 1 Case 1.}
\end{figure}

	\begin{table*}	
	\scriptsize
	\caption{ The comparison between the GFDM and the LGI FEM [23] for Example 1 Case 1}
	\begin{tabular}{ccccccccccc}
		\hline
		\multirow{1}{*}{$ $}&\multicolumn{1}{c}{$N_x$} &\multicolumn{1}{c}{$h$}& \multicolumn{2}{c}{$u_f$} & \multicolumn{1}{c}{$p$}& \multicolumn{2}{c}{$\phi$} &\multicolumn{1}{c}{$CPU(s)$}\\
		\hline
		
		$ $&$ $&$ $& $L_{2,relative}$      &  $H^1_{relative}$  
		& $L_{2,relative}$   
		&  $L_{2,relative}$  
		&  $H^1_{relative}$ \\ 
		\hline		
		
		2nd GFDM&$8$ & $ $&$8.00\times10^{-2}$ & $2.77\times10^{-1}$ & $2.23\times10^{-1}$ &  $7.36 \times10^{-2}$ & $1.11\times10^{-1}$& $0.32$\\
		
		$m=40$ &$16$  &$ $& $1.58\times10^{-3}$ & $2.03 \times10^{-2}$ & $8.32 \times10^{-3}$ &  $1.60\times10^{-2}$ & $1.63\times10^{-2}$&$0.20$\\
		
		$ $ &$32$  &$ $& $4.30 \times10^{-4}$ & $4.61\times10^{-3}$ & $2.30\times10^{-3}$ &  $4.92\times10^{-3}$ &  $2.78\times10^{-3}$&$0.32$\\
		
		$ $ &$64$  &$ $& $1.10\times10^{-4}$ & $1.09\times10^{-3}$ & $5.92 \times10^{-4}$ &  $1.31\times10^{-3}$ &  $5.61\times10^{-4}$&$0.87$\\
				
		$ $ &$128$  &$ $& $2.77\times10^{-5}$ & $2.65\times10^{-4}$ & $1.48\times10^{-4}$ &  $3.35\times10^{-4}$ &  $1.27\times10^{-4}$&$4.22$\\
		
		$ $ &$Order$  &$ $& $2.71$ & $3.08$ & $2.64$ &  $1.94$ &  $2.44$\\
		\hline
		LGI FEM[23]&$ $ & $\frac{1}{8}$&$1.55\times10^{-2}$ & $3.97\times10^{-1}$ & $1.15\times10^{-1}$ &  $1.14\times10^{-2}$ & $5.50\times10^{-1}$\\
		
		$ $ &$ $  &$\frac{1}{16}$& $3.92\times10^{-3}$ & $1.98\times10^{-1}$ & $3.67\times10^{-2}$ &  $2.88\times10^{-3}$ & $2.75\times10^{-1}$ \\
		
		$ $ &$ $  &$\frac{1}{32}$& $9.87\times10^{-4}$ & $9.91\times10^{-2}$ & $1.16\times10^{-2}$ &  $7.20\times10^{-4}$ & $1.38\times10^{-1}$ \\
		
		$ $ &$ $  &$\frac{1}{64}$& $2.47\times10^{-4}$ & $4.95\times10^{-2}$ & $3.71\times10^{-3}$ &  $1.80\times10^{-4}$ & $6.88\times10^{-2}$ \\
		
		$ $ &$ $  &$\frac{1}{128}$& $6.20\times10^{-5}$ & $2.48\times10^{-2}$ & $1.22\times10^{-3}$ &  $4.50\times10^{-5}$ & $3.44\times10^{-2}$\\
		
		$ $ &$$  &$Order$& $1.99$ & $1.00$ & $1.64$ &  $2.00$ & $1.00$ & \\
		\hline
		
	\end{tabular}
\end{table*}

\begin{table*}	
	\scriptsize
	\caption{ The $L_{2,relative}$ errors when the different order GFDMs are adopted for Example 1 Case 1}
	\begin{tabular}{cccccccccccc}
		\hline
		\multirow{1}{*}{$ $}&\multicolumn{1}{c}{$N_x$} & \multicolumn{2}{c}{$u_f$} &\multicolumn{2}{c}{$p$}& \multicolumn{2}{c}{$\phi$}& \multicolumn{2}{c}{$u_p$} & \multicolumn{1}{c}{$CPU(s)$}\\
		\hline
		
		$ $&$ $& $L_{2,relative}$      &  $Order$  
		& $L_{2,relative}$   &  $Order$     &  $L_{2,relative}$  &  $Order$  
		&  $L_{2,relative}$  &  $Order$   \\ 
		\hline		
		
		2nd order GFDM &$16$ & $1.58\times10^{-3}$ & $-$ & $8.32 \times10^{-3}$ &  $-$ & $1.60\times10^{-2}$ & $-$& $1.28\times10^{-2}$ & $-$&$0.13$\\
		
		$m=40$ &$32$ & $4.30\times10^{-4}$ & $1.88 $ & $2.30\times10^{-3}$ &  $1.85$ & $4.92\times10^{-3}$ & $1.70$& $2.44 \times10^{-3}$ & $2.39$&$0.27$\\
		
		$ $ &$64$ & $1.10\times10^{-4}$ & $1.97$ & $5.92\times10^{-4}$ &  $1.96 $ & $1.31\times10^{-3}$ & $1.91 $& $5.33\times10^{-4}$ & $ 2.19$&$0.89$\\
		
		$ $ &$128$ & $2.77\times10^{-5}$ & $1.99$ & $1.48\times10^{-4}$ &  $2.00$ & $3.35\times10^{-4}$ & $1.97$& $1.25\times10^{-4}$ & $2.09$&$4.14$\\
		
		\hline

		4th order GFDM &$32$ & $2.74\times10^{-5}$ & $-$ & $1.32 \times10^{-4}$ &  $-$ & $2.87\times10^{-4}$ & $-$& $2.48\times10^{-4}$ &$-$& $1.23$\\
		
		$m=210$ &$64$ & $1.34\times10^{-6}$ & $4.35$ & $8.50\times10^{-6}$ &  $ 3.96$ & $2.05 \times10^{-5}$ & $3.81 $& $1.07 \times10^{-5}$ & $4.53 $&$66.2$\\
		
		$ $ &$128$ & $8.86\times10^{-8}$ & $3.92$ & $ 5.75\times10^{-7}$ &  $3.89$ & $ 1.40\times10^{-6}$ & $3.88$& $5.06\times10^{-7}$ & $4.40$&$47.1$\\
		
		\hline
		6th order GFDM&$32$ &$1.63\times10^{-6}$ & $-$ & $1.11\times10^{-5}$ &  $-$ & $2.53\times10^{-5}$ & $-$& $1.79\times10^{-5}$ & $-$&$5.06$\\

         $m=450$ &$64$ & $3.05\times10^{-8}$ & $5.74$ & $2.06\times10^{-7}$ &  $5.74$ & $4.16\times10^{-7}$ & $5.92 $& $2.33\times10^{-7}$ & $6.26$&$88.6$\\
	
       $ $ &$128$ & $5.40\times10^{-10}$ & $5.82$ & $3.70\times10^{-9}$ &  $5.80 $ & $ 7.24\times10^{-9}$ & $5.85$& $3.16\times10^{-9}$ & $6.21$&$183.2$\\
	
	\hline
	\end{tabular}
\end{table*}

\begin{table*}	
	\scriptsize
	\caption{ The $H^1_{relative}$ errors when the different order GFDMs are adopted for Example 1 Case 1}
	\begin{tabular}{cccccccccccc}
		\hline
		\multirow{1}{*}{$ $}&\multicolumn{1}{c}{$N_x$} & \multicolumn{2}{c}{$u_f$} &\multicolumn{2}{c}{$p$}& \multicolumn{2}{c}{$\phi$}& \multicolumn{2}{c}{$u_p$} & \multicolumn{1}{c}{$CPU(s)$}\\
		\hline
		
		$ $&$ $& $H^1_{relative}$      &  $Order$  
		& $H^1_{relative}$    &  $Order$     &  $H^1_{relative}$    &  $Order$  
		&  $H^1_{relative}$    &  $Order$   \\ 
		\hline

	2nd order GFDM&$8$ &$2.77\times10^{-1}$ & $-$ & $3.43\times10^{-1}$ &  $-$ & $1.11\times10^{-1}$ & $-$& $3.86\times10^{-1}$ & $-$&$ 0.17$\\
	
	$m=40$ &$16$ & $2.03\times10^{-2}$ & $3.78  $ & $1.84 \times10^{-2}$ &  $4.22  $ & $1.63\times10^{-2}$ & $2.77 $& $1.14 \times10^{-1}$ & $1.76  $&$0.14$\\
	
	$ $ &$32$ & $4.61\times10^{-3}$ & $ 2.13 $ & $4.88 \times10^{-3}$ &  $1.92$ & $2.78\times10^{-3}$ & $ 2.55 $& $ 3.70\times10^{-2}$ & $1.63$&$0.28$\\
	
	$ $ &$64$ & $1.09\times10^{-3}$ & $2.08$ & $1.30\times10^{-3}$ &  $1.91$ & $ 5.61\times10^{-4}$ & $2.31 $& $1.25\times10^{-2}$ & $1.57$&$0.91$\\
	
	$ $ &$128$ & $2.65\times10^{-4}$ & $2.04$ & $ 3.52\times10^{-4}$ &  $1.89$ & $1.27\times10^{-4}$ & $2.14$& $ 4.31\times10^{-3}$ & $1.53$&$4.22$\\
	
	\hline
	
	4th order GFDM&$32$ &$1.80\times10^{-4}$ & $-$ & $3.27\times10^{-4}$ &  $-$ & $2.56\times10^{-4}$ & $-$& $2.04\times10^{-3}$ & $-$&$ 1.27$\\
	
	$m=210$ &$64$ & $ 6.23\times10^{-6}$ & $4.85$ & $ 1.96 \times10^{-5}$ &  $4.06$ & $1.08 \times10^{-5}$ & $4.56$& $1.75\times10^{-4}$ & $3.55$&$65.5$\\
	
	$ $ &$128$ & $2.68\times10^{-7}$ & $4.54$ & $1.20 \times10^{-6}$ &  $ 4.03$ & $5.09\times10^{-7}$ & $ 4.41$& $ 1.49\times10^{-5}$ & $3.56$&$45.6$\\

	\hline
	6th order GFDM&$32$ &$3.79\times10^{-6}$ & $-$ & $8.18\times10^{-6}$ &  $-$ & $1.80\times10^{-5}$ & $-$& $1.41\times10^{-4}$ & $-$&$5.17$\\
	
	$m=450$ &$64$ & $7.27\times10^{-8}$ & $5.70 $ & $1.12\times10^{-7}$ &  $6.19$ & $2.37\times10^{-7}$ & $6.25$& $3.60 \times10^{-6}$ & $5.29$&$79.8$\\
	
	$ $ &$128$ & $1.23\times10^{-9}$ & $5.88$ & $4.31\times10^{-9}$ &  $4.70$ & $ 3.18\times10^{-9}$ & $6.22$& $7.87\times10^{-8}$ & $5.52$&$171.3$\\
	
	\hline
	\end{tabular}
\end{table*}
\begin{figure}
	\centering
	\includegraphics[scale=.6]{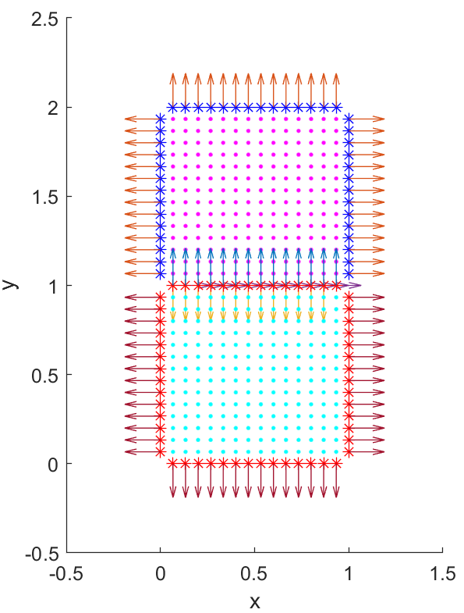}
	\raisebox{0.5\height}{\includegraphics[scale=.6]{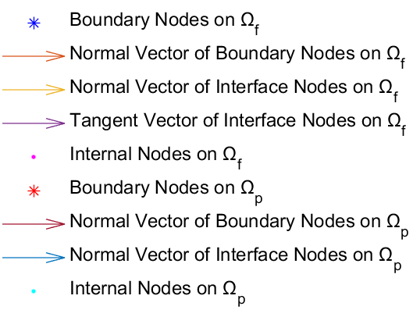}}
	\caption{ The point collocation for Example 1 Case 2.}
\end{figure}
\begin{figure}
	\centering
	\includegraphics[scale=.4]{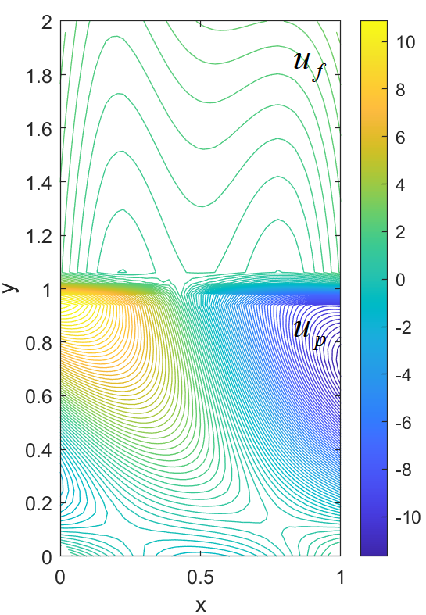}
	\hspace{1cm}
	\includegraphics[scale=.4]{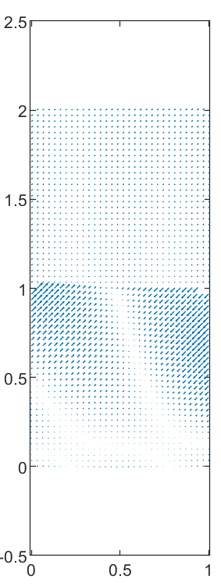}
	
	\caption{ The contour (left) and the vector (right) of the numerical solution for Example 1 Case 2.}
\end{figure}

	\begin{table*}	
	\scriptsize
	\caption{ The comparison between the GFDM and the FEDDM [20] for Example 1 Case 2}
	\begin{tabular}{ccccccccccc}
		\hline
		\multirow{1}{*}{$ $}&\multicolumn{1}{c}{$N_x$} &\multicolumn{1}{c}{$h$}& \multicolumn{2}{c}{$u_f$} & \multicolumn{1}{c}{$p$}& \multicolumn{1}{c}{$\phi$} & \multicolumn{2}{c}{$u_p$}&\multicolumn{1}{c}{$CPU(s)$}\\
		\hline
		
		$ $&$ $&$ $& $L_{2,relative}$      &  $H^1_{relative}$  
		& $L_{2,relative}$      &  $L_{2,relative}$  
		&  $L_{2,relative}$  
		&  $H^1_{relative}$ \\ 
		\hline		
		
		4th GFDM&$8$ & $ $&$1.09\times10^{-1}$ & $1.29\times10^{-1}$ & $3.56\times10^{-1}$ &  $2.03\times10^{-0}$ & $2.76 \times10^{-2}$ & $5.34\times10^{-2}$&$0.15$\\
		
		$m=40$ &$16$  &$ $& $4.45\times10^{-3}$ & $4.94\times10^{-3}$ & $1.39\times10^{-2}$ &  $7.88\times10^{-2}$ & $1.20\times10^{-3}$ & $ 5.17\times10^{-3}$&$0.19$\\
		
		$ $ &$32$  &$ $& $2.26\times10^{-4}$ & $2.35\times10^{-4}$ & $6.97 \times10^{-4}$ &  $3.92\times10^{-3}$ & $ 5.26\times10^{-5}$ & $3.82\times10^{-4}$&$0.41$\\
		
		$ $ &$64$  &$ $& $1.27\times10^{-5}$ & $1.26\times10^{-5}$ & $3.90\times10^{-5}$ &  $2.18\times10^{-4}$ & $2.80\times10^{-6}$ & $2.97\times10^{-5}$&$1.44$\\
		
		\hline
		FEDDM[20]&$ $ & $\frac{1}{8}$&$3.25\times10^{-2}$ & $7.33\times10^{-1}$ & $2.02\times10^{-1}$ &  $1.42\times10^{-1}$ & $1.34\times10^{-1}$ & $ 2.80\times10^{-0}$\\
		
		$ $ &$ $  &$\frac{1}{16}$& $8.14\times10^{-3}$ & $3.65\times10^{-1}$ & $6.64\times10^{-2}$ &  $7.04\times10^{-2}$ & $3.44\times10^{-2}$ & $ 1.41\times10^{-0}$&\\
		
		$ $ &$ $  &$\frac{1}{32}$& $2.04\times10^{-3}$ & $1.82\times10^{-1}$ & $2.27\times10^{-2}$ &  $3.51\times10^{-2}$ & $8.68\times10^{-3}$ & $ 7.04\times10^{-1}$&\\
		
		$ $ &$ $  &$\frac{1}{64}$& $5.10\times10^{-4}$ & $9.11\times10^{-2}$ & $7.90\times10^{-3}$ &  $1.75\times10^{-2}$ & $2.18\times10^{-3}$ & $ 3.52\times10^{-1}$&\\
		
		\hline
		
	\end{tabular}
\end{table*}

	\begin{table*}	
	\scriptsize
	\caption{ The $L_{2,relative}$ errors when the different order GFDMs are adopted for Example 1 Case 2}
	\begin{tabular}{cccccccccccc}
		\hline
		\multirow{1}{*}{$ $}&\multicolumn{1}{c}{$N_x$} & \multicolumn{2}{c}{$u_f$} &\multicolumn{2}{c}{$p$}& \multicolumn{2}{c}{$\phi$}& \multicolumn{2}{c}{$u_p$} & \multicolumn{1}{c}{$CPU(s)$}\\
		\hline
		
		$ $&$ $& $L_{2,relative}$      &  $Order$  
		& $L_{2,relative}$   &  $Order$     &  $L_{2,relative}$  &  $Order$  
		&  $L_{2,relative}$  &  $Order$   \\ 
		\hline		
		
		4th-order GFDM&$16$ &$4.45\times10^{-3}$ & $-$ & $1.39\times10^{-2}$ &  $-$ & $7.88\times10^{-2}$ & $-$& $1.20\times10^{-3}$ & $-$&$ 0.18$\\
		
		$m=40$ &$32$ & $2.26\times10^{-4}$ & $4.30 $ & $6.97 \times10^{-4}$ &  $4.32$ & $3.92\times10^{-3}$ & $4.33$& $5.26\times10^{-5}$ & $4.51 $&$0.40$\\
		
		$ $ &$64$ & $1.27\times10^{-5}$ & $4.15$ & $3.90\times10^{-5}$ &  $4.16$ & $2.18\times10^{-4}$ & $4.17 $& $2.80\times10^{-6}$ & $4.23 $&$1.56$\\
		
		$ $ &$128$ & $7.55\times10^{-7}$ & $4.08$ & $ 2.31\times10^{-6}$ &  $4.08$ & $1.28\times10^{-5}$ & $4.09$& $1.68\times10^{-7}$ & $4.06$&$14.9$\\
		
		\hline
        6th-order GFDM&$16$ &$2.55\times10^{-4}$ & $-$ & $8.02\times10^{-4}$ &  $-$ & $5.33\times10^{-3}$ & $-$& $2.08\times10^{-4}$ & $-$&$0.61$\\

        $m=140$ &$32$ & $2.15\times10^{-6}$ & $6.89 $ & $6.76\times10^{-6}$ &  $6.89$ & $4.33\times10^{-5}$ & $6.95 $& $1.86\times10^{-6}$ & $6.81$&$4.26$\\

        $ $ &$64$ & $2.44\times10^{-8}$ & $6.47 $ & $ 7.61 \times10^{-8}$ &  $6.47 $ & $4.16\times10^{-7}$ & $6.70 $& $2.61\times10^{-8}$ & $6.15 $&$23.0$\\

        $ $ &$128$ & $1.94\times10^{-10}$ & $6.97$ & $7.63\times10^{-10}$ &  $6.64$ & $3.47\times10^{-9}$ & $6.91$& $5.17\times10^{-10}$ & $5.66$&$102.0$\\

\hline

	\end{tabular}
\end{table*}

	\begin{table*}	
	\scriptsize
	\caption{ The $H^1_{relative}$ errors when the different order GFDMs are adopted for Example 1 Case 2}
	\begin{tabular}{cccccccccccc}
		\hline
		\multirow{1}{*}{$ $}&\multicolumn{1}{c}{$N_x$} & \multicolumn{2}{c}{$u_f$} &\multicolumn{2}{c}{$p$}& \multicolumn{2}{c}{$\phi$}& \multicolumn{2}{c}{$u_p$} & \multicolumn{1}{c}{$CPU(s)$}\\
		\hline
		
		$ $&$ $& $H^1_{relative}$      &  $Order$  
		& $H^1_{relative}$    &  $Order$     &  $H^1_{relative}$    &  $Order$  
		&  $H^1_{relative}$    &  $Order$   \\ 
		\hline		
		
		4th-order GFDM&$16$ &$4.91\times10^{-3}$ & $-$ & $9.60\times10^{-3}$ &  $-$ & $1.41 \times10^{-3}$ & $-$& $6.67\times10^{-3}$ & $-$&$ 0.21$\\
		
		$m=50$ &$32$ & $2.25\times10^{-4}$ & $4.45 $ & $5.69\times10^{-4}$ &  $4.08$ & $6.31 \times10^{-5}$ & $4.49$& $4.57\times10^{-4}$ & $3.87$&$0.44$\\
		
		$ $ &$64$ & $1.19\times10^{-5}$ & $4.24$ & $3.78\times10^{-5}$ &  $3.91 $ & $ 3.12\times10^{-6}$ & $4.33$& $3.07\times10^{-5}$ & $3.90$&$1.84$\\
		
		$ $ &$128$ & $6.84\times10^{-7}$ & $4.13$ & $ 2.65\times10^{-6}$ &  $3.84$ & $ 1.73 \times10^{-7}$ & $4.18$& $2.21\times10^{-6}$ & $ 3.80$&$11.1$\\
		
		\hline
		6th-order GFDM&$16$ &$2.88\times10^{-4}$ & $-$ & $6.19\times10^{-4}$ &  $-$ & $2.07\times10^{-4}$ & $-$& $8.55\times10^{-4}$ & $-$&$ 0.75$\\
		
		$m=140$ &$32$ & $2.69\times10^{-6}$ & $6.75$ & $7.88\times10^{-6}$ &  $6.30$ & $ 1.97\times10^{-6}$ & $6.71$& $2.27\times10^{-5}$ & $5.23$&$2.79$\\
		
		$ $ &$64$ & $3.31\times10^{-8}$ & $6.34$ & $1.21 \times10^{-7}$ &  $6.03 $ & $2.60\times10^{-8}$ & $6.25$& $3.78\times10^{-7}$ & $5.91  $&$12.1$\\
		
		$ $ &$128$ & $3.29\times10^{-10}$ & $ 6.65$ & $2.46\times10^{-9}$ &  $ 5.62$ & $5.13\times10^{-10}$ & $5.66$& $5.76\times10^{-9}$ & $6.04$&$82.0$\\
		
		\hline
		
	\end{tabular}
\end{table*}

Fig.3 and Fig. 5 present the point collocations, in which we can see that the point collocation in each subdomain and the distributions of the normal vector and tangent vector for both Case 1 and Case 2. The contour (left) and the vector (right) of the numerical solution are shown in Fig.4 and Fig. 6, we can see that the distributions of the numerical solutions. The comparison between the GFDM and the LGI FEM[23]  is given in Table. 1. The comparison between the GFDM and the FEDDM[20]  is given in Table. 4. From these tables, we can see that our method is more accurate than the LGI FEM[23] and FEDDM [20]highly efficient because the CPU time is very small. Furthermore, our errors are all relative errors which can illustrate the numerical solution are accurate and stable. Table.2 and Table.3 show the $L_2, $$H^1_{relative}$ errors when the 2nd GFDM, 4th GFDM and 6th GFDM are adopted. Table.5 and Table.6 show the $L_2, $$H^1_{relative}$ errors when the 4th GFDM and 6th GFDM are adopted. We can see that the 2nth order covergence can be obtained, when the 2nth order derivatives of Taylar series are truncated, (n=1,2,3). 
\subsection{Example 2: The Stokes-Darcy coupled problem with a varying physical coefficients}
In this example, we consider the Stokes-Darcy coupled problem with a linear interface (see Fig.5). The exact solutions are\\

$\textbf{u}_f=(u_1,u_2)^T=
\left(
\begin{matrix}
	[y^2-2y+1] \\
	[x^2-x] \\
\end{matrix}
\right)$,\\

$p=2v(x+y-1)+\frac{1}{3K},$\\ 

$\phi=\frac{1}{K}[x(1-x)(y-1)+\frac{1}{3}y^3-y^2+y]+2vx,$\\

$\textbf{u}_p=-K\nabla \phi.$\\

The biggest feature of these exact solutions is that the $p$ and $\phi$ are all related to the kinematic viscosity $v$ of the fluid and the element of the hydraulic conductivity tensor $K$.
In order to examine the proposed meshless GFDM is also applicable with the realistic physical parameters is of practical interest, which is proposed by Ref[20]. We consider the Stokes-Darcy coupled problem with a varying kinematic viscosity $v$ of the fluid and the element of the hydraulic conductivity tensor $K$. In this example, the computational subdomains are $\Omega_f=[0,1]\times [1,2]$, $\Omega_p=[0,1]\times [0,1]$ and $\Gamma=(0,1)\times[1]$ and we adopt the initial coefficients $v=1,$ $K=1,$ $g=1$ and $\alpha=1$.  

\begin{figure}
	\centering
	\includegraphics[scale=.4]{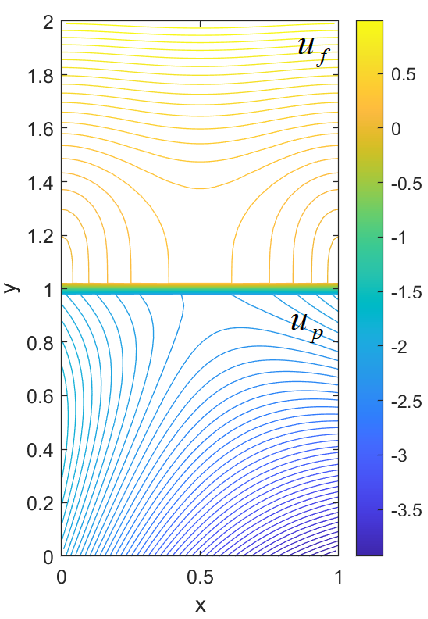}
	\hspace{1cm}
	\includegraphics[scale=.4]{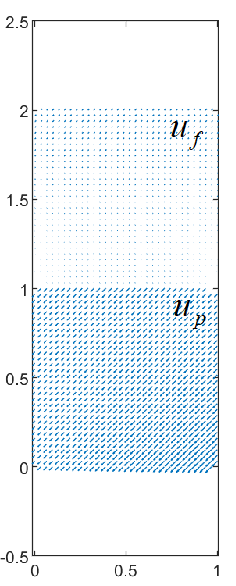}
	\caption{ The contour(left) and the vector (right) of numerical solution for Example 2.}
\end{figure}
\begin{figure}
	\centering
	\includegraphics[scale=.4]{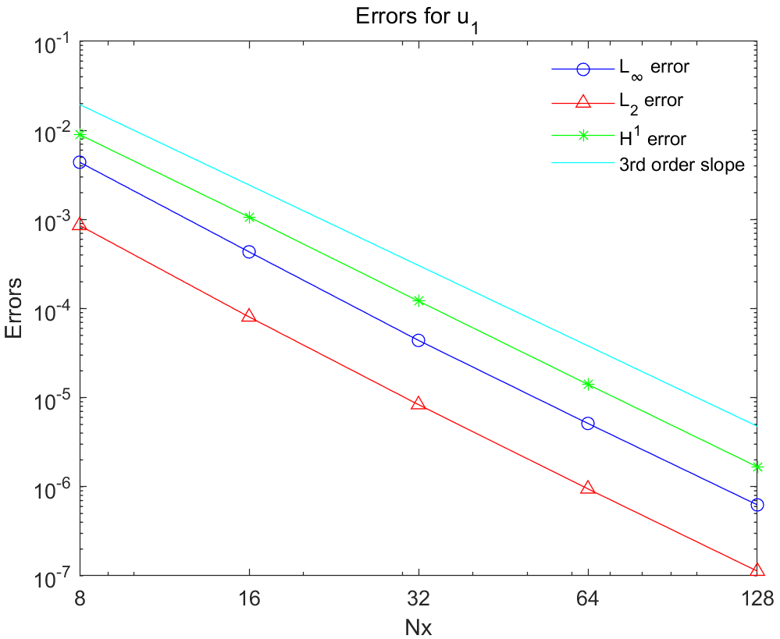}
	\includegraphics[scale=.4]{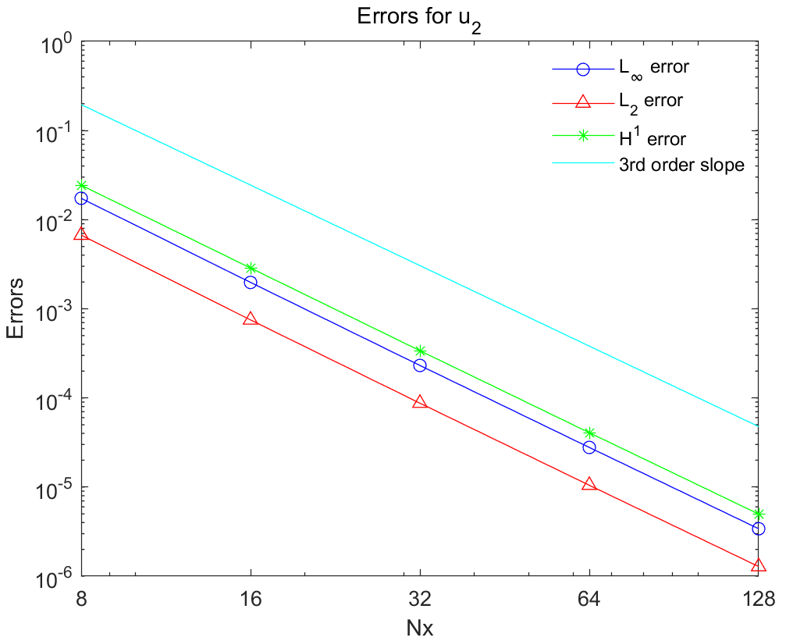}
	\includegraphics[scale=.4]{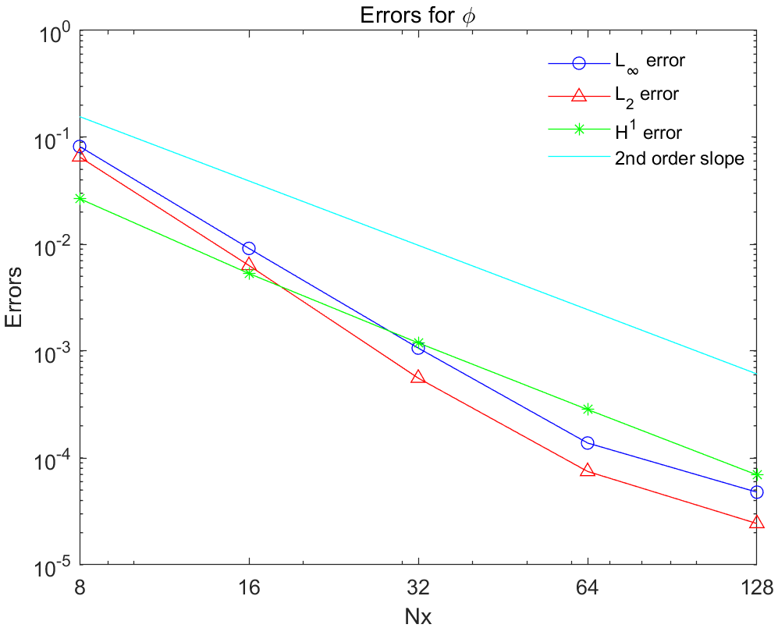}
	\caption{ $L_{\infty}$,$L_2$ and $H^1$ errors of the component functions for Example 2.}
\end{figure}
\begin{figure}
	\centering
	\includegraphics[scale=.4]{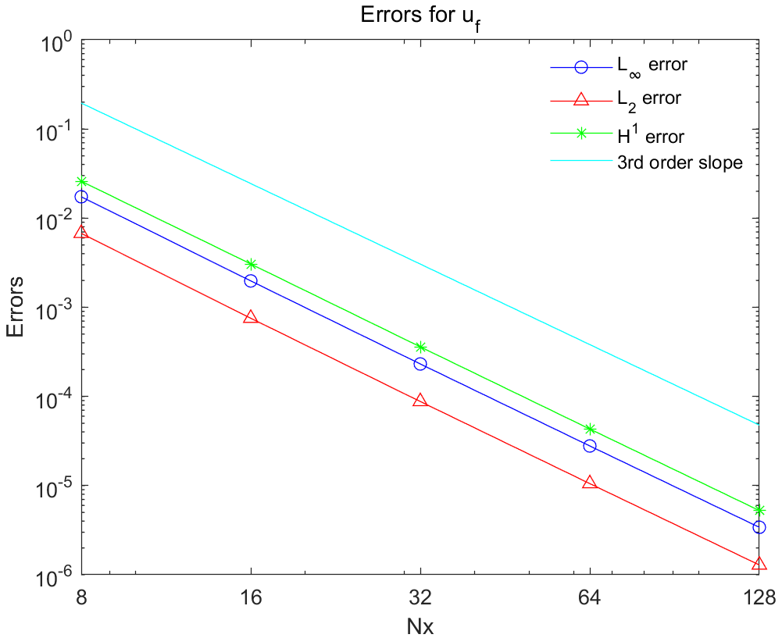}
	\includegraphics[scale=.4]{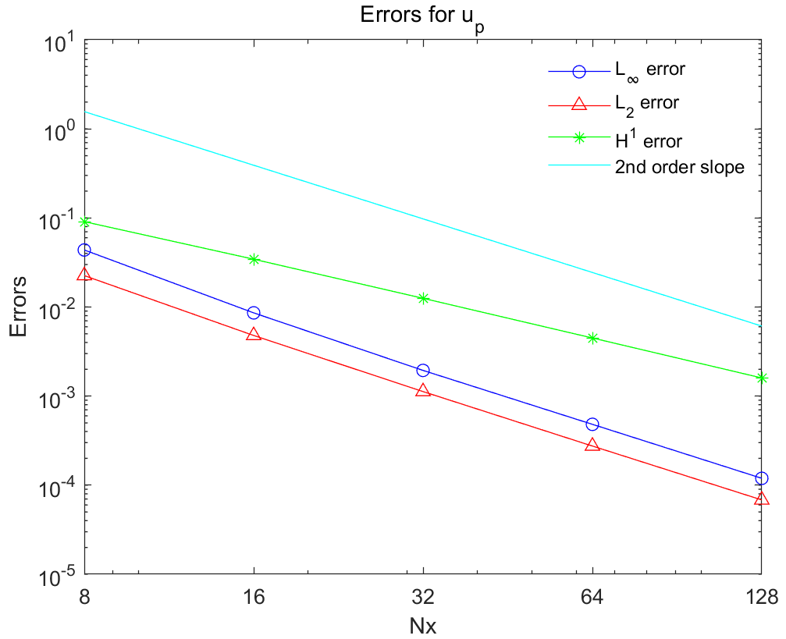}
	\includegraphics[scale=.4]{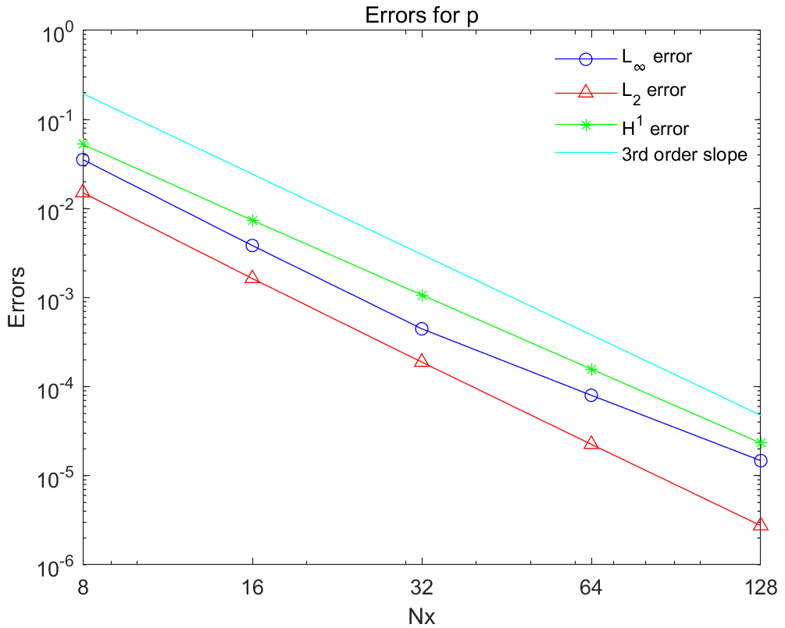}
	\caption{ $L_{\infty}$,$L_2$ and $H^1$ errors of fluid velocity $\textbf{u}_f$ (left), Darcy velocity $\textbf{u}_p$ (middle) and pressure $p$ (right) for Example 2.}
\end{figure}
	\begin{table*}	
	\scriptsize
	\caption{ The $L_2$ and $H^1$ relative errors and CPU (s) of the GFDM for Example 2}
	\begin{tabular}{ccccccccccc}
		\hline
		\multirow{1}{*}{$ $}&\multicolumn{1}{c}{$N_x$}& \multicolumn{2}{c}{$u_f$} & \multicolumn{1}{c}{$p$}& \multicolumn{1}{c}{$\phi$} & \multicolumn{1}{c}{$u_p$}&\multicolumn{1}{c}{$CPU(s)$}\\
		\hline
		
		$ $&$ $& $L_{2,relative}$      &  $H^1_{relative}$  
		& $L_{2,relative}$      &  $L_{2,relative}$  
		&  $L_{2,relative}$   \\ 
		\hline		
		
		2nd GFDM&$8$ &$3.79\times10^{-2}$ & $3.98\times10^{-2}$ & $6.03\times10^{-3}$ &  $4.97\times10^{-2}$ & $8.81\times10^{-3}$&$0.15$ \\
		
		$ $ &$16$ & $4.20\times10^{-3}$ & $4.83\times10^{-3}$ & $6.59\times10^{-4}$ &  $4.79\times10^{-3}$ & $1.88\times10^{-3}$&$0.10$ \\
		
		$ $ &$32$ & $4.85\times10^{-4}$ & $5.81 \times10^{-4}$ & $7.61\times10^{-5}$ &  $4.24 \times10^{-4}$ & $4.42\times10^{-4}$&$0.17$ \\
		
		$ $ &$64$& $5.79\times10^{-5}$ & $7.06\times10^{-5}$ & $9.11\times10^{-6}$ &  $5.69 \times10^{-5}$ & $1.08\times10^{-4}$&$1.14$ \\
		
		$ $ &$128$& $ 7.05\times10^{-6}$ & $8.67\times10^{-6}$ & $1.11\times10^{-6}$ &  $1.85\times10^{-5}$ & $2.69\times10^{-5}$&$6.18$ \\
		\hline
		
	\end{tabular}
\end{table*}

\begin{figure}
	\centering
	\includegraphics[scale=.4]{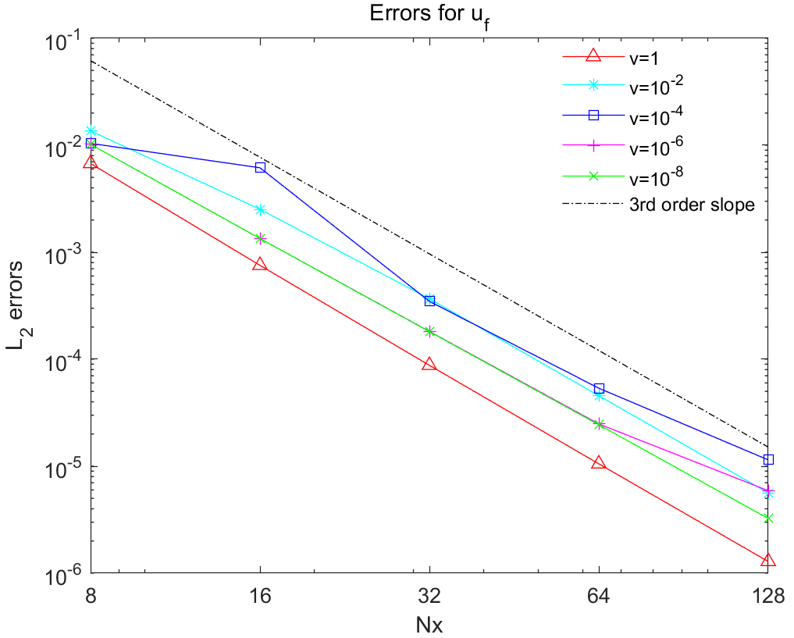}
	\includegraphics[scale=.4]{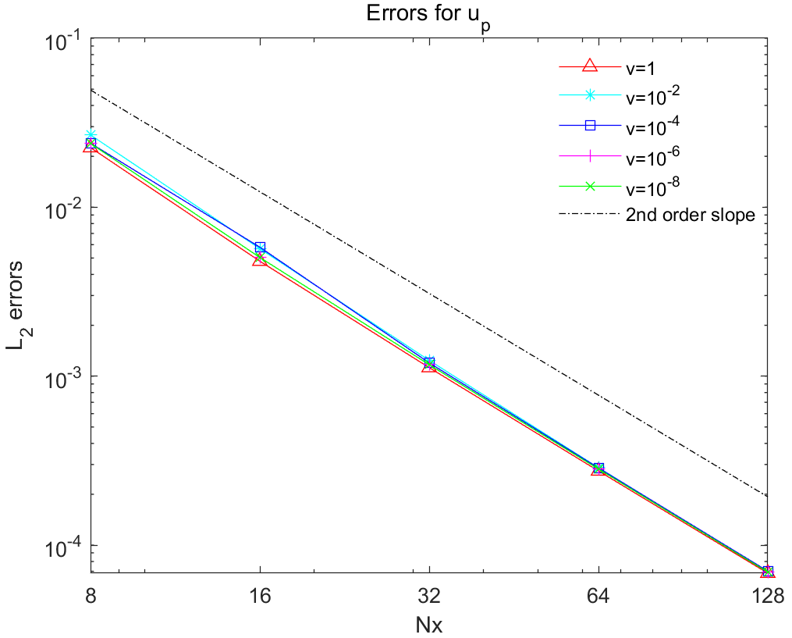}
	\includegraphics[scale=.4]{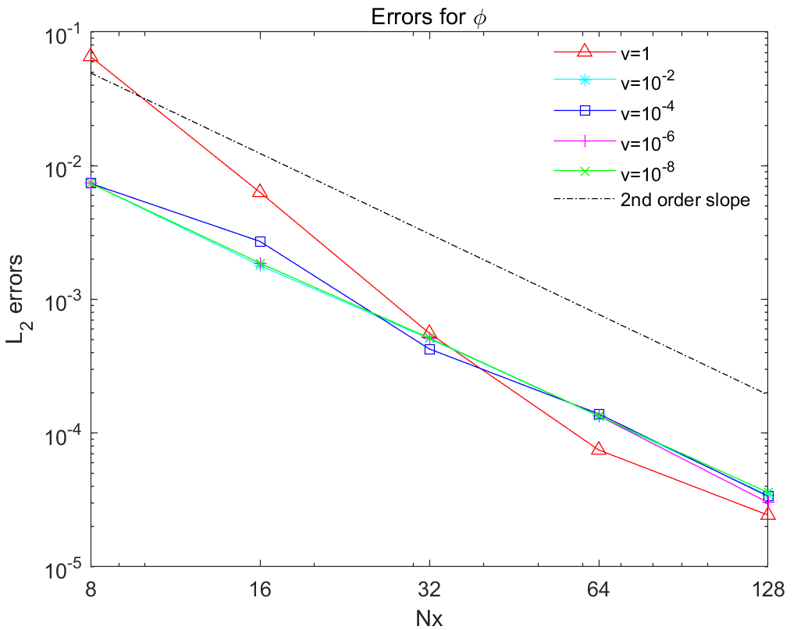}
	\caption{The $L_2$ errors for the numerical velocity of Stokes $\textbf{u}_f$ (left),  Darcy $\textbf{u}_p$ (middle) and $\phi$ (right) with varying the kinematic viscosity of fluid $v$ for Example 2.}
\end{figure}

\begin{figure}
	\centering
	\includegraphics[scale=.4]{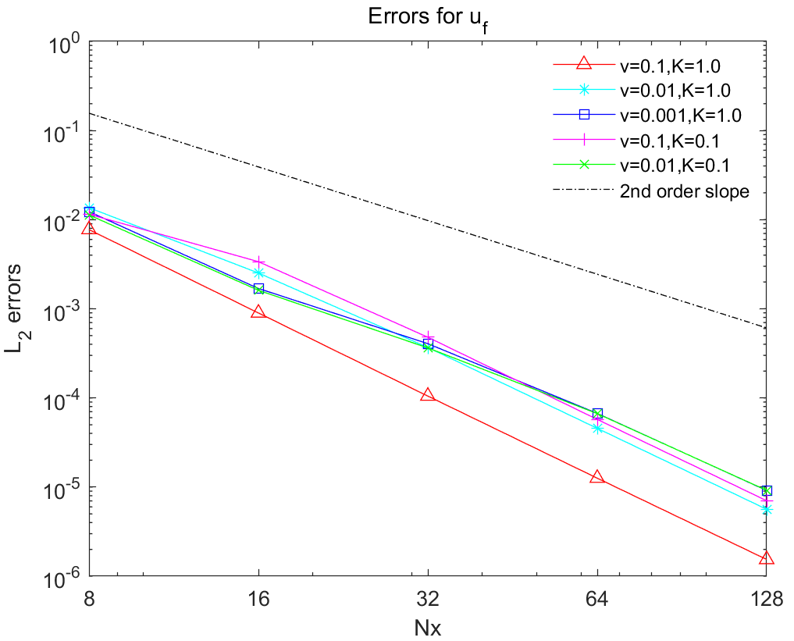}
	\includegraphics[scale=.4]{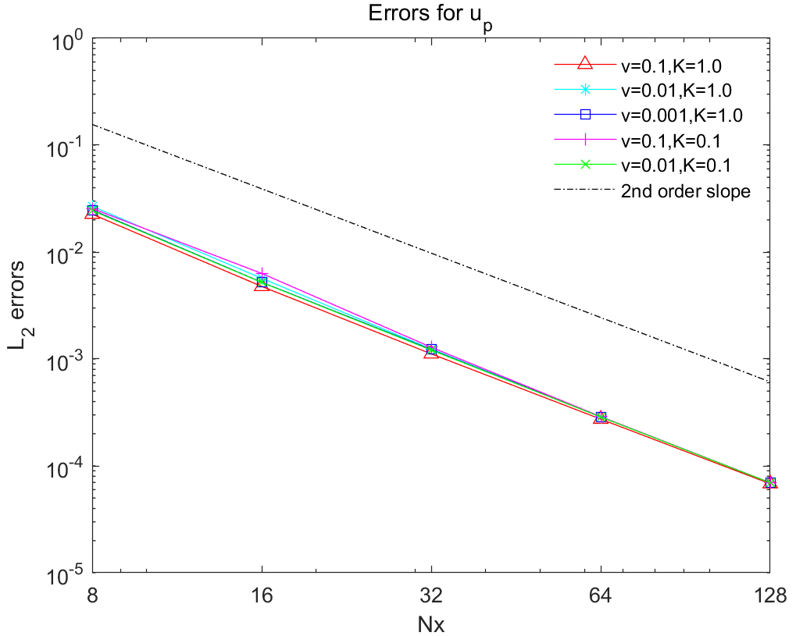}
	\includegraphics[scale=.4]{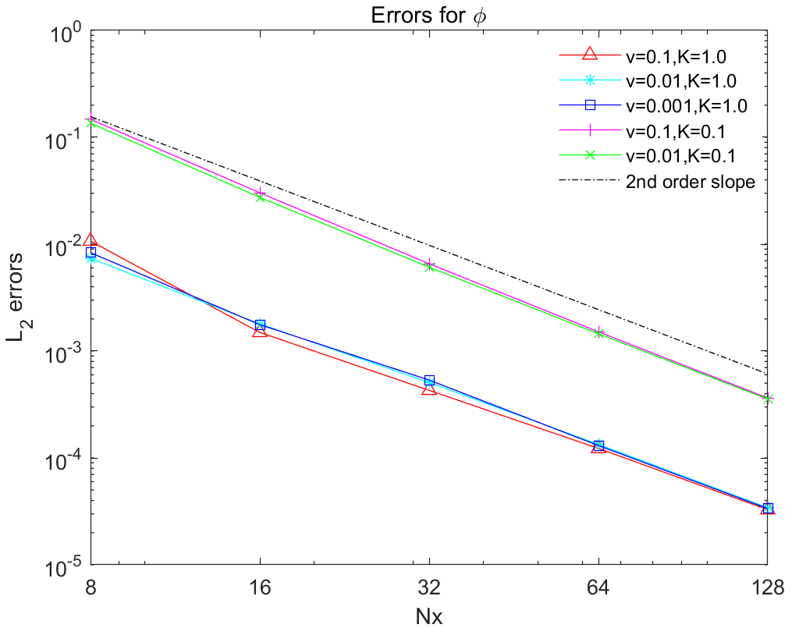}
	\caption{The $L_2$ errors for the numerical velocity of Stokes $\textbf{u}_f$ (left), Darcy $\textbf{u}_p$ (middle) and $\phi$ (right) with different parameters for Example 2.}
\end{figure}
\begin{table*}	
	\scriptsize
	\caption{ $L_{\infty}$,$L_2$ and $H^1$ errors of the GFDM with different $m$  when $N_x=32$ for Example 2}
	\begin{tabular}{cccccccccccc}
		\hline
		\multirow{1}{*}{$m$} & \multicolumn{3}{c}{$u_f$} & \multicolumn{3}{c}{$u_p$}& \multicolumn{2}{c}{$p$}&
		\multicolumn{2}{c}{$\phi$} \\
		\hline
		
		& $L_{\infty}$      &  $L_2 $   &   $H^1$
		
		& $L_{\infty}$      &  $L_2$   &   $H^{1,relative}$
		
		&  $L_2$   &   $H^1$
		&  $L_2$   &   $H^1$\\
		\hline
		$12$ & $1.60\times10^{-4}$ & $6.06\times10^{-5}$ & $2.49\times10^{-4}$ &  $1.58\times10^{-3}$ & $1.26\times10^{-3}$ & $1.68\times10^{-2}$&
		$1.31\times10^{-4}$ & $7.44\times10^{-4}$ & $3.54\times10^{-4}$&$1.31\times10^{-3}$ \\

		$14$ & $1.88\times10^{-4}$ & $7.10\times10^{-5}$ & $2.91\times10^{-4}$ &  $1.80\times10^{-3}$ & $1.24\times10^{-3}$ & $1.49\times10^{-2}$&
		$1.53\times10^{-4}$ & $8.75\times10^{-4}$ & $4.11\times10^{-4}$ &$1.30\times10^{-3}$\\
		
		$16$ & $2.25\times10^{-4}$ & $8.53\times10^{-5}$ & $3.50\times10^{-4}$ &  $2.07\times10^{-3}$ & $1.26\times10^{-3}$ & $1.41\times10^{-2}$&
		$1.84\times10^{-4}$ & $1.04\times10^{-3}$ & $5.13\times10^{-4}$&$1.32\times10^{-3}$ \\
		
		$18$ & $2.25\times10^{-4}$ & $8.53\times10^{-5}$ & $3.50\times10^{-4}$ &  $1.93\times10^{-3}$ & $1.24\times10^{-3}$ & $1.38\times10^{-2}$&
		$1.84\times10^{-4}$ & $1.04\times10^{-3}$ & $5.14\times10^{-4}$&$1.31\times10^{-3}$ \\
		
		$20$ & $2.31\times10^{-4}$ & $8.78\times10^{-5}$ & $3.60\times10^{-4}$ &  $1.94\times10^{-3}$ & $1.12\times10^{-3}$ & $1.25\times10^{-2}$&
		$1.89\times10^{-4}$ & $1.07\times10^{-3}$ & $5.58\times10^{-4}$ &$1.19\times10^{-3}$ \\
		
		$22$ & $2.39\times10^{-4}$ & $9.05\times10^{-5}$ & $3.72\times10^{-4}$ &  $1.99\times10^{-3}$ & $9.44\times10^{-4}$ & $1.18\times10^{-2}$&
		$1.95\times10^{-4}$ & $1.10\times10^{-3}$ & $6.10\times10^{-4}$ & $1.01\times10^{-3}$ \\
		\hline
	\end{tabular}
\end{table*}

The contour (left) and the vector (right) of the numerical solution are shown in Fig. 7. From Fig. 8 and Fig. 9, we can see that all results are accurate, stable and maintain 2nd order convergence. In particular, the $H^1$ and $H^1 relative$ errors for the Darcy velocity $\textbf{u}_p$ are also accurate and can keep 1.5 order convergence. It means our method can also perform well in the second order derivatives and have no error accumulation, which is in good agreement with the above conclusion that it can show the advantage of the GFDM for the derivative functions. In Table 7, the $L_2$ and $H^1$ relative errors of the GFDM are presented. From this table, we can see that our method is accurate and high efficiency due to the CPU time is very small. The varying kinematic viscosity of fluid $v$ can be set $v=1.0,10^{-2},10^{-4},10^{-6},10^{-8}.$ The $L_2$ errors for the numerical velocity of Stokes $\textbf{u}_f$ (left) and Darcy $\textbf{u}_p$ (right) with varying $v$ are shown in Fig. 10. From these figures, we can see the numerical results support that the meshless GFDM are accurate and keep 2nd order convergence. Furthermore, some numerical results with the different realistic parameters $v=0.1,0.01,0,001$ and $K=1.0,0.1$ are shown in Fig. 11. We can see that the numerical results are also accurate and keep 2nd order convergence. It means that our meshless method can do well in the Stokes-Darcy coupled problem in realistic applications. In this example we also test the stability about the 'm'. Note that our numerical errors are all kept at the same accuracy and our results are all stable in the range from 12 to 22. In this example, we adopt $m=20.$   
\subsection{Example 3: The Stokes-Darcy coupled problem with closed complex interfaces.}
In this example, we consider the first problem with a circle interface (see Fig.12) $\Gamma:$$\varphi=x^2+y^2-0.25=0.$ This problem is inspired by Li[24], which consders the Stokes-Darcy coupled problem with closed interface and makes a channel flow outside the circle and the porous medis flow inside the interface. Furthermore, in order to show the advantage of the GFDM for complex interface, we consider the following complex interfaces:\\
(1) Two-petaled interface:
\begin{equation}
	r(\theta)=\left\{
	\begin{array}{cl}
		{x(\theta)=0.5cos(\theta)+0.2sin(2\theta)cos(\theta)},\\{y(\theta)=0.5sin(\theta)+0.2sin(2\theta)sin(\theta)},
	\end{array}
	0\leq\theta\leq2\pi.
	\right.
\end{equation} 
(2) Flower shaped interface:
\begin{equation}
	r(\theta)=\left\{
	\begin{array}{cl}
		{x(\theta)=0.5cos(\theta)+0.2sin(5\theta)cos(\theta)},\\{y(\theta)=0.5sin(\theta)+0.2sin(5\theta)sin(\theta)},
	\end{array}
	0\leq\theta\leq2\pi.
	\right.
\end{equation} 
(3) Heart shaped interface:
\begin{equation}
	r(\theta)=\left\{
	\begin{array}{cl}
		{x(\theta)=0.3(1-sin(\theta))cos(\theta)},\\{y(\theta)=0.3(1-sin(\theta))sin(\theta)+0.2},
	\end{array}
	0\leq\theta\leq2\pi.
	\right.
\end{equation} 
 \begin{figure}
	\centering
	\includegraphics[scale=.6]{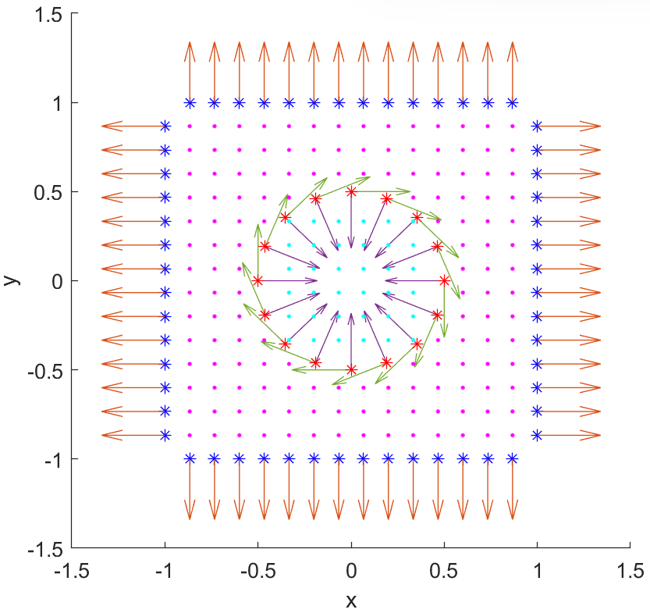}
	\raisebox{0.5\height}{\includegraphics[scale=.6]{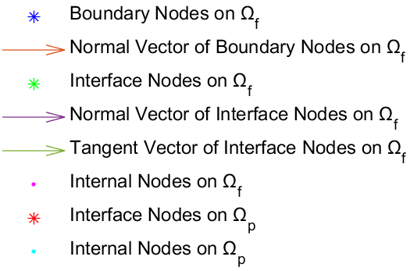}}
	\caption{ The point collocation for Example 3.}
\end{figure}
\begin{figure}
	\centering
		\includegraphics[scale=.4]{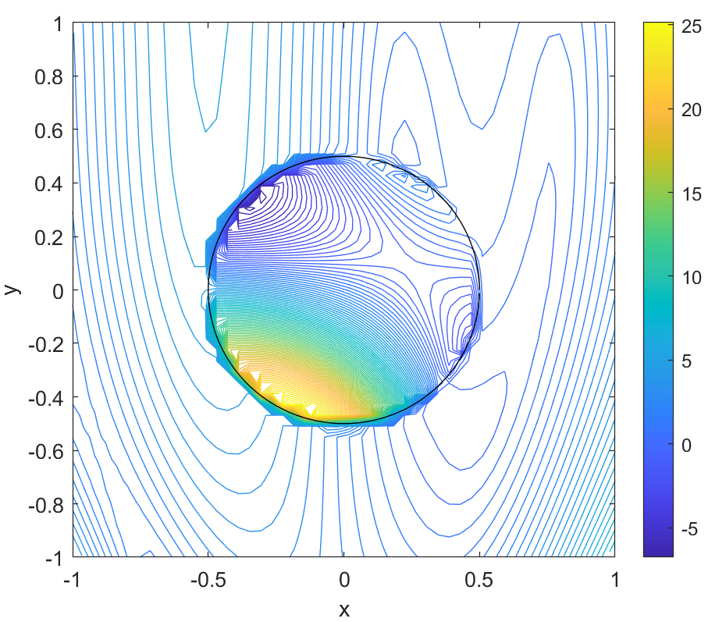}
	\hspace{1cm}
	\raisebox{0.1\height}{\includegraphics[scale=.4]{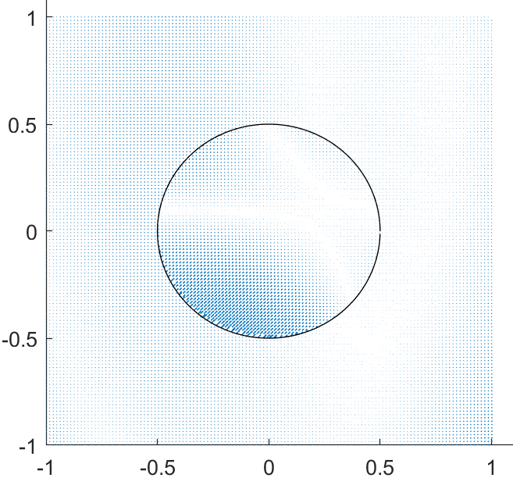}}
	
	\caption{ The contour (left) and vector (right) of numerical solution for Example 3.}
\end{figure}
\begin{figure}
	\centering
	\includegraphics[scale=.4]{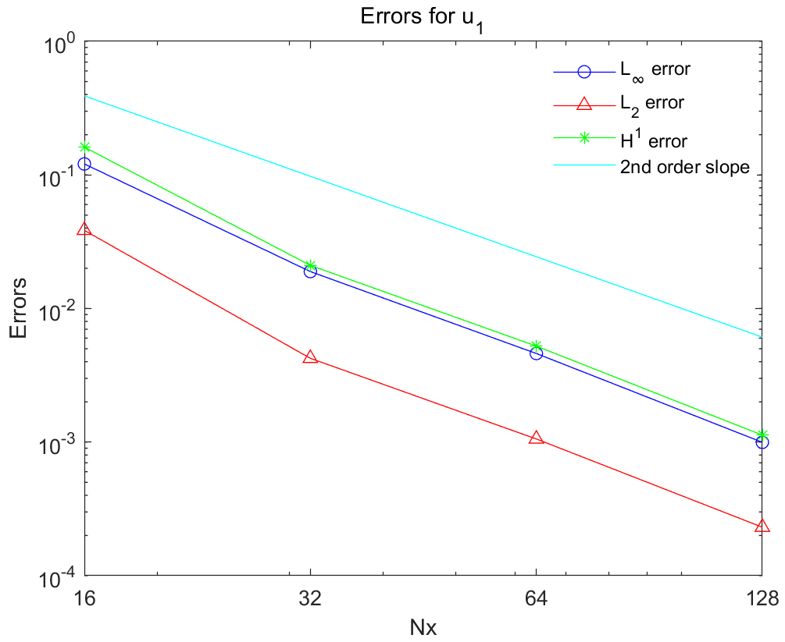}
	\includegraphics[scale=.4]{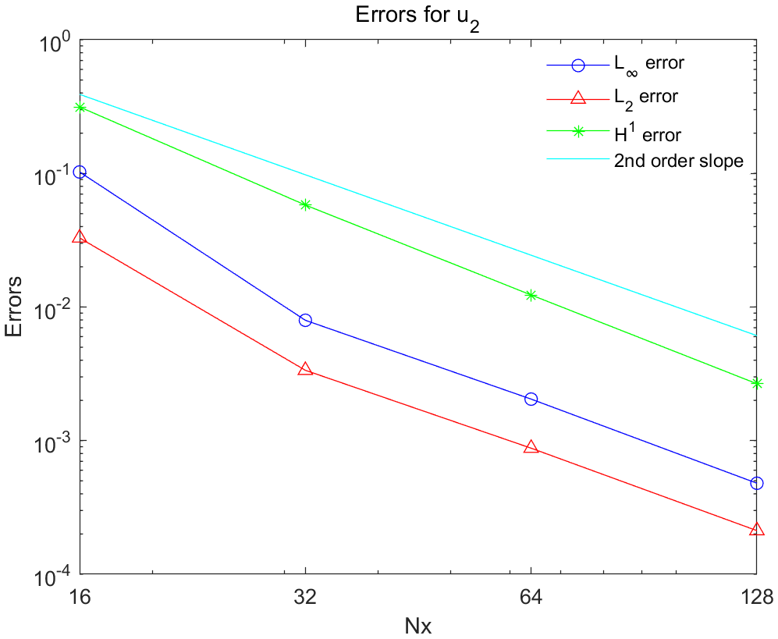}
	\includegraphics[scale=.4]{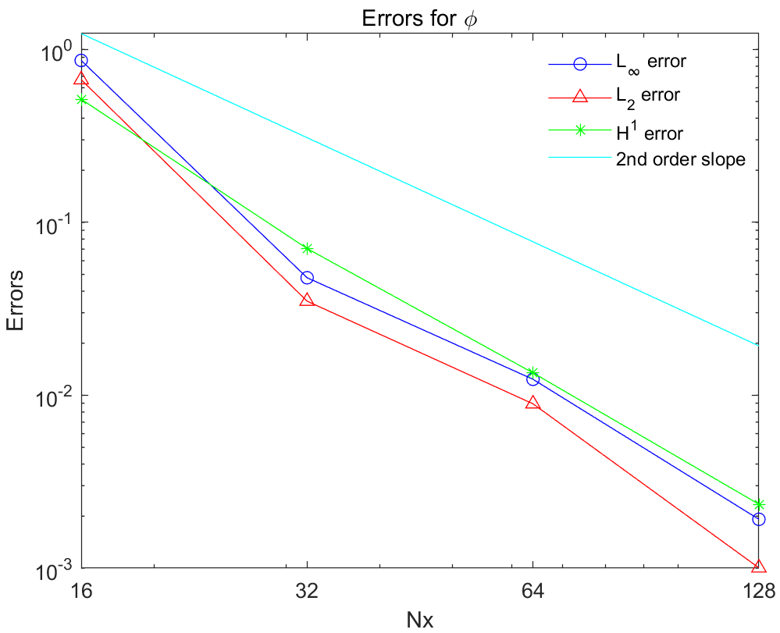}
	\caption{ $L_{\infty}$, $L_2$ and $H^1$ errors of the component functions for Example 3.}
\end{figure}

\begin{figure}
	\centering
	\includegraphics[scale=.4]{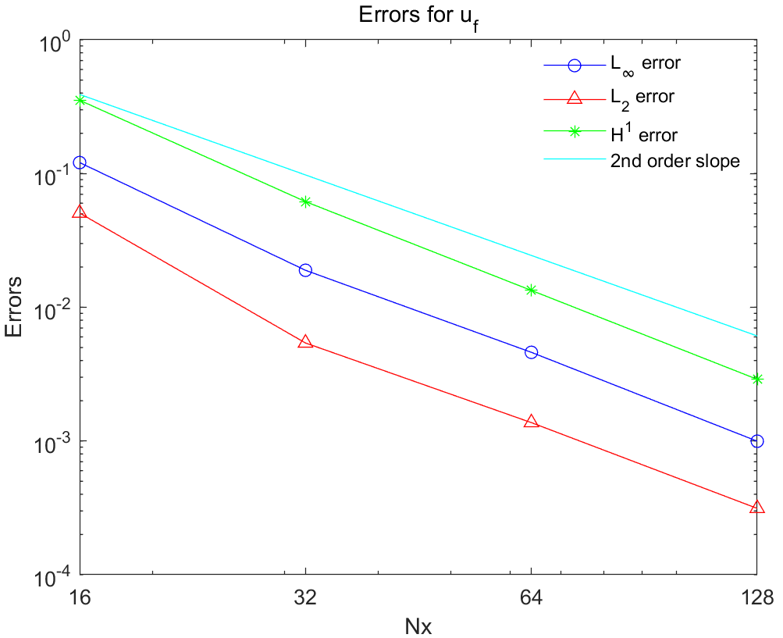}
	\includegraphics[scale=.4]{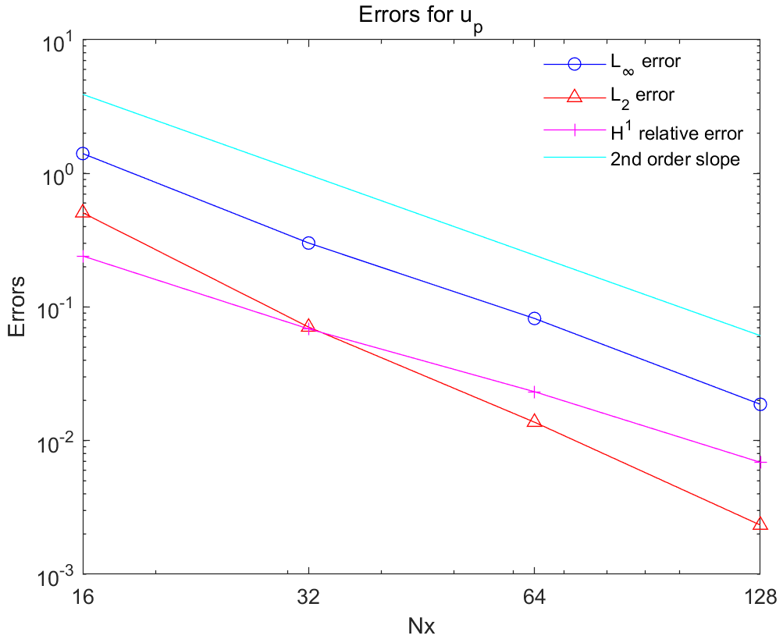}
	\includegraphics[scale=.4]{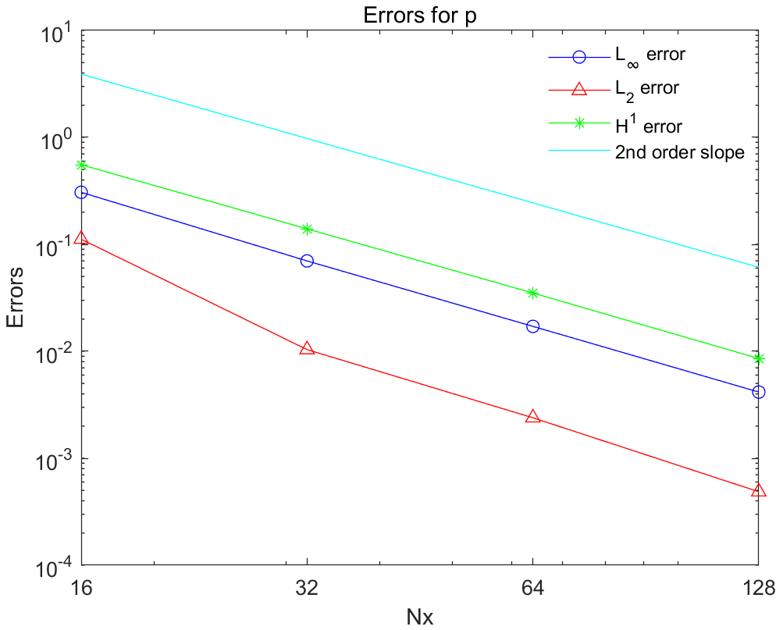}
	\caption{ $L_{\infty}$, $L_2$ and $H^1$ errors of fluid velocity $\textbf{u}_f$ (left), Darcy velocity $\textbf{u}_p$ (middle) and pressure $p$ (right) for Example 3.}
\end{figure}
	\begin{table*}	
	\scriptsize
	\caption{ The $L_2$ and $H^1$ relative errors and CPU(s) of the GFDM for Example 3}
	\begin{tabular}{ccccccccccc}
		\hline
		\multirow{1}{*}{$ $}&\multicolumn{1}{c}{$N_x$}& \multicolumn{2}{c}{$u_f$} & \multicolumn{1}{c}{$p$}& \multicolumn{1}{c}{$\phi$} & \multicolumn{2}{c}{$u_p$}
		&\multicolumn{1}{c}{$CPU(s)$}\\
		\hline
		
		$ $&$ $& $L_{2,relative}$      &  $H^1_{relative}$  
		& $L_{2,relative}$      &  $L_{2,relative}$  
		&  $L_{2,relative}$ &$H^1_{relative}$   \\ 
		\hline		
		
		2nd GFDM&$16$ &$2.74\times10^{-2}$ & $5.93\times10^{-2}$ & $3.24\times10^{-2}$ &  $1.88 \times10^{-1}$ & $3.40 \times10^{-2}$ &$2.41 \times10^{-1}$&$0.22$ \\
		
		$ $ &$32$ & $8.25\times10^{-3}$ & $1.48\times10^{-2}$ & $8.30\times10^{-3}$ &  $8.32\times10^{-2}$ & $6.69\times10^{-3}$& $8.13 \times10^{-2}$&$0.23$ \\
		
		$ $ &$64$& $1.63 \times10^{-3}$ & $2.85 \times10^{-3}$ & $1.66\times10^{-3}$ &  $ 1.82\times10^{-2}$ & $1.82 \times10^{-3}$& $3.11\times10^{-2}$&$0.52$ \\
		
	$ $ &$128$& $4.36\times10^{-4}$ & $6.93\times10^{-4}$ & $4.26\times10^{-4}$ &  $5.14\times10^{-3}$ & $3.62\times10^{-4}$& $1.00\times10^{-2}$&$2.10$ \\

		\hline
		
	\end{tabular}
\end{table*}
\begin{table*}	
	\scriptsize
	\caption{ $L_{\infty}$, $L_2$ and $H^1$ errors of the GFDM with different $m$  when $N_x=32$ for Example 3}
	\begin{tabular}{cccccccccccc}
		\hline
		\multirow{1}{*}{$m$} & \multicolumn{3}{c}{$u_f$} & \multicolumn{3}{c}{$u_p$}& \multicolumn{2}{c}{$p$}&
		\multicolumn{2}{c}{$\phi$} \\
		\hline
		
		& $L_{\infty}$      &  $L_2 $   &   $H^1$
		
		& $L_{\infty}$      &  $L_2$   &   $H^{1,relative}$
		
		&  $L_2$   &   $H^1$
		&  $L_2$   &   $H^1$\\
		\hline
		$12$ & $3.36\times10^{-2}$ & $7.13\times10^{-3}$ & $6.73\times10^{-2}$ &  $6.89\times10^{-1}$ & $1.21\times10^{-1}$ & $9.26\times10^{-2}$&
		$1.42\times10^{-2}$ & $1.12\times10^{-1}$ & $8.11\times10^{-2}$&$1.20\times10^{-1}$ \\

		$14$ & $2.73\times10^{-2}$ & $1.15\times10^{-2}$ & $7.43\times10^{-2}$ &  $4.62\times10^{-1}$ & $1.02\times10^{-1}$ & $7.84\times10^{-2}$&
		$2.48\times10^{-2}$ & $1.15\times10^{-1}$ & $1.65\times10^{-1}$ &$1.01\times10^{-1}$\\
		
		$16$ & $2.17\times10^{-2}$ & $9.24\times10^{-3}$ & $6.89\times10^{-2}$ &  $3.66\times10^{-1}$ & $8.74\times10^{-2}$ & $7.27\times10^{-2}$&
		$1.92\times10^{-2}$ & $1.09\times10^{-1}$ & $1.23\times10^{-1}$&$8.68\times10^{-2}$ \\
		
		$18$ & $1.97\times10^{-2}$ & $5.77\times10^{-3}$ & $6.31\times10^{-2}$ &  $3.06\times10^{-1}$ & $7.16\times10^{-2}$ & $6.57\times10^{-2}$&
		$1.06\times10^{-2}$ & $1.26\times10^{-1}$ & $4.49\times10^{-2}$&$7.18\times10^{-2}$ \\
		
		$20$ & $1.89\times10^{-2}$ & $5.41\times10^{-3}$ & $6.17\times10^{-2}$ &  $3.01\times10^{-1}$ & $7.08\times10^{-2}$ & $6.86\times10^{-2}$&
		$1.04\times10^{-2}$ & $1.40\times10^{-1}$ & $3.50\times10^{-2}$&$7.06\times10^{-2}$ \\
		
		$22$ & $1.89\times10^{-2}$ & $5.88\times10^{-3}$ & $6.28\times10^{-2}$ &  $3.13\times10^{-1}$ & $7.44\times10^{-2}$ & $7.34\times10^{-2}$&
		$1.06\times10^{-2}$ & $1.37\times10^{-1}$ & $2.99\times10^{-2}$&$7.39\times10^{-2}$ \\
		\hline
	\end{tabular}
\end{table*}
 \begin{figure}
	\centering
	\includegraphics[scale=.4]{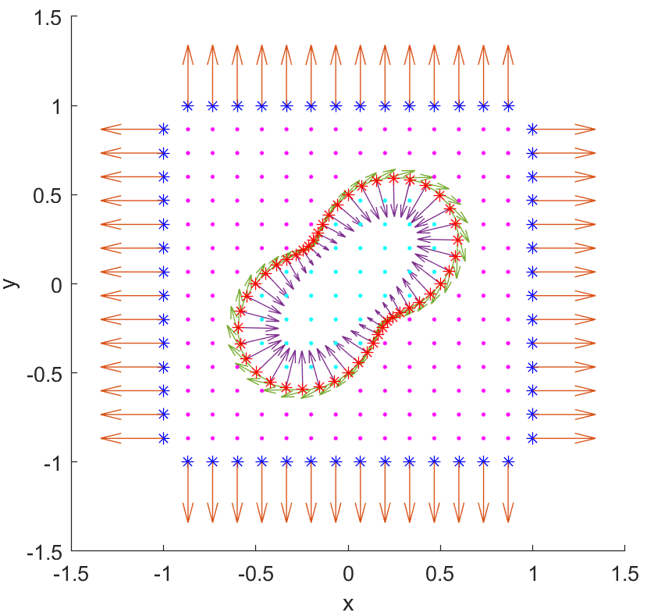}
		\includegraphics[scale=.4]{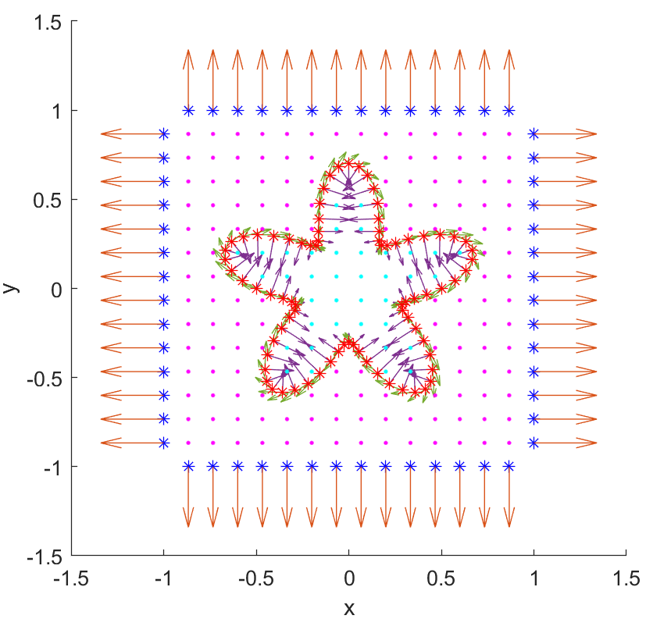}
			\includegraphics[scale=.4]{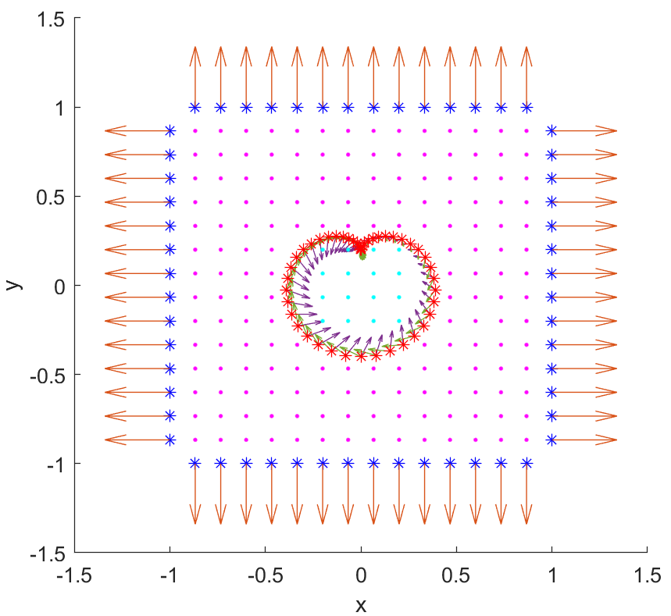}
	\caption{ The point collocations under different complex interfaces for Example 3.}
\end{figure}
\begin{figure}
	\centering
	\includegraphics[scale=.4]{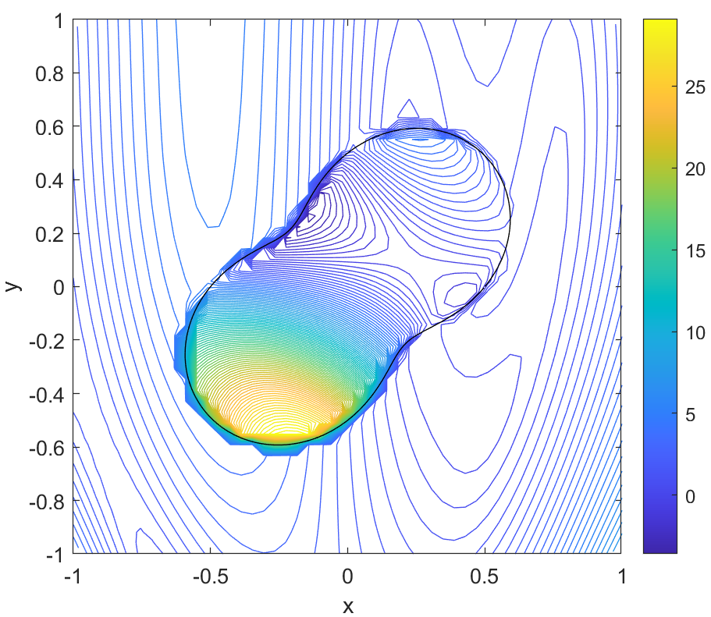}
	\includegraphics[scale=.4]{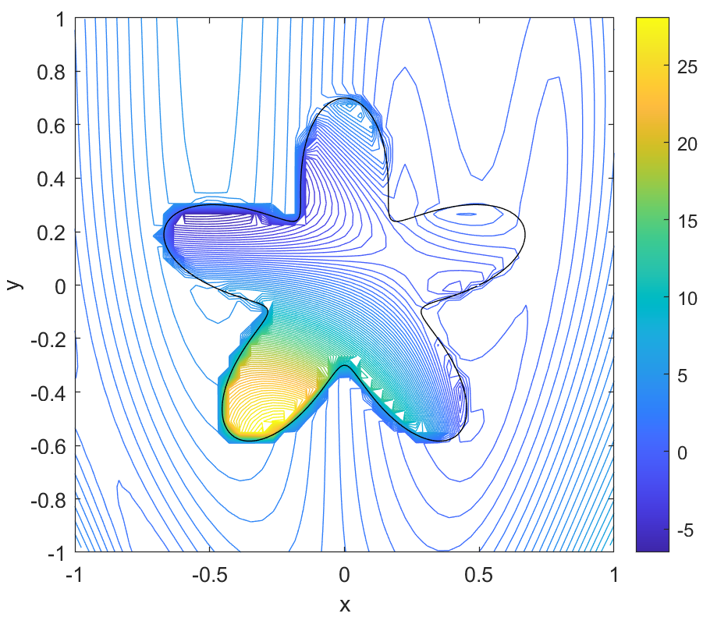}
	\includegraphics[scale=.4]{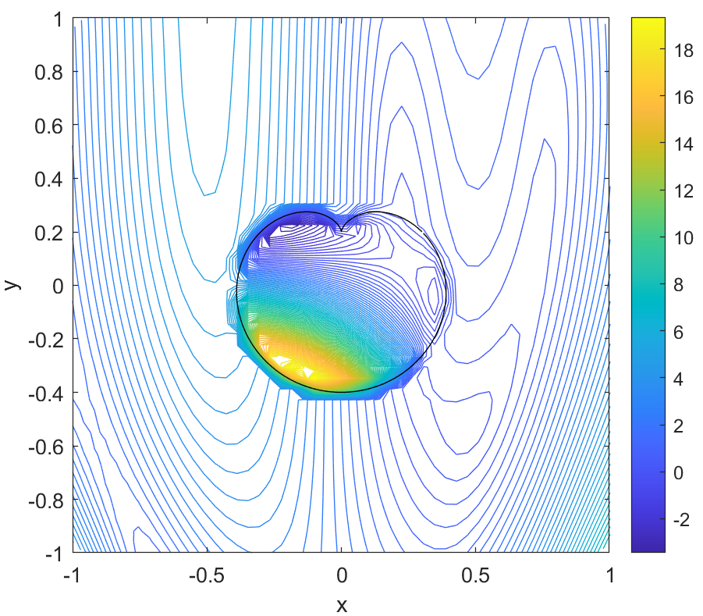}
	
	\caption{ The contour of numerical solutions under different complex interfaces for Example 3.}
\end{figure}
\begin{figure}
	\centering
	\includegraphics[scale=.4]{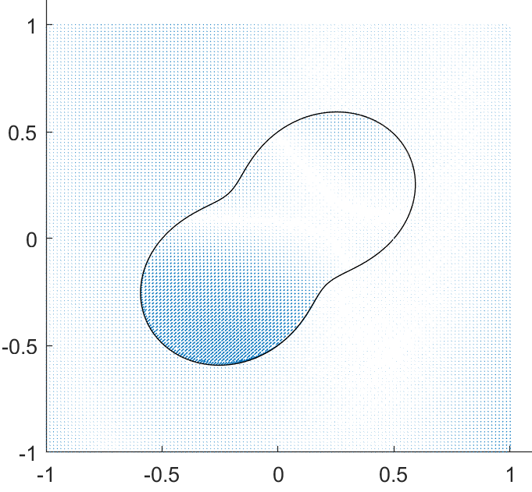}
		\hspace{1cm}
	\includegraphics[scale=.4]{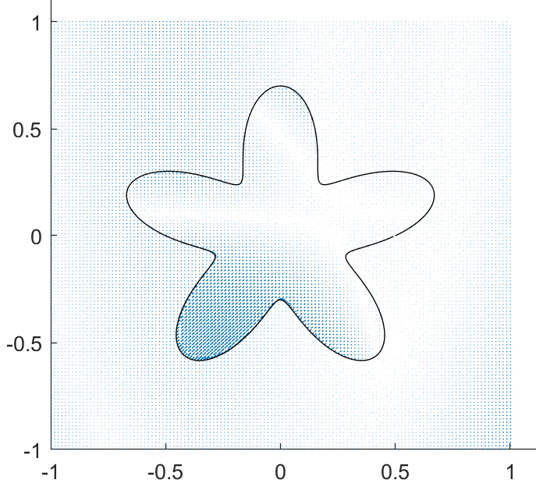}
		\hspace{1cm}
	\includegraphics[scale=.4]{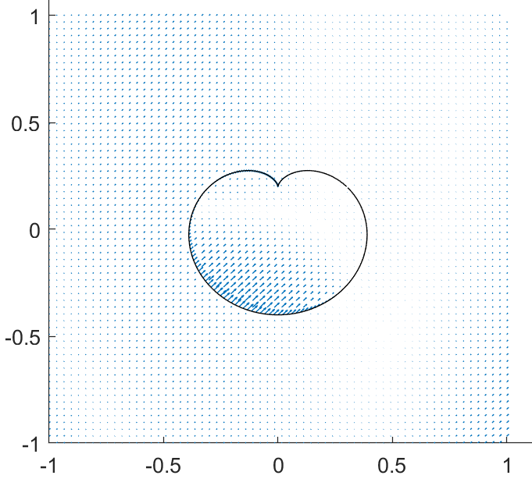}	
	\caption{ The vector of numerical solutions under different complex interfaces for Example 3.}
\end{figure}

\begin{figure}
	\centering
	\includegraphics[scale=.4]{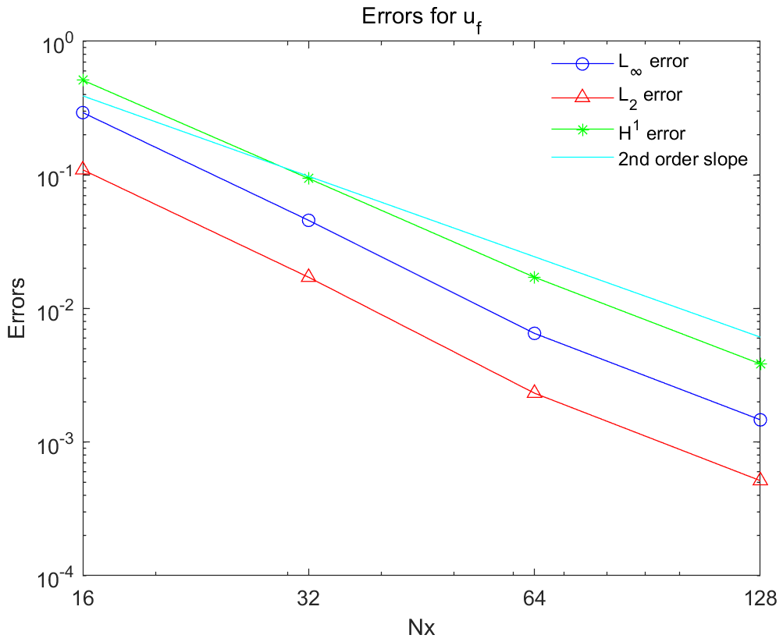}
	\includegraphics[scale=.4]{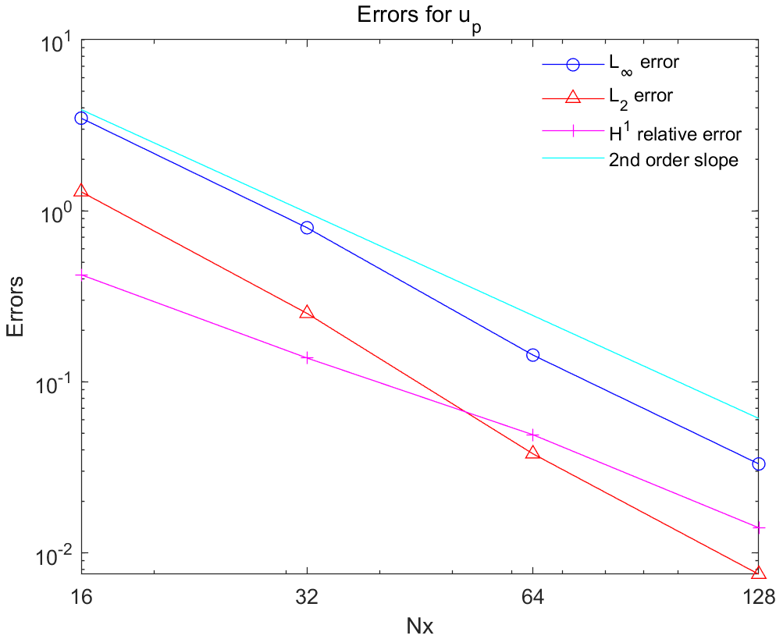}
	\includegraphics[scale=.4]{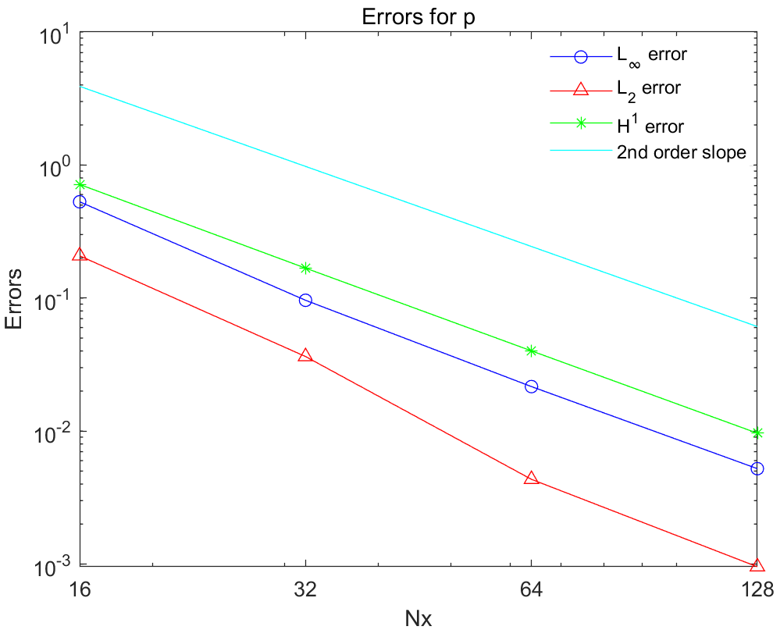}
	\caption{ $L_{\infty}$, $L_2$ and $H^1$ errors of fluid velocity $\textbf{u}_f$ (left), Darcy velocity $\textbf{u}_p$ (middle) and pressure $p$ (right) under the two-petaled interface for Example 3.}
\end{figure}
\begin{figure}
	\centering
	\includegraphics[scale=.4]{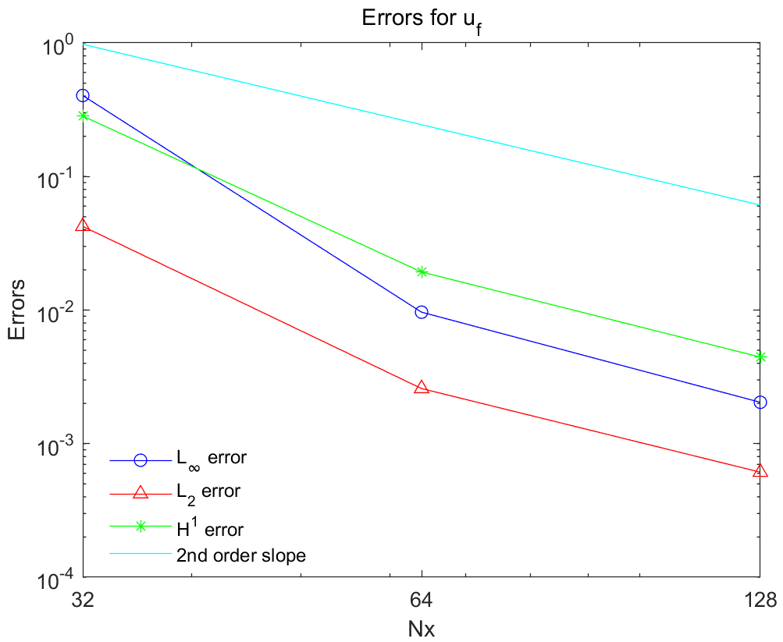}
	\includegraphics[scale=.4]{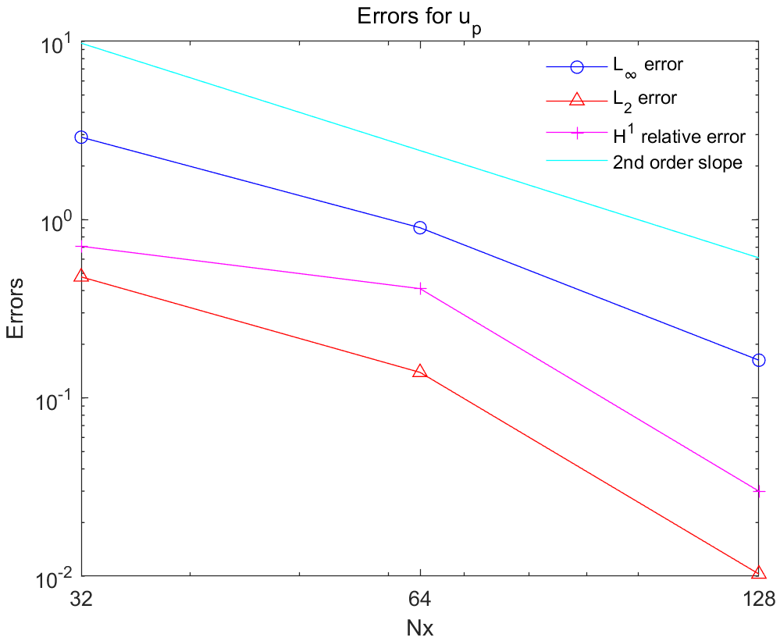}
	\includegraphics[scale=.4]{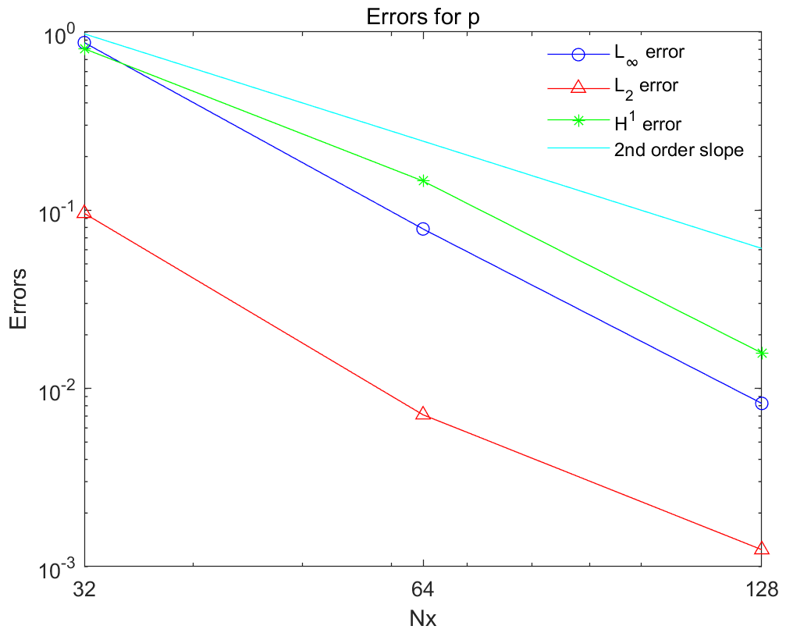}
	\caption{ $L_{\infty}$, $L_2$ and $H^1$ errors of fluid velocity $\textbf{u}_f$ (left), Darcy velocity $\textbf{u}_p$ (middle) and pressure $p$ (right) under flower shaped interface for Example 3.}
\end{figure}
\begin{figure}
	\centering
	\includegraphics[scale=.4]{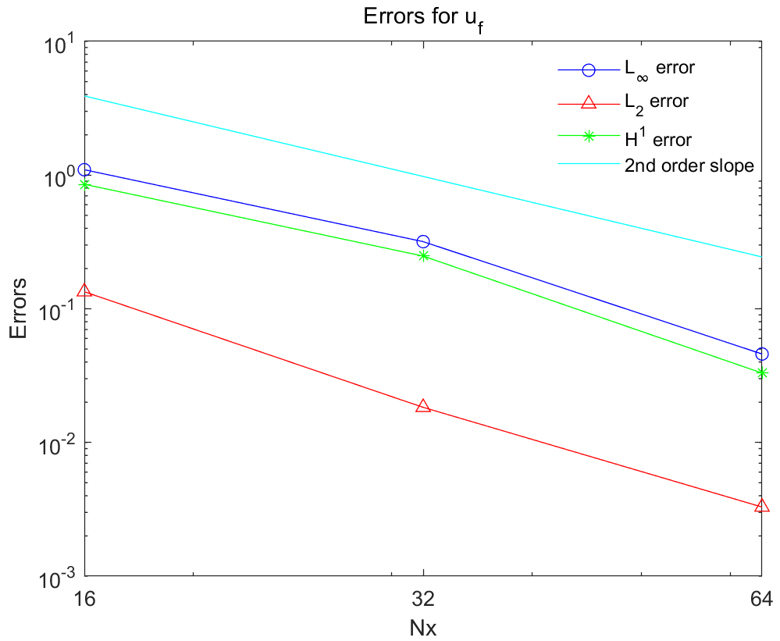}
	\includegraphics[scale=.4]{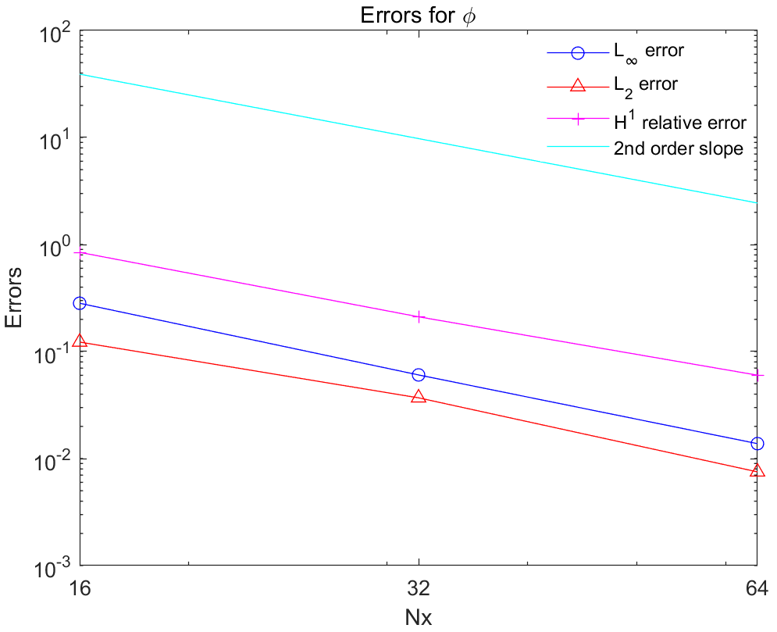}
	\includegraphics[scale=.4]{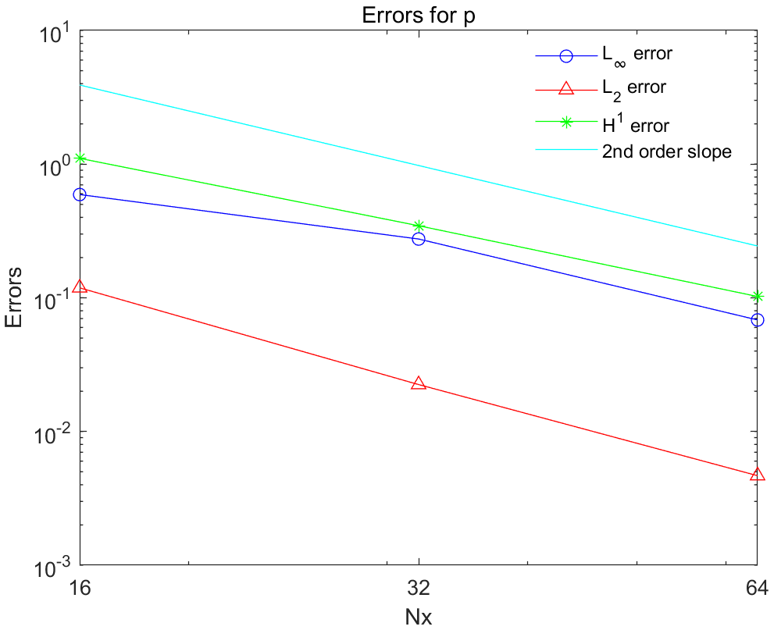}
	\caption{ $L_{\infty}$, $L_2$ and $H^1$ errors of fluid velocity $\textbf{u}_f$ (left), Darcy velocity $\textbf{u}_p$ (middle) and pressure $p$ (right) under heart shaped interface for Example 3.}
\end{figure}
\begin{table*}	
	\scriptsize
	\caption{ $L_{\infty}$, $L_2$ and $H^1$ errors of the GFDM with different $m$  when $N_x=32$ for Example 3}
	\begin{tabular}{cccccccccccc}
		\hline
		\multirow{1}{*}{$ $} &
		\multirow{1}{*}{$m$} & \multicolumn{2}{c}{$u_f$} & \multicolumn{1}{c}{$u_p$}& \multicolumn{2}{c}{$p$}&
		\multicolumn{2}{c}{$\phi$} \\
		\hline
		
		&$ $  &  $L_2 $ &  $H^1$  
		
        &  $L_2$
		
		&  $L_2$ &$H^1$
		&  $L_2$   &   $H^1$\\
		\hline
		Two-petaled interface &$22$ &  $1.61\times10^{-2}$ &
		$8.11\times10^{-2}$ & $2.12\times10^{-1}$&
		$3.69\times10^{-2}$ & $1.65\times10^{-1}$ & $9.13\times10^{-2}$&$1.65\times10^{-1}$ \\

		$ $&$26$ &  $1.75\times10^{-2}$ &
		$8.90\times10^{-2}$ & $2.34\times10^{-1}$&
		$3.85\times10^{-2}$ & $1.65\times10^{-1}$ & $9.28\times10^{-2}$ &$1.79\times10^{-1}$\\
		
		$$&$30$ &  $1.66\times10^{-2}$ &
		$8.70\times10^{-2}$ & $2.29\times10^{-1}$&
		$3.59\times10^{-2}$ & $1.63\times10^{-1}$ & $8.61\times10^{-2}$&$1.79\times10^{-1}$ \\
		
		$$&$32$ &  $1.81\times10^{-2}$ &
		$9.10\times10^{-2}$ & $2.47\times10^{-1}$&
		$3.89\times10^{-2}$ & $1.64\times10^{-1}$ & $9.43\times10^{-2}$&$1.93\times10^{-1}$ \\
		
		$$&$36$ &  $1.71\times10^{-2}$ &
		$9.40\times10^{-2}$ & $2.51\times10^{-1}$&
		$3.63\times10^{-2}$ & $1.68\times10^{-1}$ & $8.87\times10^{-2}$&$1.96\times10^{-1}$ \\
		
		$$&$40$ &  $2.08\times10^{-2}$ &
		$1.09\times10^{-1}$ & $2.85\times10^{-1}$&
		$4.36\times10^{-2}$ & $1.85\times10^{-1}$ & $1.06\times10^{-1}$&$2.23\times10^{-1}$ \\
		\hline
		
		Flower shaped interface &$22$ &  $1.42\times10^{-2}$ &
		$8.81\times10^{-2}$ & $1.68\times10^{-1}$&
		$2.84\times10^{-2}$ & $2.80\times10^{-1}$ & $8.39
		\times10^{-2}$&$1.68\times10^{-1}$ \\

	$ $&$26$ &  $1.40\times10^{-2}$ &
	$8.90\times10^{-2}$ & $1.99\times10^{-1}$&
	$2.90\times10^{-2}$ & $2.68\times10^{-1}$ & $8.57\times10^{-2}$ &$1.98\times10^{-1}$\\
	
	$$&$30$ 
	& $3.65\times10^{-2}$ &
	$2.52\times10^{-1}$ & $1.30\times10^{-0}$&
	$8.46\times10^{-2}$ & $1.49\times10^{-0}$ & $2.52
	\times10^{-1}$&$1.30\times10^{-0}$ \\
	
	$$&$32$ 
	& $3.97\times10^{-2}$ &
	$2.62\times10^{-1}$ & $1.34\times10^{-0}$&
	$8.47\times10^{-2}$ & $1.48\times10^{-0}$ & $2.47\times10^{-1}$&$1.33\times10^{-0}$ \\
	
	$$&$36$ &
	 $4.23\times10^{-2}$ &
	 $2.84\times10^{-1}$ & $4.77\times10^{-1}$&
	$9.60\times10^{-2}$ & $8.06\times10^{-1}$ & $1.21\times10^{-1}$&$5.05\times10^{-1}$ \\
	
	$$&$40$ & 
	 $3.07\times10^{-2}$ &
	 $2.10\times10^{-1}$ & $5.28\times10^{-1}$&
	$7.64\times10^{-2}$ & $6.52\times10^{-1}$ & $1.33\times10^{-1}$&$5.17\times10^{-1}$ \\
	\hline
	\end{tabular}
\end{table*}

Fig.12 presents the point collocation, we can see the point distribution of the circle interface. From Fig. 13, we can see the distribution of the numerical solution. In Fig.14 and Fig. 15, we can see the closed circle interface doesn't have any influence on the numerical results of the Stokes-Darcy coupled problem. In particular, the $H^1 relative$ errors for the Darcy velocity $\textbf{u}_p$ are also accurate and can keep 2nd order convergence. It means that the meshless GFDM is efficient for this circle interface problem. Table 9 presents the $L_2$ and $H^1$ relative errors of the GFDM and the numerical results are accurate and highly efficient because the CPU time is very small. The interface shape may need more 'm', therefore, we also test the stability about the 'm'. Note that our numerical errors are all kept at almost the same accuracy and our results are all stable in the range from 12 to 22. In this example, we adopt $m=20.$

In order to test the stability of the complexity of the interface shape, we consider three types of complex interfaces, the point collocation of each one are shown in Fig.16, the explanation of points and vectors is the same as the Fig. 12, we can see that the interface shape and the distributions of the normal vector and the tangent. The contour and vector of the numerical solution are shown in Fig.17 and Fig. 18, we observe that the regular distribution of the numerical solutions. The $L_{\infty}$, $L_2$ and $H^1$ errors of fluid velocity $\textbf{u}_f$ (left), Darcy velocity $\textbf{u}_p$ (middle) and pressure $p$ (right) under these three types of interfaces are shown in Fig.19-Fig.21. Note that the numerical errors are all efficient and keep 2nd order convergence, but the $H^1$ error of the $\textbf{u}_p$ has a little oscillation, we don't show the lines. Therefore, there is little impact on the numerical results of the GFDM for the Stokes-Darcy coupled problem with the complex interfaces. The interface complexity may need more 'm', Therefore, we test the stability about the 'm' in the range 22 to 40 in Table. 11. Note that our numerical errors are  kept at almost the same accuracy and our results are all stable, even in the number of points is small. In this subsection, we adopt $m=36.$
\subsection{Example 4: The Stokes-Darcy coupled problem with the moving interface.}
In this example, we also consider the first problem with two moving circle interfaces (see Fig.22) from Ref.[36], in which we can see the specific expression and consider $n=5.$ This problem is also inspired by Li[24], which makes a channel flow outside the circle and the porous medis flow inside the interface, and Xing and Zheng [36], which considers the anisotropic elliptic problem with moving interface. In this example, we test the influence of the interface location for the Stokes-Darcy coupled problem.  
\begin{figure}
	\centering
	\includegraphics[scale=.4]{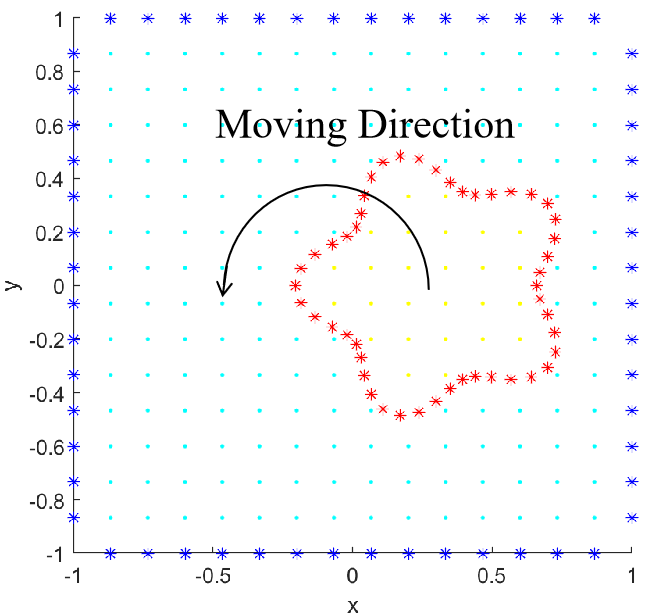}
	\hspace{1cm}
	\includegraphics[scale=.4]{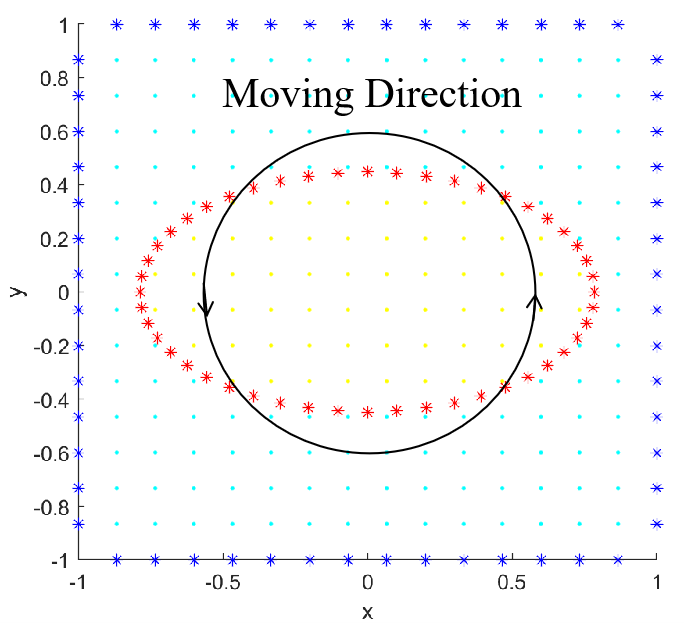}

	\caption{ The point collocations and moving directions of the proposed problems with the petagon interface and ellipse interface for Example 4.}
\end{figure}

\begin{figure}
	\centering
	\includegraphics[scale=.4]{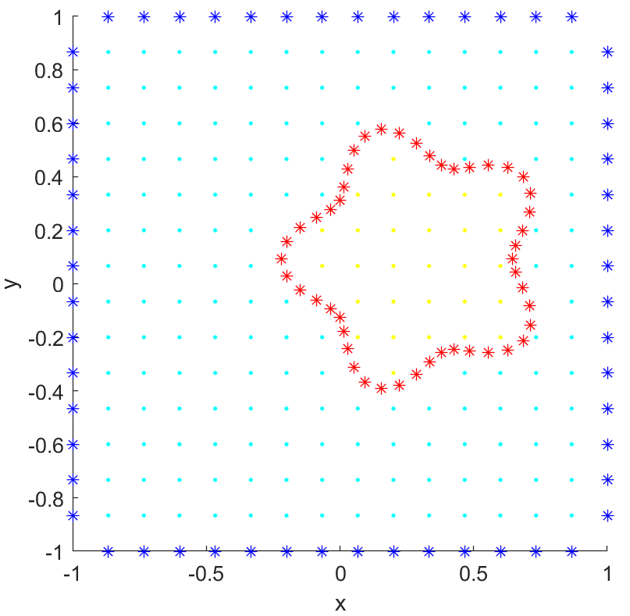}
		\hspace{1cm}
\includegraphics[scale=.4]{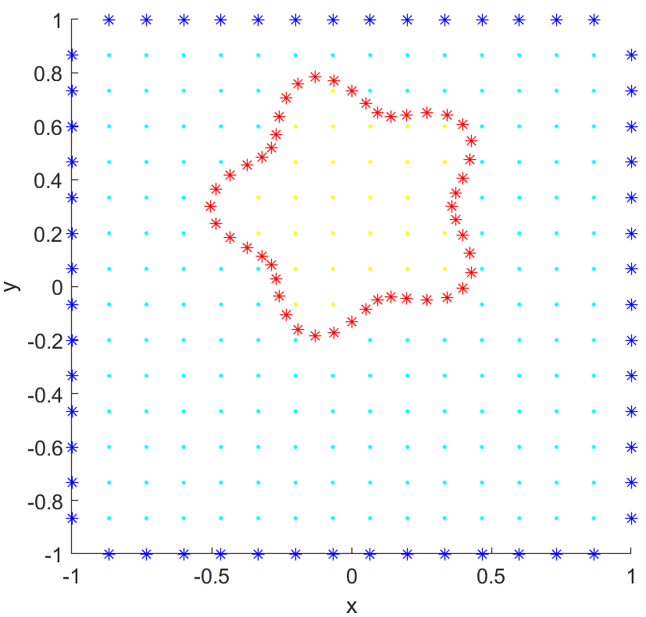}
	\hspace{1cm}
\includegraphics[scale=.4]{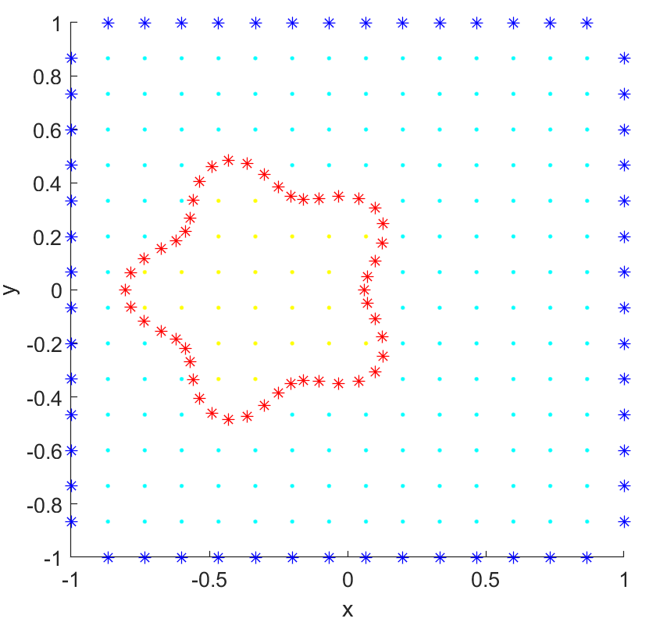}
	
	\caption{ The point collocation when $t=0.1$(left), $t=0.5$ (middle) and $t=1$ (right) under the petagon interface for Example 4.}
\end{figure}
\begin{figure}
	\centering
	\includegraphics[scale=.35]{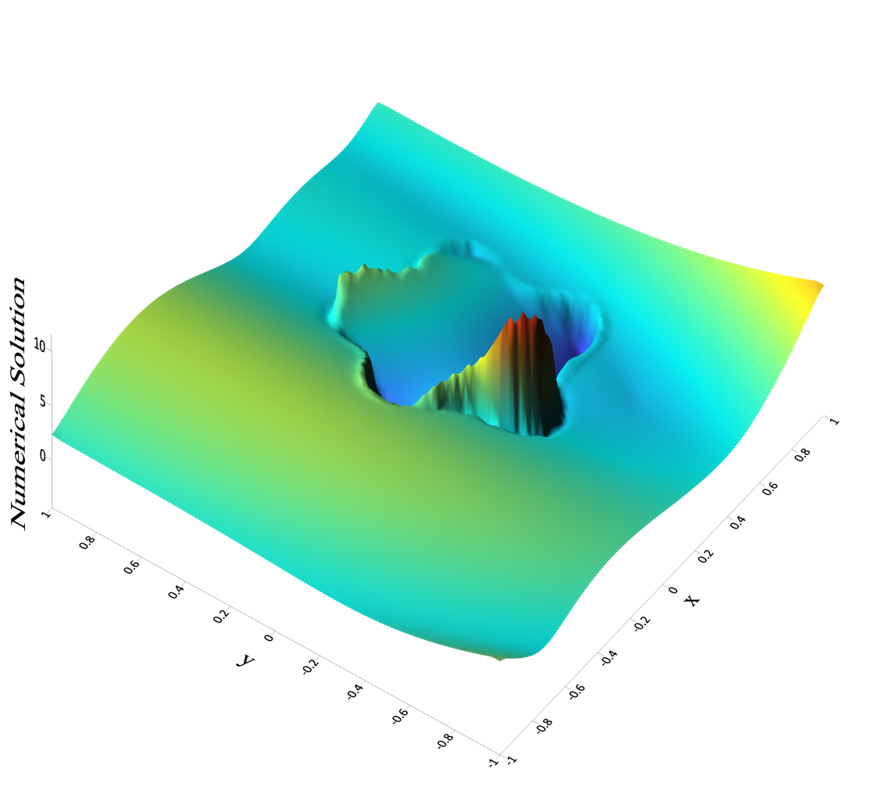}
	\includegraphics[scale=.35]{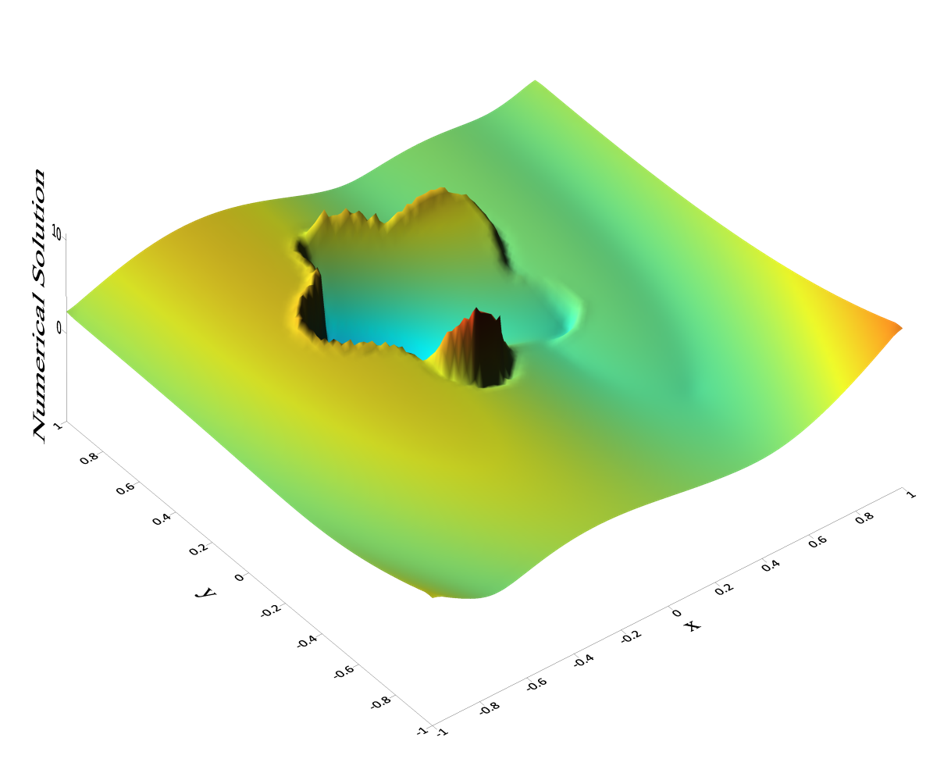}
	\includegraphics[scale=.35]{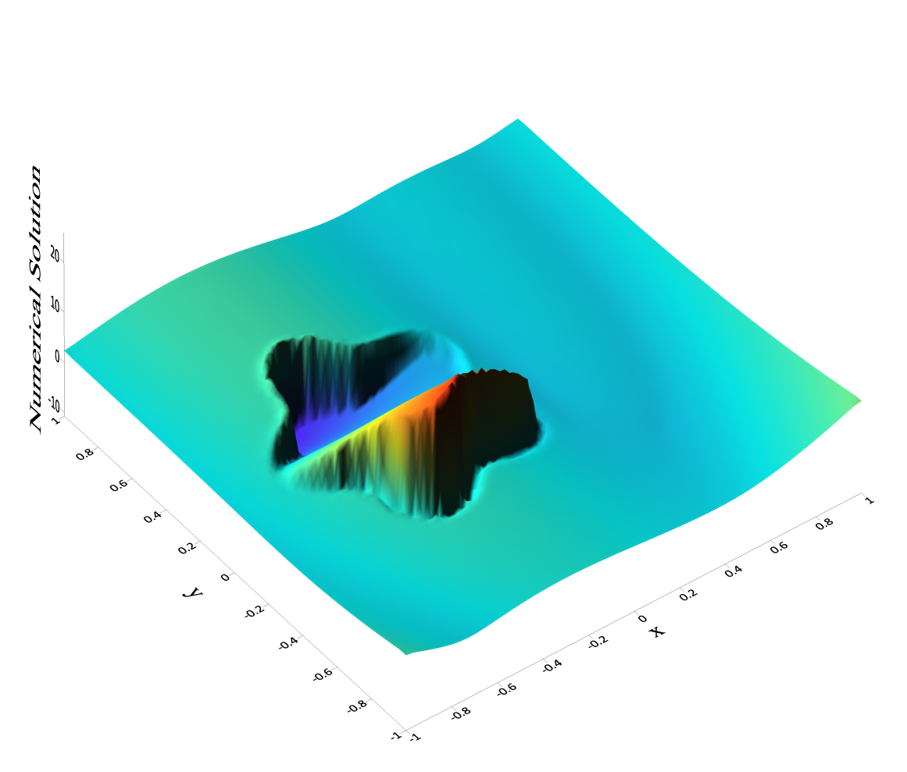}
	
	\caption{ The numerical solution when $t=0.1$(left), $t=0.5$ (middle) and $t=1$ (right) under the petagon interface for Example 4}
\end{figure}
\begin{figure}
	\centering
	\includegraphics[scale=.4]{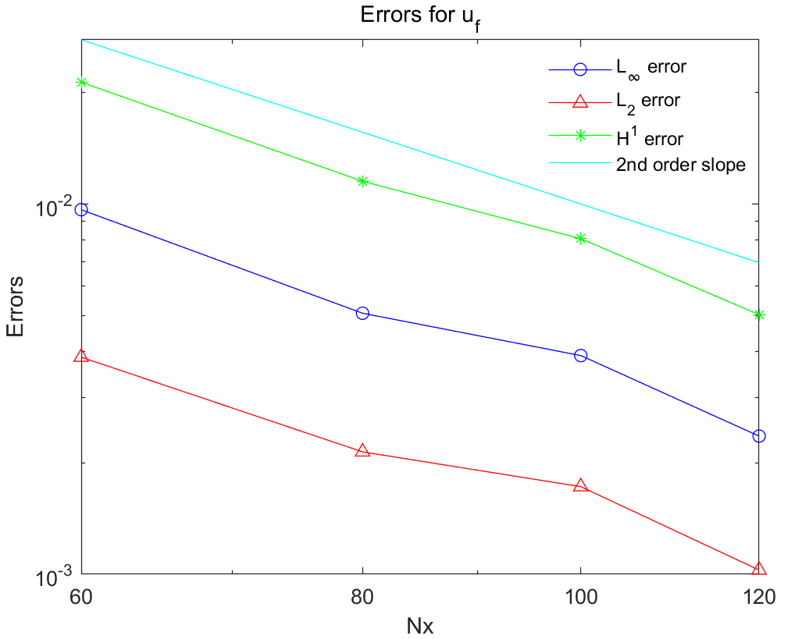}
	\includegraphics[scale=.4]{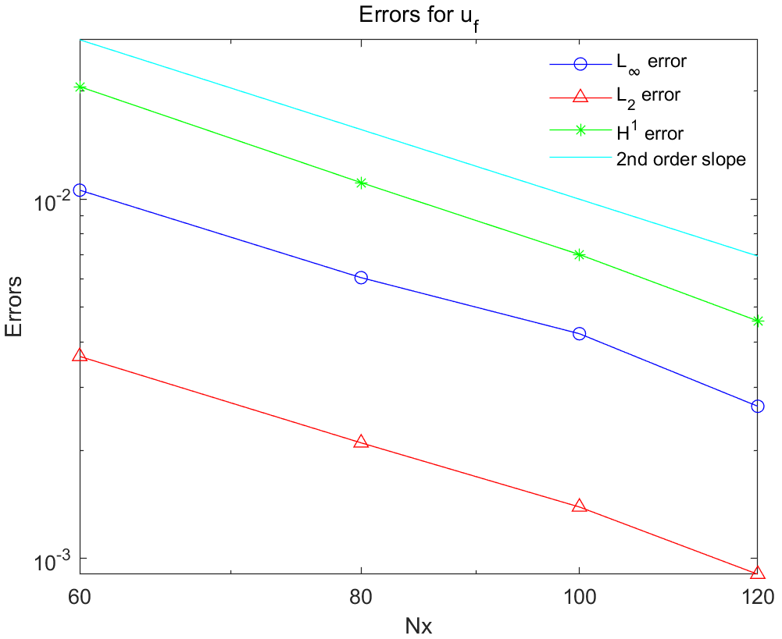}
	\includegraphics[scale=.4]{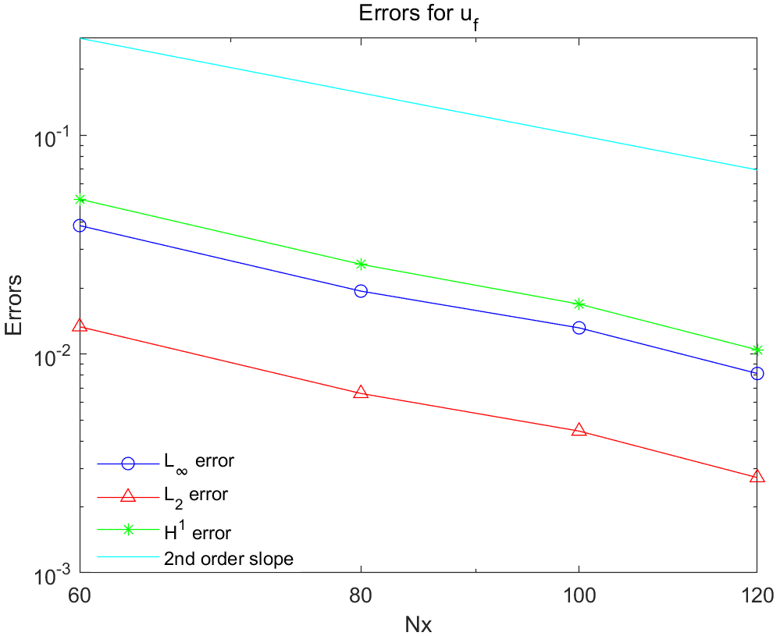}
	
	\caption{ $L_{\infty}$,$L_2$ and $H^1$ errors of $\textbf{u}_f$ when $t=0.1$(left), $t=0.5$ (middle) and $t=1$ (right) under the petagon interface for Example 4.}
\end{figure}

\begin{figure}
	\centering
	\includegraphics[scale=.4]{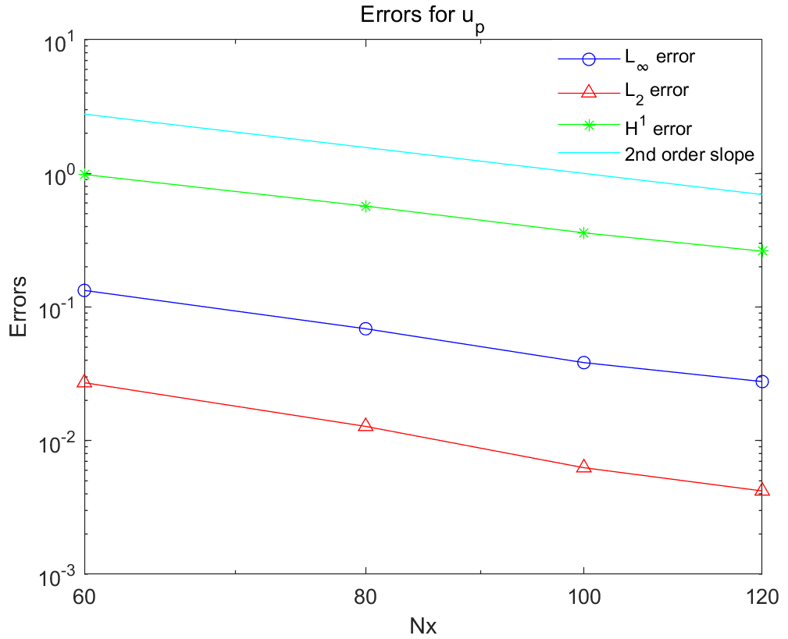}
	\includegraphics[scale=.4]{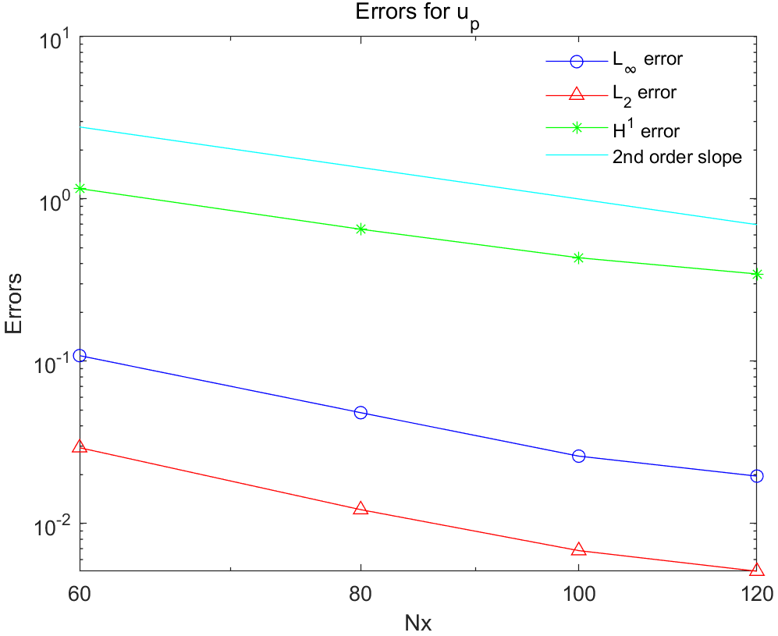}
	\includegraphics[scale=.4]{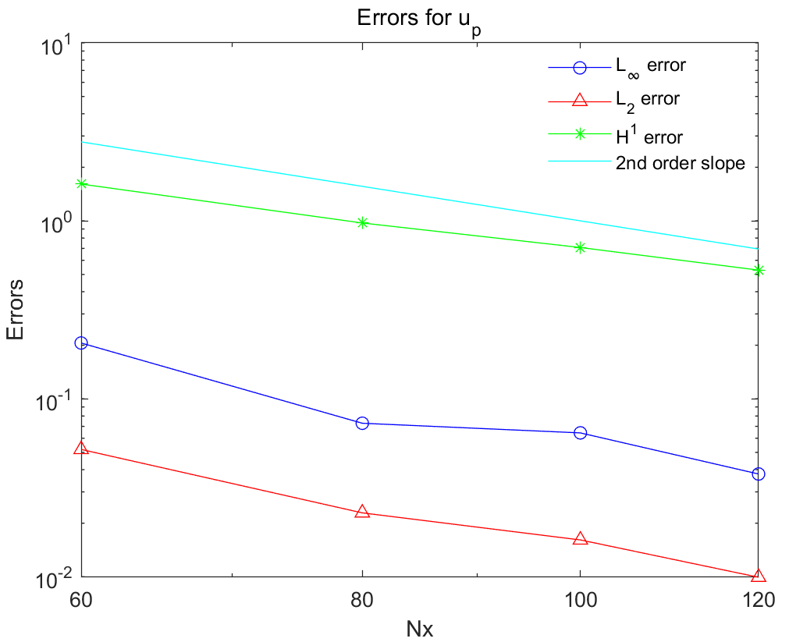}
	
	\caption{$L_{\infty}$,$L_2$ and $H^1$ errors of $\textbf{u}_p$ when $t=0.1$(left), $t=0.5$ (middle) and $t=1$ (right) under the petagon interface for Example 4.}
\end{figure}

\begin{figure}
	\centering
	\includegraphics[scale=.4]{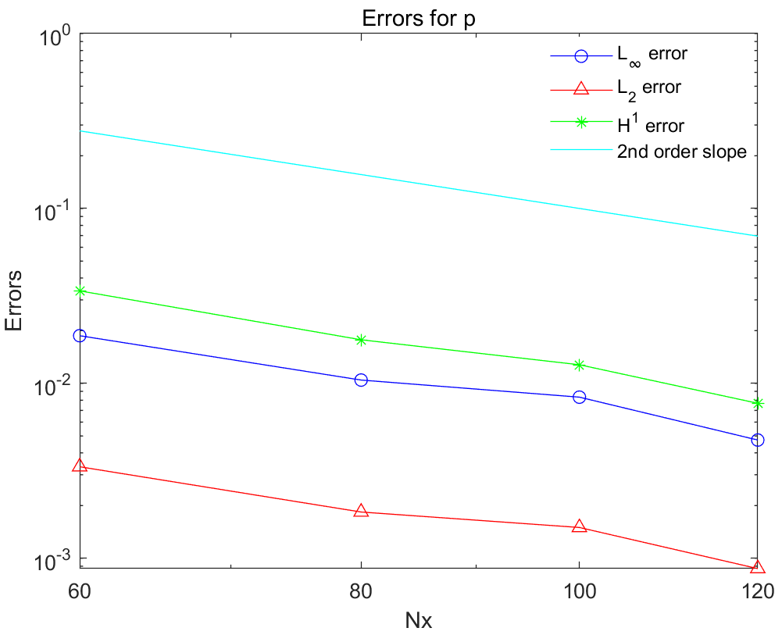}
	\includegraphics[scale=.4]{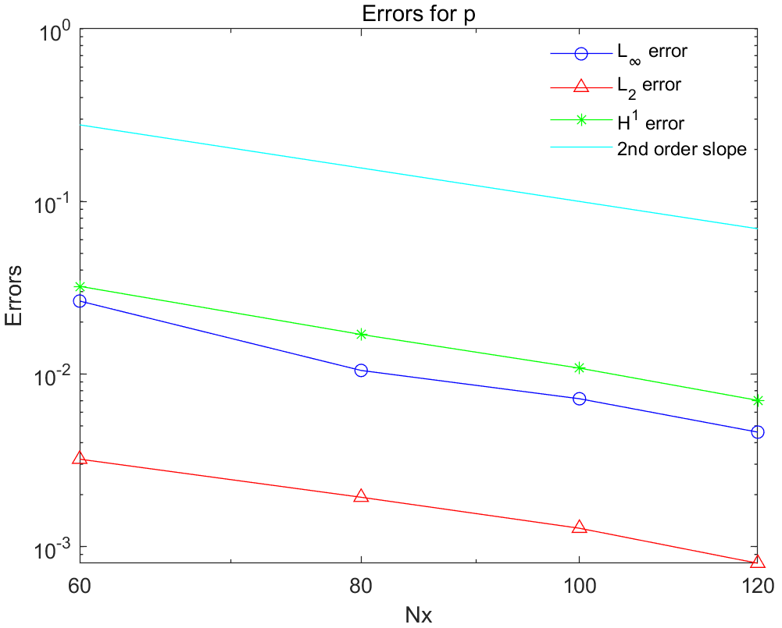}
	\includegraphics[scale=.4]{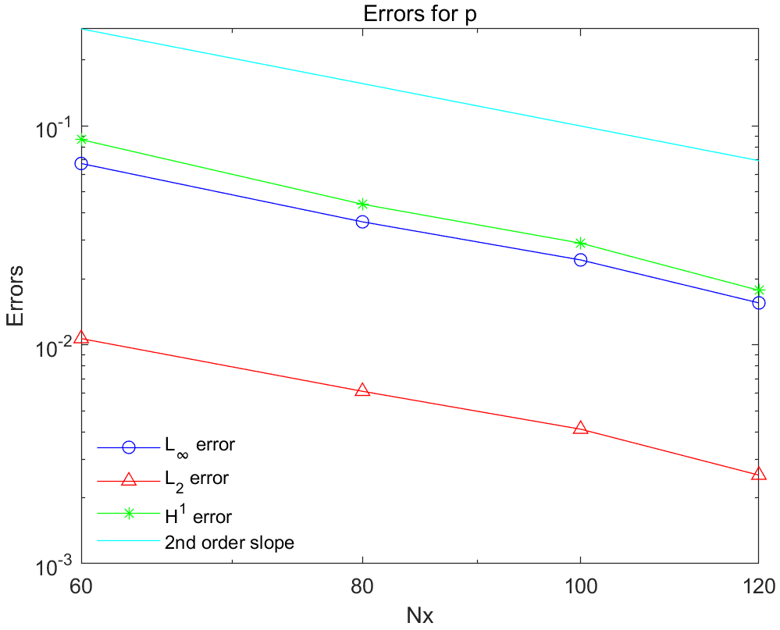}
	
	\caption{$L_{\infty}$,$L_2$ and $H^1$ errors of $p$ when $t=0.1$(left), $t=0.5$ (middle) and $t=1$ (right) under the petagon interface for Example 4.}
\end{figure}

\begin{figure}
	\centering
	\includegraphics[scale=.4]{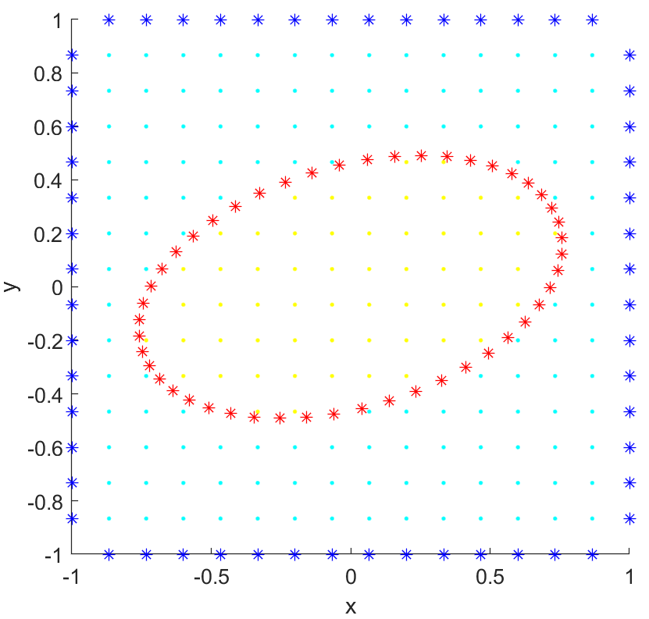}
		\hspace{1cm}
	\includegraphics[scale=.4]{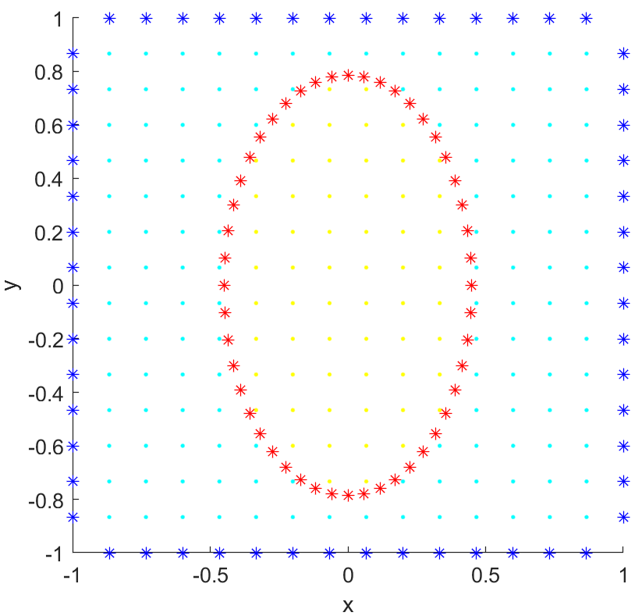}
		\hspace{1cm}
	\includegraphics[scale=.4]{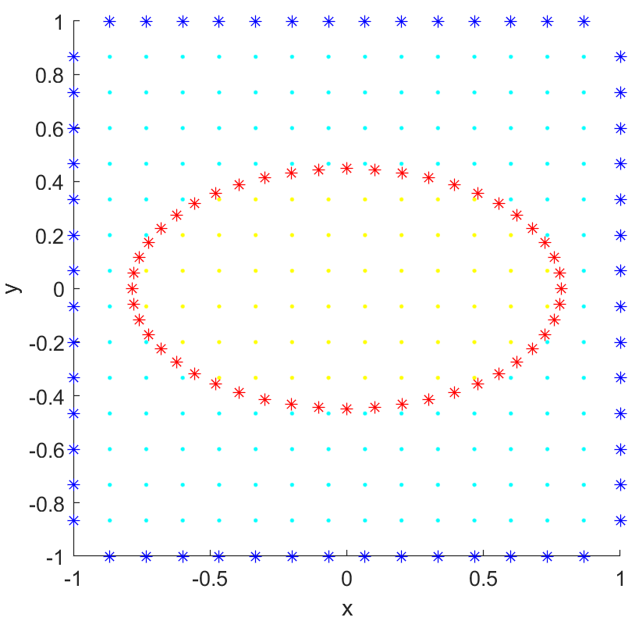}
	
	\caption{ The point collocation when $t=0.1$(left), $t=0.5$ (middle) and $t=1$ (right) under the ellipse interface for Example 4.}
\end{figure}
\begin{figure}
	\centering
	\includegraphics[scale=.35]{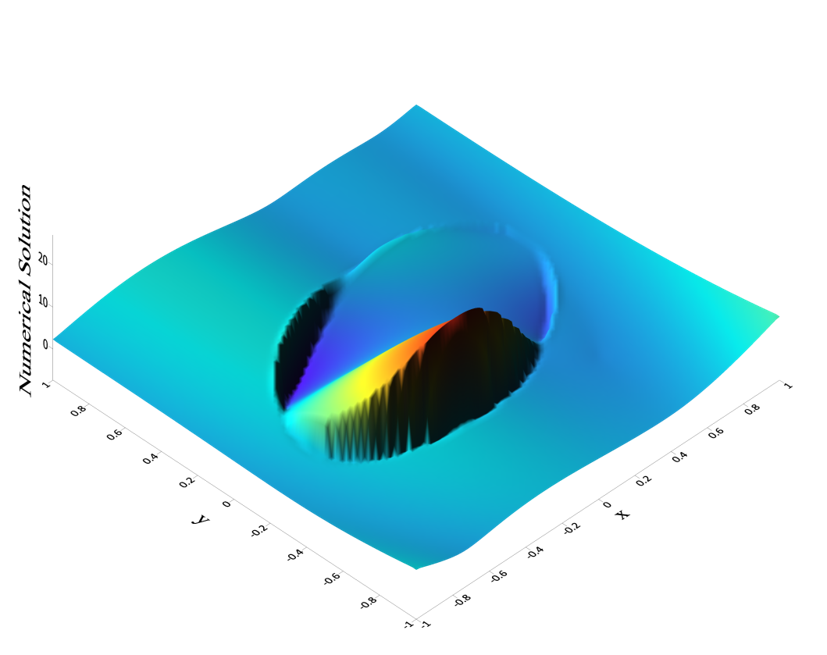}
	\includegraphics[scale=.35]{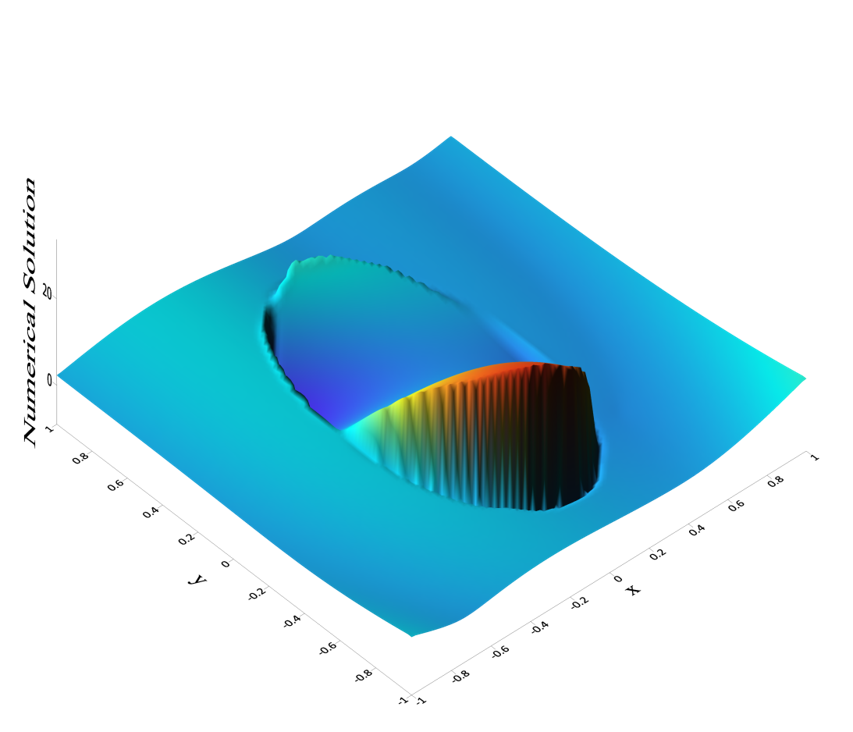}
	\includegraphics[scale=.35]{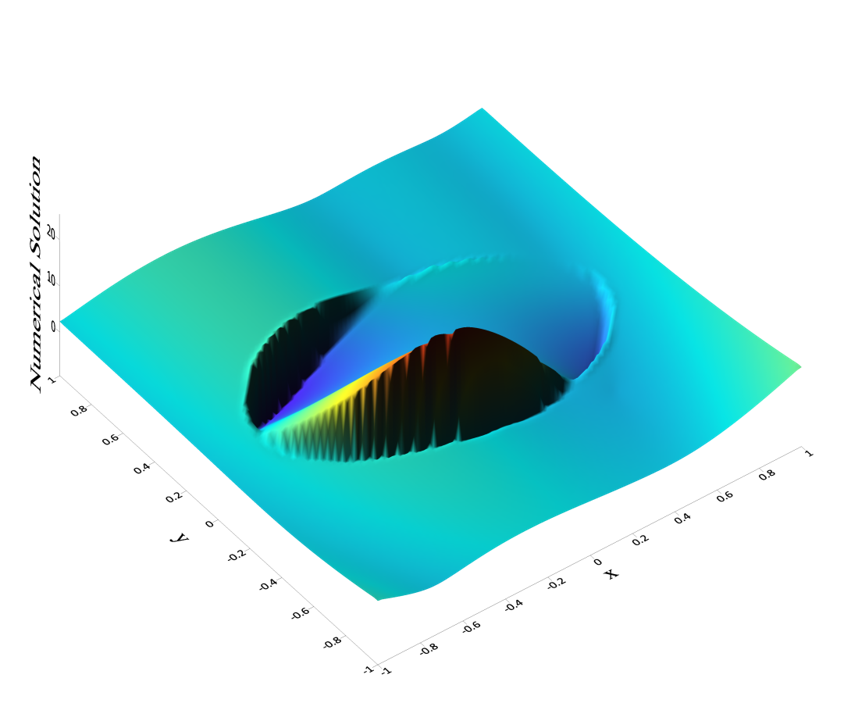}
	
	\caption{The numerical solution when $t=0.1$(left), $t=0.5$ (middle) and $t=1$ (right) under the ellipse interface for Example 4.}
\end{figure}

 \begin{figure}
 	\centering
 	\includegraphics[scale=.4]{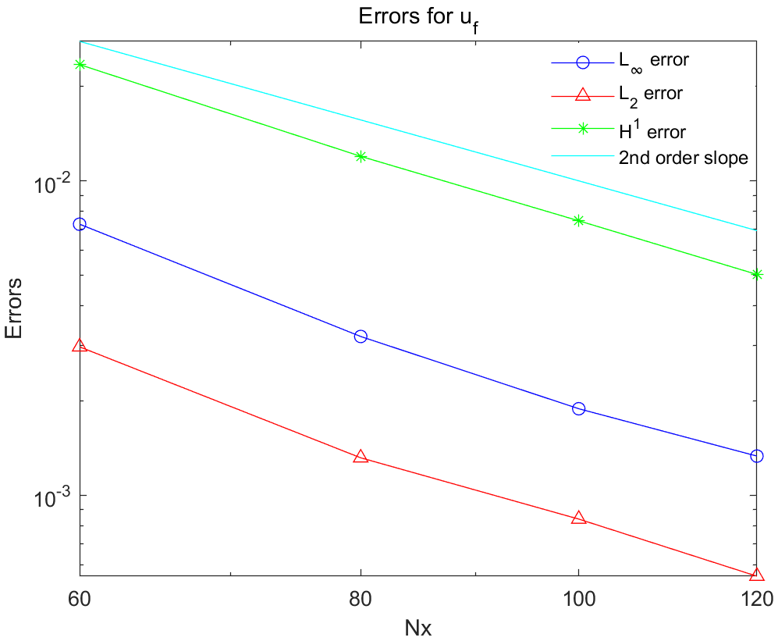}
 	\includegraphics[scale=.4]{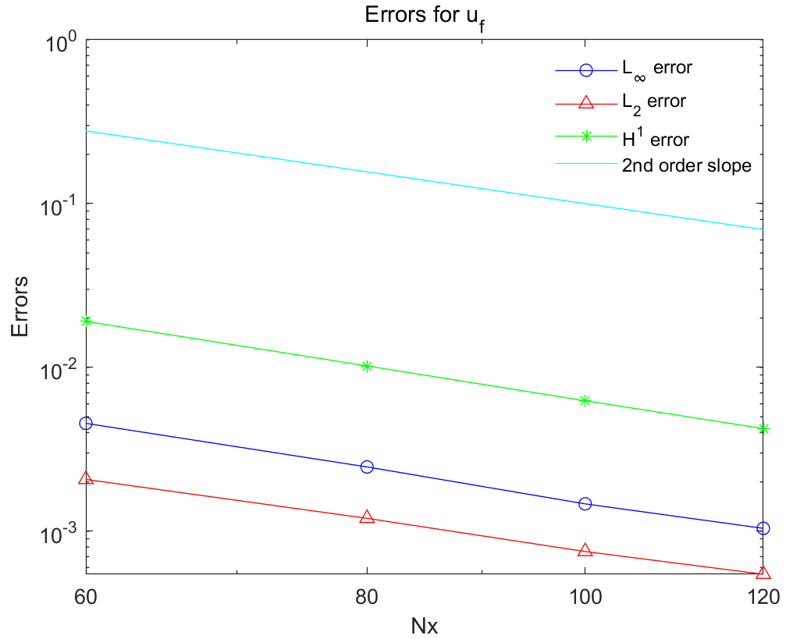}
 	\includegraphics[scale=.4]{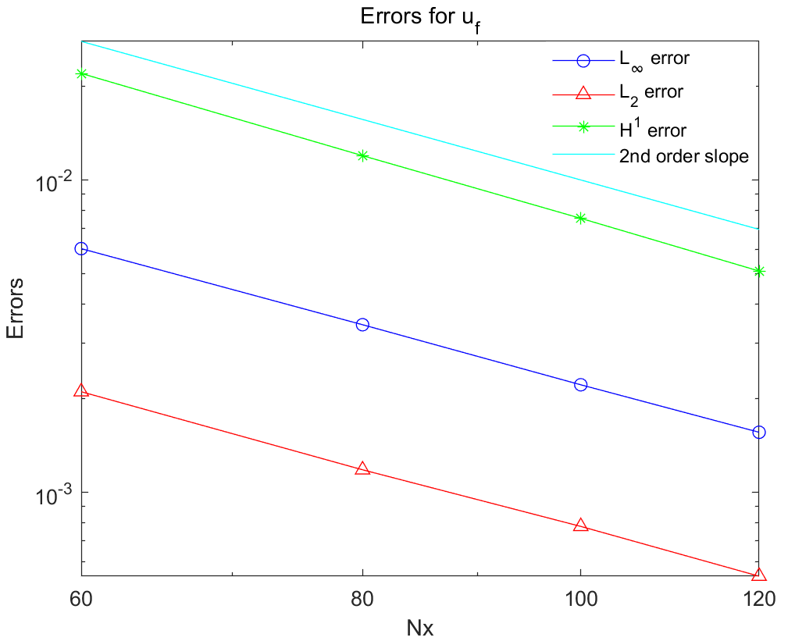}
 	
 	\caption{$L_{\infty}$, $L_2$ and $H^1$ errors of  $\textbf{u}_f$ when $t=0.1$(left), $t=0.5$ (middle) and $t=1$ (right) under the ellipse interface for Example 4.}
 \end{figure}

 \begin{figure}
 	\centering
 	\includegraphics[scale=.4]{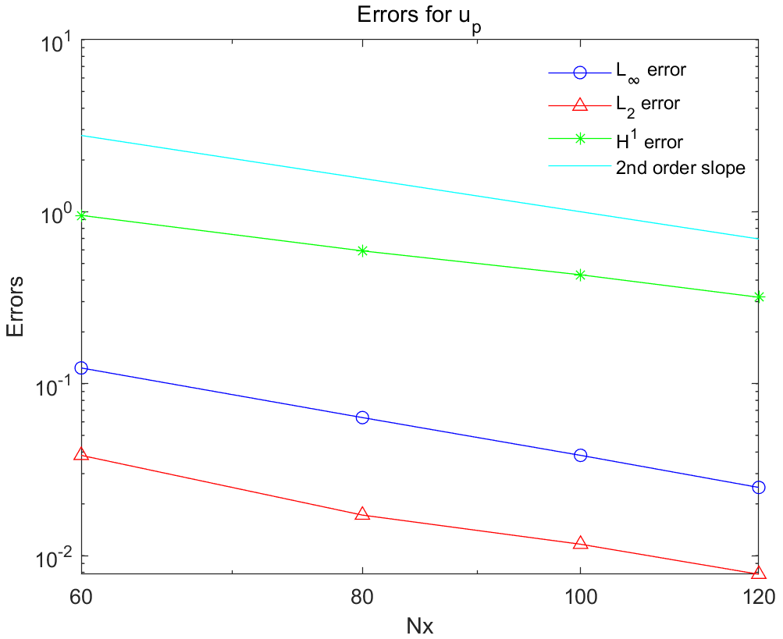}
 	\includegraphics[scale=.4]{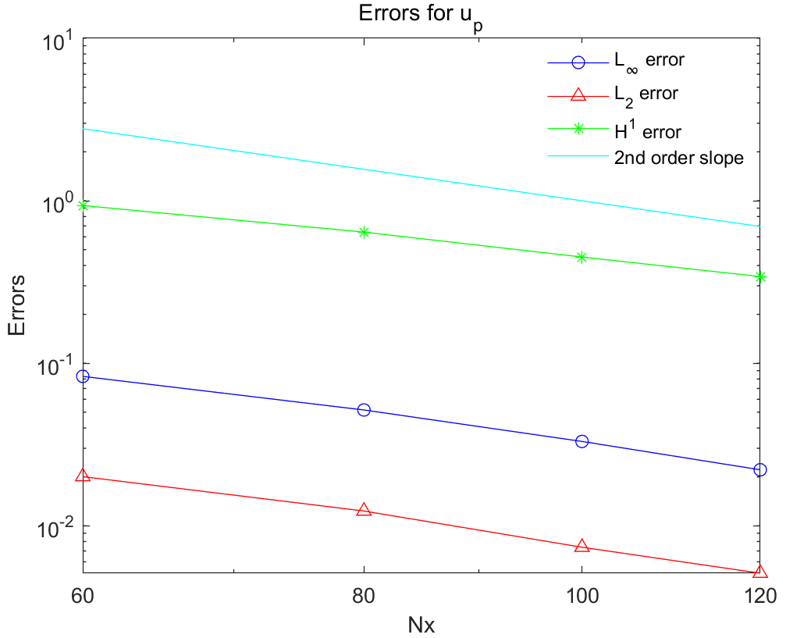}
 	\includegraphics[scale=.4]{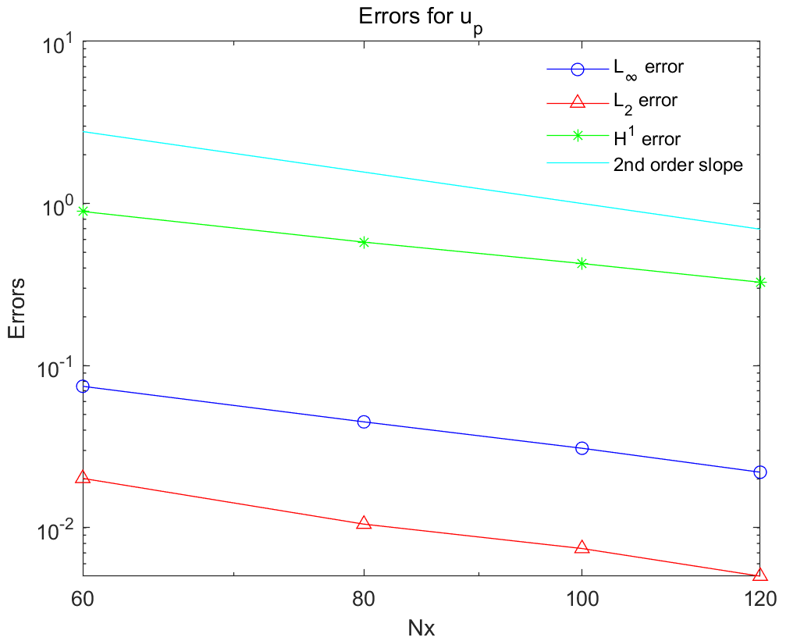}
 	
 	\caption{$L_{\infty}$, $L_2$ and $H^1$ errors of $\textbf{u}_p$ when $t=0.1$(left), $t=0.5$ (middle) and $t=1$ (right) under the ellipse interface for Example 4.}
 \end{figure}

 \begin{figure}
 	\centering
 	\includegraphics[scale=.4]{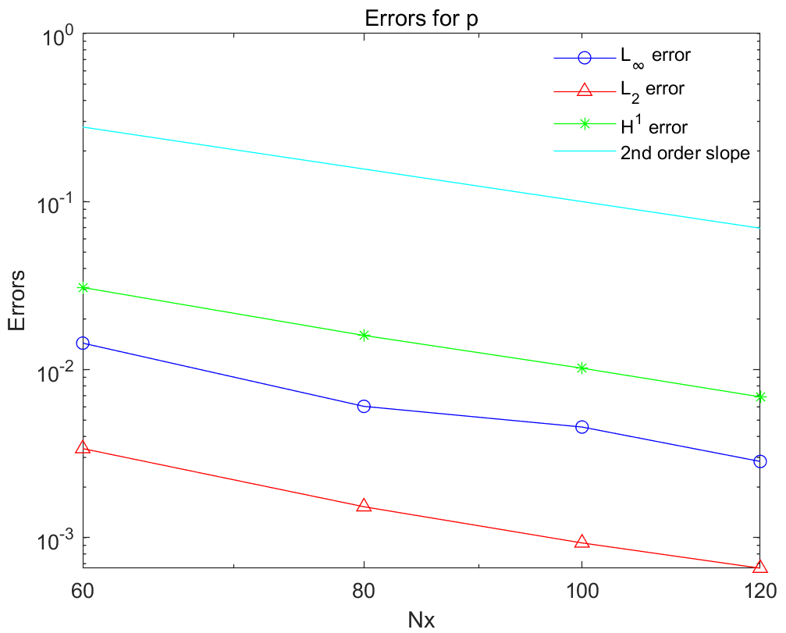}
 	\includegraphics[scale=.4]{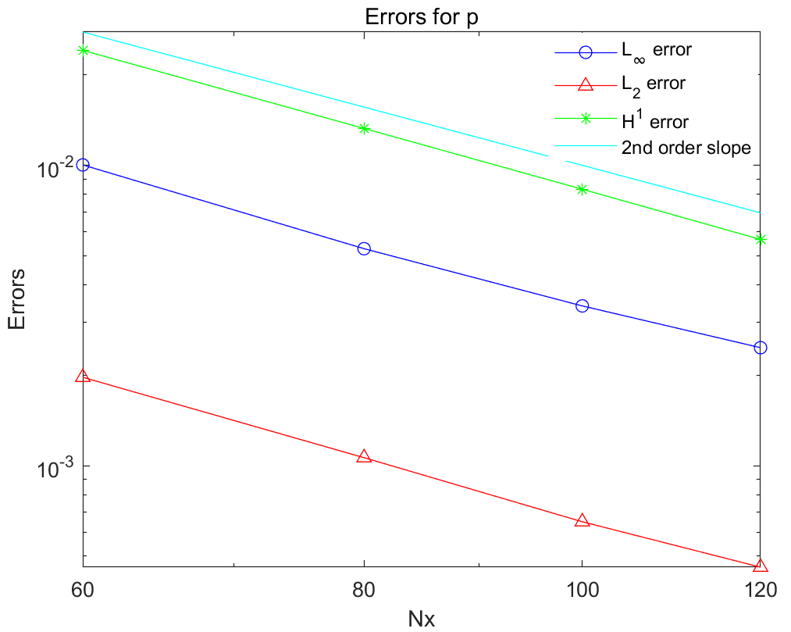}
 	\includegraphics[scale=.4]{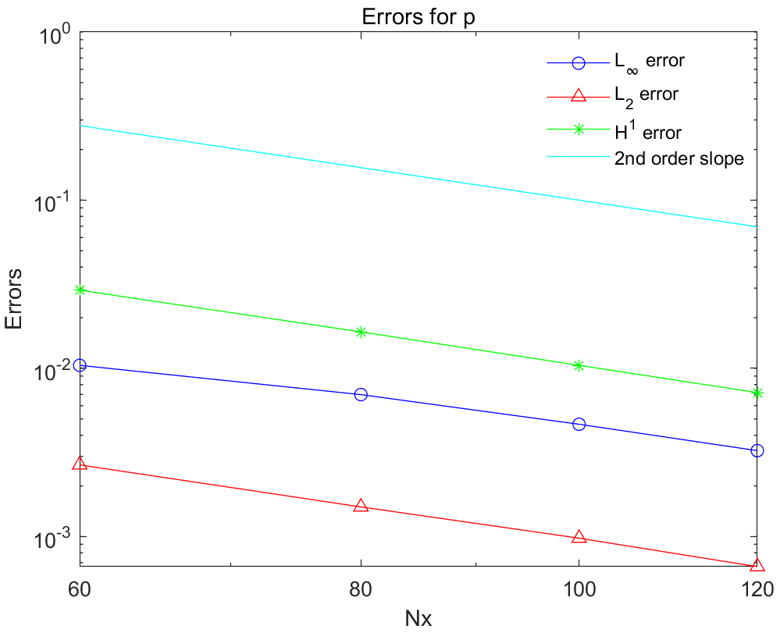}
 	
 	\caption{$L_{\infty}$, $L_2$ and $H^1$ errors of $p$ when $t=0.1$(left), $t=0.5$ (middle) and $t=1$ (right) under the ellipse interface for Example 4.}
 \end{figure}

	\begin{table*}	
	\scriptsize
	\caption{ The CPU times of the GFDM when $N_x=60, t=1$ for Example 3 and Example 4}
	\begin{tabular}{ccccccccccc}
		\hline
		\multirow{1}{*}{ }&\multicolumn{4}{c}{Example 3}& \multicolumn{2}{c}{Example 4}\\
		\hline		
			Interface&Circle interface& Two-petaled interface      &  Flower shaped interface 
		& Heart shaped interface      &  Petagon interface  
		&  Ellipse interface \\ 
		\hline
			CPU(s)&$0.48$& $1.39$&$1.07$ & $1.23$ & $1.25$ &  $1.06$ \\
			\hline
	\end{tabular}
\end{table*}
\begin{table*}	
	\scriptsize
	\caption{ $L_{\infty}$, $L_2$ and $H^1$ errors of the GFDM with different $m$  when $N_x=60, t=1$ and the petagon interface is used for Example 4}
	\begin{tabular}{cccccccccccc}
		\hline
		\multirow{1}{*}{$m$} & \multicolumn{3}{c}{$u_f$} & \multicolumn{3}{c}{$u_p$}& \multicolumn{2}{c}{$p$}&
		\multicolumn{2}{c}{$\phi$} \\
		\hline
		
		& $L_{\infty}$      &  $L_2 $   &   $H^1$
		
		& $L_{\infty}$      &  $L_2$   &   $H^{1,relative}$
		
		&  $L_2$   &   $H^1$
		&  $L_2$   &   $H^1$\\
		\hline
		$22$ & $9.95\times10^{-3}$ & $3.46\times10^{-3}$ & $2.12\times10^{-2}$ &  $7.92\times10^{-2}$ & $2.64\times10^{-2}$ & $3.28\times10^{-2}$&
		$4.39\times10^{-3}$ & $2.78\times10^{-2}$ & $2.00\times10^{-2}$&$2.64\times10^{-2}$ \\

		$26$ & $1.05\times10^{-2}$ & $3.64\times10^{-3}$ & $2.25\times10^{-2}$ &  $8.29\times10^{-2}$ & $2.78\times10^{-2}$ & $3.46\times10^{-2}$&
		$4.57\times10^{-3}$ & $2.97\times10^{-2}$ & $2.07\times10^{-2}$ &$2.78\times10^{-2}$\\
		
		$30$ & $1.12\times10^{-2}$ & $3.84\times10^{-3}$ & $2.45\times10^{-2}$ &  $8.75\times10^{-2}$ & $3.16\times10^{-2}$ & $3.76\times10^{-2}$&
		$4.93\times10^{-3}$ & $3.31\times10^{-2}$ & $2.16\times10^{-2}$&$3.15\times10^{-2}$ \\
		
		$36$ & $6.03\times10^{-3}$ & $2.10\times10^{-3}$ & $2.19\times10^{-2}$ &  $7.45\times10^{-2}$ & $2.02\times10^{-2}$ & $3.90\times10^{-2}$&
		$2.67\times10^{-3}$ & $2.92\times10^{-2}$ & $9.29\times10^{-3}$&$1.95\times10^{-2}$ \\
		
		$40$ & $6.39\times10^{-3}$ & $2.28\times10^{-3}$ & $2.53\times10^{-2}$ &  $8.69\times10^{-2}$ & $2.11\times10^{-2}$ & $4.08\times10^{-2}$&
		$2.69\times10^{-3}$ & $3.24\times10^{-2}$ & $8.96\times10^{-3}$&$2.04\times10^{-2}$ \\
		\hline
	\end{tabular}
\end{table*}

Fig. 23 and Fig. 28 present the point collocations and moving directions of the proposed problems with the petagon interface and the ellipse interface. we can see the moving directions and the interface location at the specific time. The numerical solutions at each time are shown in Fig.24 and Fig.29. We can see the distributions of the numerical solutions. Figs.25-27 present the numerical results of the GFDM for the Stokes-Darcy coupled problems with petagon interface when $t=0.1$ (left), $t=0.5$ (middle) and $t=1$ (right). Note that the $L_{\infty}$,$L_2$ and $H^1$ errors are all accurate, stable and keep 2nd order convergence. It means that the interface location has little influence on the numerical results of the proposed GFDM for the proposed problem with petagon interface. Figs.30-32 present the numerical results of the GFDM for the Stokes-Darcy coupled problems with ellipse interface when $t=0.1$ (left), $t=0.5$ (middle) and $t=1$ (right). We can see the same conclusion. It means that the interface location has little influence on the numerical results of the proposed GFDM for the proposed problem with ellipse interface.
The CPU times of the GFDM when $N_x=60,t=1$ for Example 3 and Example 4 is presented in Table. 12. From this table we can see the high efficiency of the proposed GFDM for the Stokes-Darcy coupled problem with the proposed all types of closed interfaces. The interface location may need different 'm', therefore, we also test the stability about the 'm' for the proposed problem with petagon interface in Table.13. Note that our numerical errors are all keep at almost the same accuracy and our results are all stable in the range from 22 to 40. The 'm' can be tested for the Stokes-Darcy peoblem with ellipse interface and we can get almost the same conclusion. In this example, we adopt $m=36.$
\section{Conclusion}
   In this paper, a meshless generalized finite difference method is proposed to solve the Stokes-Darcy coupled problem with BJS interface conditions. We adopt some high order GFDMs to show the high order accuracy and the convergence order of the GFDMs for Stokes-Darcy coupled problem. Some Stokes-Darcy coupled problems with closed interface which has more complex geometric shape are given to show the accuracy and stability of the GFDM. The interface location has been changed to show that the little influence of the interface location for the Stokes-Darcy coupled problem. Numerical results show that the better accuracy and stability can be obtained by GFDM compared with other high efficient numerical methods.
\section*{Acknowledgments}
 
 This work is partially supported by Science and Technology Commission of Shanghai Municipality (Grant Nos.  22JC1400900, 22DZ2229014).\\

\section*{Appendix A}
In 2D case, minimizing the residual function of the second-order partial derivatives in the second-order generalized finite difference method, a linear algebraic equation system can be obtained,
\begin{equation}
	AD_{2,u}=b_2,
\end{equation}
let 
\begin{equation}
	\begin{aligned}
		&\left.{A=}\right.\\
		&\left.{\left (
			\begin{array}{lllll}
				{\sum_{k=1}^m{h_k^2\omega_k^2}} &{\sum_{k=1}^m{h_kl_k\omega_k^2}}&{\sum_{k=1}^m{\frac{h_k^3}{2}\omega_k^2}} &{\sum_{k=1}^m{\frac{l_k^2h_k}{2}\omega_k^2}}&{\sum_{k=1}^m{h_k^2l_k\omega_k^2}} \\
				{\sum_{k=1}^m{h_kl_k\omega_k^2}} &{\sum_{k=1}^m{l_k^2\omega_k^2}}&{\sum_{k=1}^m{\frac{h_k^2l_k}{2}\omega_k^2}} &{\sum_{k=1}^m{\frac{l_k^3}{2}\omega_k^2}}&{\sum_{k=1}^m{\frac{l_k^2h_k}{2}\omega_k^2}} \\
				{\sum_{k=1}^m{\frac{h_k^3}{2}\omega_k^2}} &{\sum_{k=1}^m{\frac{h_k^2l_k}{2}\omega_k^2}}&{\sum_{k=1}^m{\frac{h_k^4}{4}\omega_k^2}} &{\sum_{k=1}^m{\frac{l_k^2h_k^2}{4}\omega_k^2}}&{\sum_{k=1}^m{\frac{h_k^3l_k}{2}\omega_k^2}} \\
				{\sum_{k=1}^m{\frac{l_k^2h_k}{2}\omega_k^2}} &{\sum_{k=1}^m{\frac{l_k^3}{2}\omega_k^2}}&{\sum_{k=1}^m{\frac{l_k^2h_k^2}{4}\omega_k^2}} &{\sum_{k=1}^m{\frac{l_k^4}{4}\omega_k^2}}&{\sum_{k=1}^m{\frac{l_k^3h_k}{2}\omega_k^2}} \\
				{\sum_{k=1}^m{h_k^2l_k\omega_k^2}} &{\sum_{k=1}^m{\frac{l_k^2h_k}{2}\omega_k^2}}&{\sum_{k=1}^m{\frac{h_k^3l_k}{2}\omega_k^2}} &{\sum_{k=1}^m{\frac{l_k^3h_k}{2}\omega_k^2}}&{\sum_{k=1}^m{l_k^2h_k^2\omega_k^2}} \\
			\end{array}\right)}\right.\\
		&\left.{=\left (
			\begin{array}{lllll}
				{a_{11}} &{a_{12}}&{a_{13}} &{a_{14}}&{a_{15}} \\
				{a_{21}}&{a_{22}} &{a_{23}}&{a_{24}}&{a_{25}} \\
				{a_{31}}&{a_{32}} &{a_{33}}&{a_{34}}&{a_{35}} \\
				{a_{41}}&{a_{42}} &{a_{43}}&{a_{44}} &{a_{45}}\\
				{a_{51}}&{a_{52}} &{a_{53}}&{a_{54}} &{a_{55}}\\
			\end{array}\right)}\right.
	\end{aligned}
\end{equation}
then
\begin{equation}
	A\left (
	\begin{array}{l}
		{\frac{\partial u_0}{\partial x}}\\
		{\frac{\partial u_0}{\partial y}}\\
		{\frac{\partial^2 u_0}{\partial x^2}}\\
		{\frac{\partial^2 u_0}{\partial y^2}}\\
		{\frac{\partial^2 u_0}{\partial x \partial y}}\\
	\end{array}\right)
	=\left (
	\begin{array}{l}
		{b_1}\\
		{b_2}\\
		{b_3}\\
		{b_4}\\
		{b_5}\\
	\end{array}\right)
	=\left (
	\begin{array}{l}
		{\sum_{k=1}^m{(-u_0+u_k)h_k\omega_k^2}}\\
		{\sum_{k=1}^m{(-u_0+u_k)l_k\omega_k^2}}\\
		{\sum_{k=1}^m{(-u_0+u_k)\frac{h_k^2}{2}\omega_k^2}}\\
		{\sum_{k=1}^m{(-u_0+u_k)\frac{l_k^2}{2}\omega_k^2}}\\
		{\sum_{k=1}^m{(-u_0+u_k)h_kl_k\omega_k^2}}\\
	\end{array}\right),
\end{equation}
Due to the matrix $A$ is symmetrical, it is possible to use the Cholesky method to solve the systems. The upper and lower triangular matries\cite{ref89}:
\begin{align}
	A=\left (
	\begin{array}{lllll}
		{a_{11}} &{a_{12}}&{a_{13}} &{a_{14}}&{a_{15}} \\
		{a_{21}}&{a_{22}} &{a_{23}}&{a_{24}}&{a_{25}} \\
		{a_{31}}&{a_{32}} &{a_{33}}&{a_{34}}&{a_{35}} \\
		{a_{41}}&{a_{42}} &{a_{43}}&{a_{44}} &{a_{45}}\\
		{a_{51}}&{a_{52}} &{a_{53}}&{a_{54}} &{a_{55}}\\
	\end{array}\right)=LL^T
	=\left (
	\begin{array}{lllll}
		{l_{11}} &{0}&{0} &{0}&{0} \\
		{l_{21}}&{l_{22}} &{0}&{0}&{0} \\
		{l_{31}}&{l_{32}} &{l_{33}}&{0}&{0} \\
		{l_{41}}&{l_{42}} &{l_{43}}&{l_{44}} &{0}\\
		{l_{51}}&{l_{52}} &{l_{53}}&{l_{54}} &{l_{55}}\\
	\end{array}\right)\left (
	\begin{array}{lllll}
		{l_{11}} &{l_{21}}&{l_{31}} &{l_{41}}&{l_{51}} \\
		{0}&{l_{22}} &{l_{32}}&{l_{42}}&{l_{52}} \\
		{0}&{0} &{l_{33}}&{l_{43}}&{l_{53}} \\
		{0}&{0} &{0}&{l_{44}} &{l_{45}}\\
		{0}&{0} &{0}&{0} &{l_{55}}\\
	\end{array}\right).
\end{align}
here
\begin{align}
	l_{i1}=\frac{a_{1i}}{\sqrt{a_{11}}}, i=1,2,3,4,5,\\	l_{i2}=\frac{a_{2i}-l_{21}l_{i1}}{\sqrt{l_{22}}}, i=2,3,4,5,\\
	l_{i3}=\frac{a_{3i}-l_{31}l_{i1}-l_{32}l_{i2}}{\sqrt{l_{33}}}, i=3,4,5,\\
	l_{i4}=\frac{a_{4i}-l_{41}l_{i1}-l_{42}l_{i2}-l_{43}l_{i3}}{\sqrt{l_{33}}}, i=4,5,\\
	l_{55}=\sqrt{a_{55}-(l_{51}^2+l_{52}^2+l_{53}^2+l_{54}^2)}
\end{align}
\begin{equation}
	\begin{aligned}
		&\left.{\left (
			\begin{array}{lllll}
				{l_{11}} &{0}&{0} &{0}&{0} \\
				{l_{21}}&{l_{22}} &{0}&{0}&{0} \\
				{l_{31}}&{l_{32}} &{l_{33}}&{0}&{0} \\
				{l_{41}}&{l_{42}} &{l_{43}}&{l_{44}} &{0}\\
				{l_{51}}&{l_{52}} &{l_{53}}&{l_{54}} &{l_{55}}\\
			\end{array}\right)\left (
			\begin{array}{lllll}
				{l_{11}} &{l_{21}}&{l_{31}} &{l_{41}}&{l_{51}} \\
				{0}&{l_{22}} &{l_{32}}&{l_{42}}&{l_{52}} \\
				{0}&{0} &{l_{33}}&{l_{43}}&{l_{53}} \\
				{0}&{0} &{0}&{l_{44}} &{l_{45}}\\
				{0}&{0} &{0}&{0} &{l_{55}}\\
			\end{array}\right)
			\left (
			\begin{array}{l}
				{\frac{\partial u_0}{\partial x}}\\
				{\frac{\partial u_0}{\partial y}}\\
				{\frac{\partial^2 u_0}{\partial x^2}}\\
				{\frac{\partial^2 u_0}{\partial y^2}}\\
				{\frac{\partial^2 u_0}{\partial x \partial y}}\\
			\end{array}\right)
			=\left (
			\begin{array}{l}
				{\sum_{k=1}^m{(-u_0+u_k)h_k\omega_k^2}}\\
				{\sum_{k=1}^m{(-u_0+u_k)l_k\omega_k^2}}\\
				{\sum_{k=1}^m{(-u_0+u_k)\frac{h_k^2}{2}\omega_k^2}}\\
				{\sum_{k=1}^m{(-u_0+u_k)\frac{l_k^2}{2}\omega_k^2}}\\
				{\sum_{k=1}^m{(-u_0+u_k)h_kl_k\omega_k^2}}\\
			\end{array}\right)
		}\right.
	\end{aligned}
\end{equation}
\begin{equation}
	\begin{aligned}
		&\left.{
			\left (
			\begin{array}{l}
				{\frac{\partial u_0}{\partial x}}\\
				{\frac{\partial u_0}{\partial y}}\\
				{\frac{\partial^2 u_0}{\partial x^2}}\\
				{\frac{\partial^2 u_0}{\partial y^2}}\\
				{\frac{\partial^2 u_0}{\partial x \partial y}}\\
			\end{array}\right)=
			\left (
			\begin{array}{lllll}
				{l_{11}} &{0}&{0} &{0}&{0} \\
				{l_{21}}&{l_{22}} &{0}&{0}&{0} \\
				{l_{31}}&{l_{32}} &{l_{33}}&{0}&{0} \\
				{l_{41}}&{l_{42}} &{l_{43}}&{l_{44}} &{0}\\
				{l_{51}}&{l_{52}} &{l_{53}}&{l_{54}} &{l_{55}}\\
			\end{array}\right)^{-1}
			\left (
			\begin{array}{lllll}
				{l_{11}} &{l_{21}}&{l_{31}} &{l_{41}}&{l_{51}} \\
				{0}&{l_{22}} &{l_{32}}&{l_{42}}&{l_{52}} \\
				{0}&{0} &{l_{33}}&{l_{43}}&{l_{53}} \\
				{0}&{0} &{0}&{l_{44}} &{l_{45}}\\
				{0}&{0} &{0}&{0} &{l_{55}}\\
			\end{array}\right)^{-1}
			\left (
			\begin{array}{l}
				{\sum_{k=1}^m{(-u_0+u_k)h_k\omega_k^2}}\\
				{\sum_{k=1}^m{(-u_0+u_k)l_k\omega_k^2}}\\
				{\sum_{k=1}^m{(-u_0+u_k)\frac{h_k^2}{2}\omega_k^2}}\\
				{\sum_{k=1}^m{(-u_0+u_k)\frac{l_k^2}{2}\omega_k^2}}\\
				{\sum_{k=1}^m{(-u_0+u_k)h_kl_k\omega_k^2}}\\
			\end{array}\right)
		}\right.\\
	\end{aligned}
\end{equation}
\begin{equation}
	\begin{aligned}
		&\left.{\left (
			\begin{array}{lllll}
				{l_{11}} &{0}&{0} &{0}&{0} \\
				{l_{21}}&{l_{22}} &{0}&{0}&{0} \\
				{l_{31}}&{l_{32}} &{l_{33}}&{0}&{0} \\
				{l_{41}}&{l_{42}} &{l_{43}}&{l_{44}} &{0}\\
				{l_{51}}&{l_{52}} &{l_{53}}&{l_{54}} &{l_{55}}\\
			\end{array}\right)^{-1}
			\left (
			\begin{array}{lllll}
				{l_{11}} &{l_{21}}&{l_{31}} &{l_{41}}&{l_{51}} \\
				{0}&{l_{22}} &{l_{32}}&{l_{42}}&{l_{52}} \\
				{0}&{0} &{l_{33}}&{l_{43}}&{l_{53}} \\
				{0}&{0} &{0}&{l_{44}} &{l_{45}}\\
				{0}&{0} &{0}&{0} &{l_{55}}\\
			\end{array}\right)^{-1}
			\stackrel{\triangle}{=}
			\left (
			\begin{array}{lllll}
				{c_{11}} &{c_{12}}&{c_{13}} &{c_{14}}&{c_{15}} \\
				{c_{21}}&{c_{22}} &{c_{23}}&{c_{24}}&{c_{25}} \\
				{c_{31}}&{c_{32}} &{c_{33}}&{c_{34}}&{c_{35}} \\
				{c_{41}}&{c_{42}} &{c_{43}}&{c_{44}} &{c_{45}}\\
				{c_{51}}&{c_{52}} &{c_{53}}&{c_{54}} &{c_{55}}\\
			\end{array}\right)
		}\right.
	\end{aligned}
\end{equation}
then
\begin{equation}
	\begin{aligned}
		&\left.{
			\left (
			\begin{array}{l}
				{\frac{\partial u_0}{\partial x}}\\
				{\frac{\partial u_0}{\partial y}}\\
				{\frac{\partial^2 u_0}{\partial x^2}}\\
				{\frac{\partial^2 u_0}{\partial y^2}}\\
				{\frac{\partial^2 u_0}{\partial x \partial y}}\\
			\end{array}\right)=
			\left (
			\begin{array}{lllll}
				{c_{11}} &{c_{12}}&{c_{13}} &{c_{14}}&{c_{15}} \\
				{c_{21}}&{c_{22}} &{c_{23}}&{c_{24}}&{c_{25}} \\
				{c_{31}}&{c_{32}} &{c_{33}}&{c_{34}}&{c_{35}} \\
				{c_{41}}&{c_{42}} &{c_{43}}&{c_{44}} &{c_{45}}\\
				{c_{51}}&{c_{52}} &{c_{53}}&{c_{54}} &{c_{55}}\\
			\end{array}\right)
			\left (
			\begin{array}{l}
				{\sum_{k=1}^m{(-u_0+u_k)h_k\omega_k^2}}\\
				{\sum_{k=1}^m{(-u_0+u_k)l_k\omega_k^2}}\\
				{\sum_{k=1}^m{(-u_0+u_k)\frac{h_k^2}{2}\omega_k^2}}\\
				{\sum_{k=1}^m{(-u_0+u_k)\frac{l_k^2}{2}\omega_k^2}}\\
				{\sum_{k=1}^m{(-u_0+u_k)h_kl_k\omega_k^2}}\\
			\end{array}\right)
		}\right.\\
	\end{aligned}
\end{equation}
We consider $f_2(u)=\beta_1\frac{\partial u}{\partial y}$ as a example, $\beta_1$ is a constant, and then
\begin{equation}
	\begin{aligned}
		f_2(u)=&\left.{
			-u_0(\beta_1(\sum_{k=1}^m{[c_{21}h_k\omega_k^2+c_{22}l_k\omega_k^2+c_{23}\frac{h_k^2}{2}\omega_k^2+c_{24}\frac{l_k^2}{2}\omega_k^2+c_{25}h_kl_k\omega_k^2]}))}\right.\\
		&\left.{+\sum_{k=1}^m{[\beta_1(c_{21}h_k\omega_k^2+c_{22}l_k\omega_k^2+c_{23}\frac{h_k^2}{2}\omega_k^2+c_{24}\frac{l_k^2}{2}\omega_k^2+c_{25}h_kl_k\omega_k^2)]}u_k.
		}\right.\\
	\end{aligned}
\end{equation}
According to the explicit form of the generalized finite difference method\textcolor{blue}{\eqref{Eq:100}}, the coefficients $m_0$ and $m_k$ are as follows,
\begin{align}
	m_0=\beta_1(\sum_{k=1}^m{[c_{21}h_k\omega_k^2+c_{22}l_k\omega_k^2+c_{23}\frac{h_k^2}{2}\omega_k^2+c_{24}\frac{l_k^2}{2}\omega_k^2+c_{25}h_kl_k\omega_k^2]}),\\
	m_k=\beta_1(c_{21}h_k\omega_k^2+c_{22}l_k\omega_k^2+c_{23}\frac{h_k^2}{2}\omega_k^2+c_{24}\frac{l_k^2}{2}\omega_k^2+c_{25}h_kl_k\omega_k^2),
\end{align}
that carries out condition $m_0=\sum_{k=1}^mm_k$.

\printcredits


\begin{thebibliography}{99}
\bibitem{ref24} {Y. Cao, M. Gunzburger, X. Hu, F. Hua, X. Wang, W. Zhao}, Finite element approximations for Stokes-Darcy flow with Beavers-Joseph interface conditions, SIAM.J. Numer. Anal., 47, 4239-4256 (2010).

\bibitem{ref25} {H. Rui, R. Zhang}, A unified stabilized mixed finite element method for coupling Stokes and Darcy flows, Comput. Methods Appl. Mech. Engrg., 198, 33-36(2009).

\bibitem{ref26}{J.P. Yu, Y. Z. Sun, F. Shi, H. B. Zheng}, Nitsche’s type stabilized finite element method for the fully mixed Stokes-Darcy problem with Beavers-Joseph conditions, Appl. Math. Letters., 110, 106588 (2020).

\bibitem{ref27}{MA. A. Mahbub, L. Shan, H. B. Zheng}, Uncoupling evolutionary groundwater-surface water flows: stabilized mixed methods in both porous media and fluid regions, Numer. Algor., 92(3), 1837-1874 (2023).

\bibitem{ref28}{M. Mu, J. Xu}, A two-grid method of a mixed Stokes-Darcy model for coupling fluid flow with porous media flow, SIAM J. Numer. Anal., 45, 1801-1813 (2007).

\bibitem{ref29} {M. Cai, M. Mu, J. Xu}, Numerical solution to a mixed Navier-Stokes/Darcy model by the two-grid approach, SIAM J. Numer. Anal., 47, 3325-3338 (2009).

\bibitem{ref30} {L. Y. Zuo, Y. Hou}, A decoupling two-grid algorithm for the mixed Stokes-Darcy model with the Beavers-Joseph interface condition, Numer. Methods PDEs., 30, 1066-1082 (2014).

\bibitem{ref31} {Y. Hou}, Optimal error estimates of a decoupled scheme based on two-grid finite element for mixed Stokes-Darcy model, App. Math. Letters., 57, 90-96 (2016).

\bibitem{ref34}{W. J. Layton, F. Schieweck, I. Yotov}, Coupling fluid flow with porous media flow, SIAM J. Numer. Anal., 40, 2195-2218 (2003).

\bibitem{ref38} {M. Gunzburger, X. M. He, B. Li}, On Stokes-Ritz projection and multistep backward differentiation schemes in decoupling the Stokes-Darcy model, SIAM J. Numer. Anal., 56, 397-427 (2018).

\bibitem{ref39} {L. Shan, H. B. Zheng}, Partitioned time stepping method for fully evolutionary Stokes-Darcy flow with Beavers-Joseph interface conditions, SIAM J. Numer. Anal.,
51, 813-839 (2013).

\bibitem{ref40} {L. Shan, H. B. Zheng, W. J. Layton}, A decoupling method with different subdomain time steps for the nonstationary Stokes-Darcy model, Numer. Methods PDEs., 29, 549-583 (2013).

\bibitem{ref41} {W. Layton, H. Tran,C. Trenchea}, Analysis of long time stability and error of two partitioned methods for uncoupling evolutionary groudwater-surface water flows, SIAM J. Numer. Anal., 51, 248-272 (2013).

\bibitem{ref44} {M. Discacciati, A. Quarteroni, A. Valli}, Robin-Robin domain decomposition methods for the Stokes-Darcy coupling, SIAM J. Numer. Anal., 45, 1246-1268 (2007).

\bibitem{ref45} {W. Chen, M. Gunzburger, F. Hua,X. Wang}, A parallel Robin-Robin domain decomposition method for the Stokes-Darcy system, SIAM. J. Numer. Anal., 49, 1064-1084 (2011).

\bibitem{ref47}{Y. Boubendir, S. Tlupova S}, Domain decomposition methods for solving Stokes-Darcy problems with boundary integrals, SIAM J. Sci. Comput., 35, B82-B106 (2013).

\bibitem{ref49} {X. M. He, J. Li, Y. P. Lin, J. Ming}, A domain decomposition method for the steadystate Navier-Stokes-Darcy model with the Beavers-Joseph interface condition, SIAM J. Sci. Comput., 37, S264-S290 (2015).

\bibitem{ref50} {D. Vassilev, C. Wang, I. Yotov}, Domain decomposition for coupled Stokes and Darcy flows, Comput. Methods Appl. Mech. Engrg., 268, 264-283 (2014).

\bibitem{ref51} {B. Jiang}, A parallel domain decomposition method for coupling of surface and groundwater flows, Comput. Methods Appl. Mech. Engrg., 198, 947-957 (2009).

\bibitem{ref52} {Y. Z. Sun, W. W. Sun, H. B. Zheng}, Domain decomposition method for the fully-mixed Stokes-Darcy coupled problem, Comput. Methods Appl. Mech. Engrg., 374, 113578 (2021).

\bibitem{ref53} {Y. Z. Sun, F. Shi, H. B. Zheng, H. Li, F. Wang }, Two-grid domain decomposition methods for the coupled Stokes-Darcy system, Comput. Methods Appl. Mech. Engrg., 385, 114041 (2021).

\bibitem{ref54} {W. B. Chen, M. Gunzburger, F. Hua, X. M. Wang}, Parallel Robin-Robin domain decomposition method for the Stokes-Darcy system, SIAM. J. Numer. Anal., 49(3), 1064-1084 (2011).

\bibitem{ref55} {R. Li, J. Li, Z. X. Chen, Y. L. Gao}, A stabilized finite element method based on two local Gauss integrations for a coupled Stokes-Darcy problem, J Comput. Appl. Math., 292, 92-104 (2016).

\bibitem{ref56} {Z. L. Li}, An augmented Cartesian grid method for Stokes-Darcy fluid-structure interactions, Int. J. Numer. Meth. Engng., 106, 556-575 (2016).

\bibitem{ref57} {M.Safarpoor, A.Shirzadi}, Numerical investigation based on radial basis function-finite-difference (RBF-FD) method for solving  the Stokes-Darcy equations, Eng. Comput-Germany. 37, 909-920(2021).


\bibitem{ref58} {R. L. Pu, X. L. Feng}, Physics-Informed Neural Networks for solving coupled Stokes-Darcy equation, Entropy 24(8), 1106(2022). 

\bibitem{ref59} {J. Yue, J. Li}, Efficient coupled deep neural networks for the time-dependent coupled Stokes-Darcy problems, Appl. Math. Comput.,437, 127514(2023).

\bibitem{ref60} {J. Orkisz}, Meshless finite difference method I Basic approach. Computational mechanics IACM. Idelshon Oňate, Duorkin, editors. CINME; 1998.

\bibitem{ref61} {J. J. Benito , F. Urena and L. Gavete}, Influence of several factors in the generalized finite difference method, Appl. Math. Model., 25, 1039-1053 (2011).

\bibitem{ref62} {J. J. Benito , F. Urena and L. Gavete}, An h-adaptive method in generalized finite differences, Comput. Methods App. Mech. Eng., 192(5-6), 735-759(2003).

\bibitem{ref63} {L. Gavete, F. Urena, J. J. Benito, A. Garcia, M. Urena, E. Salete}, Solving second order non-linear elliptic partial differential equations using generalized finite difference method, J.Comput Appl.Math., 318, 378-387(2017).

\bibitem{ref64} {E. Salete, J. J. Benito, F. Urena, L. Gavete, M. Urena, A. Garcia}, Stability of perfectly matched layer regions in generalized finite difference method for wave problems, J.Comput Appl.Math., 312, 231-239(2017).

\bibitem{ref65} {C. M. Fan, P. W. Li, W. C. Yeih},  Generalized finite difference method for solving two-dimensional inverse Cauchy problems, Inverse. Probl. Sci. Eng., 23(5), 737-759(2015).

\bibitem{ref66} {Y. Gu, L. Wang, W. Chen, C. Z. Zhang, X. Q. He}, Application of the meshless generalized finite difference method to inverse heat source problems, Int. J. Heat. Mass. Transf., 108, 721-729(2017).

\bibitem{ref67} {O. Davydov and M. Safarpoor}, A meshless finite difference method for elliptic interface
problems based on pivoted QR decomposition, Appl. Numer. Math., 161, 489–509 (2021).

\bibitem{ref68} {Y. N. Xing, H. B. Zheng}, A high order generalized finite difference method for solving the anisotropic elliptic interface problem in static and moving systems, Comput. Math. Appl., 166C, 1-23(2024).

\bibitem{ref69} {H. Kraus, J. Kuhnert, A. Meister and P. Suchde}, A meshfree point collocation method for elliptic interface problems, Appl. Math. Model., 113, 241–261(2023).

\bibitem{ref70} {Y. N. Xing, L.N. Song, X. M. He, C. X. Qiu}, A generalized finite difference method for solving elliptic interface problems, Math. Comput. Simulat., 178, 109-124(2020).

\bibitem{ref71} {Y. N. Xing, L. N. Song and C. M. Fan}, A generalized finite difference method for solving elasticity interface problems, Eng. Anal. Bound. Elem., 128(1), 105-117(2021).

\bibitem{ref72} {Y. N. Xing, L.N. Song and P. W. Li}, A generalized finite difference method for solving biharmonic interface problems, Eng. Anal. Bound. Elem., 135, 132-144 (2022).

\bibitem{ref86} {L. N. Song, P. W. Li, Y. Gu, C. M. Fan}, Generalized finite difference method for solving stationary 2D and 3D Stokes equations with a mixed boundary condition, Comput. Math. Appl., 80, 1726-1743(2020).

\bibitem{ref87} {M. R. Shao, L. N. Song and P. W. Li}, A generalized finite difference method for solving stokes interface problems, Eng. Anal. Bound. Elem., 132(1), 50-64(2021).

\bibitem{Shirokoff} {D. Shirokoff}, I.: A Pressure Poisson Method for the Incompressible Navier-Stokes Equations: II. Long Time Behavior of the Klein-Gordon Equations. PhD thesis, Massachusetts Institute of Technology (2011). 

\bibitem{ref89}{J. J. Benito, F. Urena and L. Gavete}, Solving parabolic and hyperbolic equations by the generalized finite difference method, J. Comput. Appl. Math., 209, 208-233(2007).


		\end{thebibliography}



\end{document}